\documentclass[11pt,letterpaper]{amsart}

\usepackage{amsmath,amsthm,amsfonts,amssymb}
\usepackage{setspace}
\usepackage[pdftex,pdfborder={0 0 0}]{hyperref}
\usepackage{fullpage}
\usepackage{enumerate}
\usepackage{datetime}
\usepackage[margin=1in]{geometry}

\usepackage[english]{babel}

\date{\today}
\hypersetup{allcolors=black}
\usepackage{titlesec}
\usepackage{xcolor}
\DeclareMathOperator{\Tor}{Tor}
\newcommand{\m}{\mathfrak m}
\newcommand{\Pone}{\mathcal P_{\mathrm{i}}}
\newcommand{\st}{\simeq_{\mathrm{st}}}

\usepackage{mathrsfs}
\usepackage{hyperref}
\usepackage{comment}
\usepackage{tikz-cd}
\usepackage{quiver}
\makeatletter
\providecommand{\@secnumpunct}{.}
\DeclareMathOperator{\heightideal}{ht}
\DeclareMathOperator{\NF}{NF}
\makeatother
\usepackage{setspace}
\usepackage{fancyhdr}
\usepackage{wallpaper}
\newcommand{\seq}{\subseteq}
\newcommand{\x}{\xrightarrow{}}

\newcommand{\sW}{\mathfrak W}

\DeclareMathOperator{\type}{type}

\newcommand{\PiC}{\mathcal P_{i}}
\newcommand{\PiiC}{\mathcal P_{ii}}

\DeclareMathOperator{\hf}{\mathscr{H}\kern-3pt\mathit{om}}
\DeclareMathOperator{\h}{Hom}
\DeclareMathOperator{\Ext}{Ext}

\DeclareMathOperator{\dep}{depth}

\DeclareMathOperator{\Ass}{Ass}
\DeclareMathOperator{\Supp}{Supp}
\DeclareMathOperator{\Min}{Min}
\DeclareMathOperator{\Spec}{Spec}
\DeclareMathOperator{\Ann}{Ann}
\DeclareMathOperator{\img}{im}

\DeclareMathOperator{\gr}{grade}

\DeclareMathOperator{\depth}{depth}

\DeclareMathOperator{\Proj}{Proj}

\DeclareMathOperator{\coker}{coker}
\DeclareMathOperator{\rank}{rank}
\DeclareMathOperator{\height}{height}
\theoremstyle{plain}
\newtheorem{thm}{Theorem}[section]
\newtheorem{lem}[thm]{Lemma}
\newtheorem{prop}[thm]{Proposition}
\newtheorem{cor}[thm]{Corollary}
\newtheorem{conjecture}[thm]{Conjecture}
\newtheorem{question}[thm]{Question}

\theoremstyle{definition}
\newtheorem{defi}[thm]{Definition}
\newtheorem{ex}[thm]{Example}

\theoremstyle{remark}
\newtheorem{remark}[thm]{Remark}

\newenvironment{recalledquestion}[1]
{%
  \par\medskip
  \noindent
  \textbf{\hyperref[#1]{Question~\ref*{#1}.}}
  \ \itshape
}
{%
  \par\medskip
}

\newenvironment{recalledtheorem}[1]
{%
  \par\vspace{7pt}
  \noindent
  \textbf{\hyperref[#1]{Theorem~\ref*{#1}.}}
  \ \itshape
}
{%
  \par\vspace{7pt}
}

\makeatletter
\def\l@section{\@tocline{1}{0pt}{1pc}{3pc}{}}
\def\l@subsection{\@tocline{2}{0pt}{1pc}{5pc}{}}
\makeatother

\usepackage{titlesec}

\titleformat{\subsection}
  {\bfseries}
  {\thesubsection.}
  {0.5em}
  {}

\title{Gorenstein Linkage of Ideals and Doubly layered filtrations of MCM modules}
\author{Siva Somasundaram}
\address{Department of Mathematics,
Purdue University,
West Lafayette, IN 47907, USA}

\email{somasus@purdue.edu}

\date{}

\begin{document}

\begin{abstract} \noindent The glicci conjecture asks whether every Cohen-Macaulay ideal in a regular ring belongs to the Gorenstein linkage class of a complete intersection. We show that the corresponding statement is false for nonregular rings. We prove that the conjecture holds for zero-dimensional ideals and for ideals of dimensions one and two in polynomial rings. Geometrically, this implies that every zero-dimensional subscheme of \(\mathbb P^n_k\), as well as curves in \(\mathbb P^n_k\), is glicci. We also prove that every Cohen–Macaulay monomial ideal is glicci. Our approach involves reducing the problem to height-two ideals and studying their \(N\)-type resolutions and the filtrations of the associated modules. We also show that every orientable maximal Cohen–Macaulay module over a two-dimensional Gorenstein ring is stably isomorphic to a doubly layered module. Finally, we prove that every ideal belongs to the Gorenstein linkage class of a Golod ideal.
\end{abstract}

\begingroup \maketitle \footnotesize \tableofcontents \endgroup

\section{ Introduction}
Linkage theory, also called liaison theory, was formally introduced by Peskine and Szpiro \cite{peskine1974liaison} and has been an important tool in the study of unmixed ideals. Let $R$ be a commutative Noetherian ring, and let $I$ and $J$ be proper ideals. We say that $I$ and $J$ are complete intersection linked if there exists a complete intersection ideal $(\alpha_1,...,\alpha_g)$ such that $J=(\alpha_1,...,\alpha_g):I$ and $I=(\alpha_1,...,\alpha_g):J$; we write $I\sim J$. We say $I$ and $J$ are in the same linkage class if there exists a finite sequence of ideals $I=I_0,..,I_n=J$ such that $I_{i-1}\sim I_{i}$. We say an ideal $I$ is \emph{licci} (in the linkage class of a complete intersection) if $I$ and a complete intersection ideal are in the same linkage class. Basic properties of complete intersection linkage can be found in \cite{peskine1974liaison}, and further properties of licci ideals in \cite{huneke1987structure, ulrich1987theory}.

Schenzel \cite{schenzel1982notes} introduced Gorenstein linkage as a generalization of complete intersection linkage, building
on ideas already present in the work of Peskine and Szpiro. Let $R$ be a Gorenstein local ring, and let $I$ and $J$ be proper ideals of $R$. We say $I$ and $J$ are Gorenstein linked if there exists a Gorenstein ideal $K$ such that $I=K:J$, and $J=K:I$; we write $I\sim_G J$. We say $I$ and $J$ are in the same Gorenstein linkage class if there exists a finite sequence of ideals $I=I_0,..,I_n=J$ such that $I_{i-1}\sim_G I_{i}$.  An ideal $I$ is called \emph{glicci} (in the Gorenstein linkage class of a complete intersection) if $I$ and a complete intersection ideal are in the same Gorenstein linkage class. Basic properties of Gorenstein linkage and glicci ideals can be found in \cite{schenzel1982notes, kustin1984deformation,kleppe2001gorenstein,migliore2002liaison,migliore2013glicci}.

A graded ring is called \emph{graded-local} if it has a unique homogeneous maximal ideal. A positively graded ring
$R=\bigoplus_{n\geq0}R_n$ is called
standard graded over $R_0$ if
$R=R_0[R_1]$. For graded-local rings and homogeneous ideals, we use the corresponding homogeneous notions. Thus, two homogeneous ideals are homogeneously complete intersection linked (respectively, homogeneously Gorenstein linked) if the linking complete intersection $(\alpha_1,\cdots,\alpha_g)$ (respectively, the linking Gorenstein ideal $K$) is homogeneous. A homogeneous ideal is homogeneously licci (respectively, homogeneously glicci) if it can be linked to a homogeneous complete intersection by homogeneous complete intersection (respectively, homogeneous Gorenstein) links. The homogeneous setting allows Hilbert functions, Hilbert series, and multiplicities to be studied
along homogeneous links.

There are several applications for the theory of linkage. Historically the properties of linkage have been studied and used since the late 1930s. For some recent developments, see \cite{guerrieri2025grade, chmiel2025structure, chong2016application, christensen2019free, guerrieri2024ade, jelisiejew2026hilbert}. One useful perspective on licci
ideals is the following: if a property or a conjecture is known for complete intersections, one can ask whether the same property holds for all licci ideals. This perspective has led to several successful applications; see \cite{herzog1984differentials, ulrich1987vanishing, planas1996vanishing, johnson2022splitting, huneke1989powers, eisenbud1995evolutions}.

An ideal $I\seq R$ is called Cohen-Macaulay if $R/I$ is a Cohen-Macaulay ring. By \cite[Corollary~3.3]{schenzel1982notes}, if $R$ is a local Gorenstein ring and  $I$ and $J$ are
ideals in the same Gorenstein linkage class, then $I$ is Cohen-Macaulay if and only if $J$ is Cohen-Macaulay. The same argument applies in the homogeneous setting. Consequently, every (homogeneously) glicci ideal is Cohen-Macaulay.

\begin{conjecture}[glicci conjecture]\label{glicciconjecture}
     Let $R$ be a regular local (graded local) ring and let $I$ be a perfect ideal (homogeneous perfect ideal) of positive height. Then $I$ is in the Gorenstein linkage class (in the homogeneous Gorenstein linkage class) of a complete intersection.
\end{conjecture}
\begin{question}[Extended glicci question]\label{extendedglicciquestion}Let $R$ be a Gorenstein local (graded local) ring and let $I$ be a perfect ideal (homogeneous perfect ideal) of positive height. Then $I$ is in the  Gorenstein linkage class (in the homogeneous Gorenstein linkage class) of a complete intersection.  
\end{question} 
In Section~3, we give examples of Cohen-Macaulay ideals that are not glicci.  Over regular rings, the glicci conjecture is known to hold for many classes of ideals. Some notable examples are Gorenstein ideals \cite[Proposition~7.1]{CasanellasDrozdHartshorne+2005+149+171}, Artinian monomial ideals in $k[[x_1,\ldots,x_n]]$ \cite[Theorem~4.2]{huneke2007liaison}, weakly vertex decomposable squarefree monomial ideals \cite[Theorem~3.3]{nagel2008glicci}, Cohen--Macaulay Borel-fixed monomial ideals \cite[Theorem~3.3]{migliore2002monomial}, determinantal ideals \cite[Corollary~3.2]{gorla2008generalized}, residual intersections of licci ideals \cite[Theorem~4.2]{hartshorne2012residual}, ladder Pfaffian ideals \cite[Corollary~2.4]{de2009g}, toric ideals of graphs for which $\mathcal U(I_G)$ consists of quadratic binomials \cite[Theorem~6.11]{cummings2023geometric}, and polarizations of stable monomial ideals \cite[Theorem~3.15]{faridi2025polarization}.

Since not every Cohen-Macaulay ideal over a Gorenstein ring is glicci, it is natural to ask what class of ideals can replace complete intersections in this setting. In Section~4, we prove that every unmixed ideal over a Gorenstein ring can be Gorenstein linked in two steps to a Golod ideal. 
\begin{recalledtheorem}{GolodGorensteinlinkage}
    Let $R$ be a Gorenstein local ring, and let $I$ be a proper unmixed ideal with $\height I=g\ge 2$. Then there exists an unmixed ideal $D$ of height g such that $I\sim_G I'\sim_G D$ and the natural projection $R\x R/D$ is a Golod homomorphism. Here $\sim_G$ stands for Gorenstein link.     
\end{recalledtheorem}

For unmixed ideals of heights one and two, we show that one can reach a Golod ideal by two complete intersection links. This raises the question of whether every unmixed ideal is in the complete intersection linkage class of some Golod ideal.

In Section~7, we prove that $\mathfrak m$-primary ideals over Gorenstein local rings of dimension at least two with an infinite residue field are glicci, and establish corresponding homogeneous results.
\begin{recalledtheorem}{mprimaryglicci}
Let $(R,\mathfrak m, k)$ be a local Gorenstein ring with an infinite residue field and $\dim R\ge 2$, or let $(R=R_0\oplus\bigoplus_{i>0}R_i,\mathfrak m=\mathfrak m_0\oplus\bigoplus_{i>0} R_i,k)$ be a standard graded-local Gorenstein ring over $R_0$, with $\gr R_nR\ge 2$ for all $n\gg 0$, $(R_0,\mathfrak m_0)$ a Artinian local ring and $|k|=\infty$.  Let $I$ be a (homogeneous) $\mathfrak m$-primary ideal. Then $I$ is (homogeneously) $2$-glicci; in particular, $I$ is (homogeneously) glicci. 
    
\end{recalledtheorem}

This gives a large class of Cohen-Macaulay ideals that are glicci without
necessarily being perfect.  It also raises the question of how glicciness behaves under deformation. In fact, the glicci conjecture is equivalent to the deformation question posed in Question~\ref{deformationquestion}. 

The proof uses techniques developed by Casanellas, Drozd, and Hartshorne, which we review in Section~6. We
give an algebraic treatment of their results and extend them to the setting considered here. Our main theorem, Theorem~\ref{2dimDL}, shows that every orientable maximal Cohen-Macaulay (MCM) module over a two-dimensional Gorenstein local ring is stably isomorphic to a doubly layered module. We then use this theorem to obtain the results in Section~9. 

A classical motivation for the study of complete intersection linkage
is its use in the classification of algebraic curves, particularly in $\mathbb P^3_k$. One other motivation for developing this approach is to provide tools for related questions about curves, including
when a curve can be realized as the set-theoretic intersection of two surfaces. Hartshorne, Eisenbud, Schreyer, Casanellas, Miró-Roig, Notari, Lesperance, Nagel, Migliore, and others have studied Gorenstein linkage of points and curves in $\mathbb{P}^n_k$ extensively; see \cite{hartshorne2001some, hartshorne2001experiments, hartshorne2008codimension, eisenbud2015twenty, migliore2013glicci, miro2000gorenstein, lesperance2001gorenstein, Hartshorne2004, nagel2000even, notari2006even, hartshorne2008liaison}. We extend all the existing results in their full generality to show that any two finite sets of reduced points are in the same Gorenstein liaison class, and any two curves are in the same Gorenstein liaison class.  

\begin{recalledtheorem}{mainsection9}
Let $S$ be a polynomial ring over an infinite field, and let $I\subsetneq S$ be a homogeneous ideal. Suppose that $S/I$ is Cohen-Macaulay.
Then $I$ is homogeneously glicci in either of the following cases:
\begin{enumerate}
\item[$(a)$] $\dim S/I=1$.
\item[$(b)$] $\dim S/I=2$.
\end{enumerate}
\end{recalledtheorem}

 These results motivate the study of glicci ideals under deformations. The behavior of licci ideals under deformations and essential deformations was studied in \cite{huneke1987structure}. Licciness also ascends under deformation. In contrast, our examples in Section~3, together with the result on $\mathfrak m$-primary ideals, show that the analogous statement for glicci ideals is false in general. 

In Section~10, we prove that Cohen-Macaulay monomial ideals are glicci. This has been a central problem in the field combinatorial commutative algebra. Linkage of monomial ideals add to the tool kit to study simplicial complexes.  
\begin{recalledtheorem}{thm:main}
Let $S=k[x_1,\ldots,x_n]$ be a standard graded polynomial ring, where $k$ is an infinite field,
and let $I\subset S$ be a proper monomial ideal such that $S/I$ is
Cohen-Macaulay and $\height I=g\geq3$. There is a homogeneous
$S$-regular sequence $f_1,\ldots,f_{g-2}\in I$ such that, on setting
\[
 K=(f_1,\ldots,f_{g-2}),\qquad R=S/K,\qquad J=I/K,
\]
the ideal $J$ is homogeneously glicci in $R$. In particular, $I$ is
homogeneously glicci in $S$.  
\end{recalledtheorem}

We obtain stronger results for Cohen-Macaulay monomial ideals that are clean or have height three: these
ideals are 2-glicci. These results extend earlier work on monomial linkage; see \cite{migliore2002monomial, migliore2013glicci, faridi2025polarization, mantero2025slightly, lyle2025ordinary, mudrak2026zero} to the general case. This raises the further question of whether every Cohen-Macaulay monomial ideal is 2-glicci.

In Section~11, we prove that any ideal is in the  Gorenstein linkage class of an ideal with deviation two. This gives a reduction in the number of generators
within the Gorenstein linkage class. This is particularly interesting because, if the glicci conjecture is false, this theorem determines the smallest possible deviation among ideals in the Gorenstein linkage class of a given ideal.
\begin{recalledtheorem}{reductionofdeviation}
    Let $S$ be a local Gorenstein ring with an infinite residue field. Let $I$ be a proper Cohen-Macaulay ideal with $\height g\ge 3$ and $\type S/I\ge 2$. If $\mu(I)\ge g+3$, then there exists a Cohen-Macaulay ideal $J$ with $\mu(J)\le g+2$ that is Gorenstein linked to $I$ by a sequence
of perfect Gorenstein ideals.
\end{recalledtheorem}

Most of the results in this paper are proved by reducing to height-two ideals
and studying their N-type resolutions. This approach expresses linkage questions in terms of filtrations of associated MCM modules. It avoids
constructing explicit Gorenstein linking ideals at each step, a task that often requires methods adapted to
particular classes of ideals. The resulting module theoretic framework brings the low-dimensional and monomial cases into a common
approach to Gorenstein linkage.

Section~2 collects the basics of various tools used in the article. We reprove some results from basic element theory in a more general setting. This is a very common tool in both algebraic geometry and commutative algebra. These results are classical; the presentation is intended to provide a convenient modern reference for future
work in commutative algebra and algebraic geometry. In Section~5, we extend Rao’s theorem \cite{rao1981liaison} to a setting in which the ambient ring need not be normal or a
domain. This version is used in Section~6. It has been a long standing open problem in commutative algebra and algebraic geometry to classify the complete intersection linkage class of varieties. This theorem classifies the complete intersection linkage class of height 2 Cohen-Macaulay ideals; and for other heights the problem remains widely open. 

Finally, in Section~12, we collect the open questions that arise throughout the paper. These questions are discussed in detail in their respective sections, together with the results that motivate them. For the convenience of the reader, we bring them together in one place to highlight the main questions that remain open and the possible directions for further research.

\subsection{Notation and conventions}
\begin{itemize}
    \item[$(i)$] In the local case, $(R,\mathfrak m,k)$ denotes a Noetherian local ring with an infinite residue field $k$. 
    \item[$(ii)$] In the graded case, $R=\bigoplus_{j\geq 0}R_j$ denotes a Noetherian positively graded ring such that $(R_0,\mathfrak m_0,k)$ is a local ring with an infinite residue field, and $\mathfrak m=\mathfrak m_0\oplus\bigoplus_{i>0}R_i$. We assume $\gr R_nR\geq 2$ for all $n\gg0$. All modules, ideals, and homomorphisms are understood to be graded or homogeneous in the graded setting.    
    \item[$(iii)$] For a finitely generated $R$-module $M$, we write $\Omega_R^i(M)$ for its $i^\text{th}$ syzygy in a minimal free resolution of $M$, $M^*=\h_R(M,R)$ for its dual, and $\mu_R(M)$ for its minimal number of generators. In the graded setting, we use the corresponding graded notions.
    
    \item[$(iv)$] An ideal $I\subseteq R$ is \emph{unmixed} if $\height \mathfrak p=\height I$ for every $\mathfrak p\in\Ass_R(R/I)$.
    
    \item[$(v)$] An ideal $I\subseteq R$ is \emph{Cohen-Macaulay} (resp.\ \emph{Gorenstein}) if $R/I$ is Cohen-Macaulay (resp.\ Gorenstein).
    
    \item[$(vi)$] If $I$ is Cohen-Macaulay of grade $g$, we write $\type(I)=\type(R/I)=\mu_R(\Ext_R^g(R/I,R))$. 
\end{itemize}

\subsection{Acknowledgments}
The author would like to thank his advisor, Bernd Ulrich, for many helpful discussions and suggestions, as well as for his guidance and support throughout the writing of this paper.  He is also grateful to Claudia Polini for valuable discussions and for carefully proofreading the manuscript.

\subsection{AI Disclosure}
All mathematical results, arguments, and proofs in this paper are due to the author. The ideas of studying doubly layered modules were original ideas of the author, who has been working on this project as part of his PhD thesis since December 2025. AI was used to find references for basic element theory, orientable modules, Golod homomorphism, and monomial co-resolution. ChatGPT was used extensively to improve the exposition and language of the manuscript, to polish and shorten the author's original proofs, and to assist in locating proofs from older articles in the literature. ChatGPT identified several gaps in the author's original proofs in Section~9, most of which were subsequently addressed by the author, in some cases with the assistance of AI.

Theorem~\ref{thm:construction} is the only result for which ChatGPT provided assistance in completing the proof. The underlying ideas (Lemma~\ref{lem:propagation}, Lemma~\ref{lem:specialize}) were developed by the author while studying deformations. ChatGPT was prompted with this and a sketch of proof similar to the height 3 case, using Shamash resolutions to get all syzygies in $\mathcal{P}_i$ inductively. It pointed out the construction in step 1 of the proof, and the rest followed immediately. Though the original answer given by it was not working, author altered the part~(b) to get it working.  

The original proofs of all other results in this paper were written by the author and later polished with the help of ChatGPT. Some proofs in Section~3 was substantially rewritten with the help of ChatGPT, this is explicitly acknowledged before the corresponding result is stated.

\section{Preliminaries} 
Let \(R\) be a Noetherian ring with total quotient ring \(Q(R)\), and let $M$ be a finitely generated $R$-module. If \(M\otimes_R Q(R)\) is a free \(Q(R)\)-module, we define $\rank_R M:=\rank_{Q(R)}(M\otimes_R Q(R))$; otherwise, \(\rank_R M\) is not defined. We say that $M$ is \emph{locally free of constant rank $r$ in depth 1} if $M_\mathfrak p$ is a free $R_\mathfrak p$ module of rank $r$ for all $\mathfrak p\in\Spec R$, with $\dep R_\mathfrak p\le 1$. Notice that $\rank_R M=r$ if and only if $M_{\mathfrak p}$ is a free $R_{\mathfrak p}$-module of rank $r$ for every $\mathfrak p\in\Ass(R)$. 

\subsection{Graded Basic Element Theory} All the results mentioned in this section are well known for local rings with an infinite residue field. To avoid repetition, we prove only the corresponding graded versions.

For the rest of this section, let \((R=R_0\oplus\bigoplus_{i>0}R_i,\mathfrak m,k)\) be a positively graded-local Noetherian ring, where \((R_0,\mathfrak m_0)\) is a local ring with an infinite residue field \(k=R_0/\mathfrak m_0\) and for all $n\gg 0$, we have $\gr (R_nR)\ge 2$, where $R_nR$ denotes the ideal generated by $R_n$. We keep the infinite residue field assumption throughout the paper to choose general elements simultaneously in the basic element argument. The assumption on $\gr(R_nR)\ge 2$ is needed to ensure that the graded version of basic element theory works; see Example~\ref{basicelementcounterexample}. Let \(M\) be a finitely generated $\mathbb Z$-graded \(R\)-module. We say that $M$ is generated in degrees at most $a$ if it admits
a finite homogeneous generating set $m_1,\ldots,m_s$ with
$\deg m_i\le a$ for every $i$. 

\begin{remark}
\label{rem:homogeneous-choice}
Let $R=\bigoplus_{n\ge0}R_n$ be a Noetherian positively graded ring,
where $(R_0,\mathfrak m_0)$ is local with an infinite residue field.
Suppose that $\operatorname{grade}(R_nR)\ge2$ for every $n\ge n_0$,
where $n_0\ge1$.
Let $N$ be a finitely generated $\mathbb{Z}$--graded $R$-module, generated in
degrees at most $a$, and fix $d\ge a+n_0$. Then
\[
(RN_d)_{\mathfrak p}=N_{\mathfrak p}
\qquad\text{whenever }\operatorname{depth}R_{\mathfrak p}\le1.
\]
Moreover, for any finite collection $\mathcal P$ of primes
with $N_{\mathfrak p}\ne0$ and $\depth R_p\le 1$, there exists $x\in N_d$ whose image in
$N_{\mathfrak p}/\mathfrak pN_{\mathfrak p}$ is nonzero for every
$\mathfrak p\in\mathcal P$. 

\end{remark}

\begin{proof}
Choose homogeneous generators $u_1,\ldots,u_m$ of $N$, with
$\deg u_i=a_i\le a$. Fix a prime $\mathfrak p$ with
$\operatorname{depth}R_{\mathfrak p}\le1$. Since $d-a_i\ge n_0$,
the ideal $R_{d-a_i}R$ 
is not contained in $\mathfrak p$. Now, for each $i$,
choose $b_i\in R_{d-a_i}\setminus\mathfrak p$. The elements
$b_iu_i$ belong to $N_d$ and generate $N_{\mathfrak p}$.
This proves the first assertion.

Now, fix the finite collection $\mathcal P$. If $\mathcal P$ is
empty, take $x=0$. Otherwise, for each $\mathfrak p\in\mathcal P$,
choose a finite set of degree-$d$ elements
whose images generate $N_\mathfrak p$. Combine
these to obtain a finite set  $v_1,\ldots,v_s\in N_d$ such that their images
generate $N_{\mathfrak p}$ for every $\mathfrak p\in\mathcal P$.

Choose an infinite set $D\subseteq R_0$ consisting of
representatives of distinct residue classes modulo
$\mathfrak m_0$. If $c,c'\in D$ are distinct, then $c-c'$ is a unit
in $R_0$, and hence in every $R_{\mathfrak p}$. Hence, distinct elements
of $D$ have distinct images in every residue field
$\kappa(\mathfrak p)=R_{\mathfrak p}/\mathfrak pR_{\mathfrak p}$.
For $c\in D$, set
\[
x(c)=v_1+cv_2+\cdots+c^{s-1}v_s\in N_d.
\]
Fix $\mathfrak p\in\mathcal P$; then there are at most $s-1$ values of
$c\in D$ for which $x(c)$ has zero image in
$N_{\mathfrak p}/\mathfrak pN_{\mathfrak p}$. Indeed, if $s$
distinct values $c_1,\ldots,c_s$ gave zero image, the resulting
Vandermonde system would have determinant
$\prod_{i<j}(\overline c_j-\overline c_i)\ne0$ in
$\kappa(\mathfrak p)$. Its invertibility would force the images
of all $v_i$ to be zero; hence $N_\mathfrak p=0$, which leads to a contradiction.

Since $\mathcal P$ is finite, only finitely many values in $D$
are excluded. We may therefore choose $x=x(c)$ with the required properties. 
\end{proof} 

 We say that \(x\in M\) is \emph{basic} (\emph{basic in depth one}) if $\mu_{R_\mathfrak p}\bigl(M_\mathfrak p/R_\mathfrak px\bigr)<\mu_{R_\mathfrak p}(M_\mathfrak p)$ for every \(\mathfrak p\in\Spec R\) with \(\depth R_\mathfrak p\leq1\). 

The following lemma is usually attributed to Bourbaki, where it first appeared in \cite{Bourbaki1961Algebre}; an English translation can be found in \cite[Chapter~7, Section~4.9, Theorem~6]{bourbaki1972commutative}. Bourbaki’s theorem is stated for normal domains; the normality assumption was later removed by Eisenbud and Evans in \cite[Theorem~A]{eisenbud1973generating}. We refer to \cite{bruns1976jede} for a proof in the local case, as it uses relatively more modern notation. See \cite{sanyal2010translation} for the English translation.  Here we prove only the graded version.

\begin{remark}\cite[Hilfssatz]{bruns1976jede}\cite[Lemma]{sanyal2010translation}\label{rem:bruns-finiteness}
Let $\mathfrak a$ be a proper ideal of a Noetherian ring $R$.
There are only finitely many primes $\mathfrak p\supseteq\mathfrak a$
such that $\dep R_{\mathfrak p}=\operatorname{grade}\mathfrak a$.
Consequently, an ideal containing an $R$-regular element is contained
in only finitely many primes of depth at most one.
\end{remark}

\begin{lem}
\label{lem:bourbaki-rank-reduction}
Let $R$ be as above, and let $M$ be a finitely generated torsion-free $\mathbb{Z}$-graded $R$-module, and locally free of constant rank $r>1$ in depth one. Suppose that $M$ is
generated in degrees at most $a$. For every $d\ge a+n_0$, there is
a graded exact sequence
\[
0\longrightarrow R(-d)\longrightarrow M\longrightarrow N
\longrightarrow0,
\]
where $N$ is torsion-free, locally free of constant rank $r-1$
in depth one,
and generated in degrees at most $a$.
\end{lem}
\begin{proof}
Fix $d\ge a+n_0$. By Remark~\ref{rem:homogeneous-choice}, choose
$x\in M_d$ whose image in
$M_{\mathfrak p}/\mathfrak pM_{\mathfrak p}$ is nonzero for every
$\mathfrak p\in\operatorname{Ass}R$.

Set $I=\operatorname{Fitt}_{r-1}(M/Rx)$. At a prime
$\mathfrak p$ with $\dep R_{\mathfrak p}\le1$, the module
$M_{\mathfrak p}$ is free of rank $r$. This quotient $(M/Rx)_{\mathfrak{p}}$ needs
$r$ generators if the residue of $x$ in $M_\mathfrak{p}/\mathfrak{p}M_\mathfrak{p}$ 
is zero, and $r-1$ otherwise.
Thus $I\subseteq\mathfrak p$ precisely when that residue is zero.
In particular, $I$ is contained in no associated prime of $R$,
so prime avoidance gives an $R$-regular element in $I$.
By Remark~\ref{rem:bruns-finiteness}, the set $B$ of depth-one
primes where $x$ has zero residue is finite.

The quotient $M/Rx$ is generated in degrees at most $a$ and is
nonzero at every prime in $\operatorname{Ass}R\cup B$, since
$r>1$. Apply Remark~\ref{rem:homogeneous-choice} to this quotient
at these primes and lift the resulting degree-$d$ element to
$y\in M_d$. At every associated prime, the residues of $x,y$
are linearly independent. At every prime in $B$, the residue
of $y$ is nonzero. Therefore, the residues of $x,y$ are never
both zero at a prime of depth at most one.

Now set $J=\operatorname{Fitt}_{r-2}(M/(Rx+Ry))$. At a prime
of depth at most one, the quotient needs at most $r-2$
generators precisely when the residues of $x,y$ are linearly
independent. Hence, $J$ is contained in no associated prime,
and the depth-one primes where these residues are dependent
are exactly those containing $J$. Hence $J$ contains an $R$-regular element, and the set of depth-one primes containing $J$ is finite by Remark~\ref{rem:bruns-finiteness}.

Remark~\ref{rem:bruns-finiteness}.

Choose an infinite set $D\subseteq R_0$ of representatives of
distinct classes modulo $\mathfrak m_0$. Distinct elements of
$D$ differ by a unit. At any of the finitely many primes where
the residues of $x,y$ are dependent, at most one $c\in D$ can
satisfy $x+cy\in\mathfrak pM_{\mathfrak p}$: two such choices,
upon subtraction, would force both $y$ and $x$ into
$\mathfrak pM_{\mathfrak p}$, a contradiction. Choose $c\in D$
avoiding these finitely many exceptions. At all remaining primes
of depth at most one, the residues of $x,y$ are independent,
so $x+cy$ also has a nonzero residue. Thus, $m=x+cy\in M_d$
satisfies $m/1\notin\mathfrak pM_{\mathfrak p}$ whenever
$\dep R_{\mathfrak p}\le1$.

Define $\varphi:R(-d)\to M$ by $\varphi(1)=m$.
If $\dep R_{\mathfrak p}\le1$, then $m/1$ extends to a basis
of $M_{\mathfrak p}$, so $\varphi_{\mathfrak p}$ is injective. So, $\ker\varphi=0$, and its cokernel $N$ satisfies
$N_{\mathfrak p}\cong R_{\mathfrak p}^{r-1}$
whenever $\dep R_{\mathfrak p}\le1$.

Finally, torsion-freeness and rank $r$ give
$M\hookrightarrow Q(R)^r$. Clearing denominators gives
an ungraded embedding $M\hookrightarrow R^r$, so
$\operatorname{depth}M_{\mathfrak p}\geq 1$ whenever
$\operatorname{depth}R_{\mathfrak p}\geq 1$. At primes
of depth at least two, the depth lemma gives $\operatorname{depth}N_{\mathfrak p}
\geq
\min\{\operatorname{depth}M_{\mathfrak p},
      \operatorname{depth}R_{\mathfrak p}-1\}
\geq 1.$ Together with the local freeness already proved, \cite[Proposition~1.4.1(a)]{bruns1998cohen} implies that $N$ is
torsionless, hence torsion-free.
\end{proof}

\begin{lem}
\label{lem:graded-bourbaki}\label{existenseofbasicelements}\label{gradedbasicelements}

Let $R=\bigoplus_{n\ge0}R_n$ be a Noetherian positively graded ring,
where $(R_0,\mathfrak m_0)$ is local with an infinite residue field.
Suppose that $\operatorname{grade}(R_nR)\ge2$ for every $n\ge n_0$,
where $n_0\ge1$.
Let $M$ be a
finitely generated torsion-free $\mathbb Z$-graded $R$-module, generated in degrees
at most $a$, that is locally free of constant rank $r\geq1$ in depth
one. Fix $d\geq a+n_0$.  Then there is a graded exact
sequence
\[
0\longrightarrow F\longrightarrow M\longrightarrow I(h)
\longrightarrow0,
\qquad F=\bigoplus_{i=1}^{r-1}R(-d),
\]
where $I$ is a homogeneous ideal of $R$ and $h\in\mathbb Z$.
\end{lem}
The result follows by induction on rank using Lemma~\ref{lem:bourbaki-rank-reduction} and the fact that, under these hypotheses, every rank-one graded torsion-free module is a shift of a homogeneous ideal.
We leave the details to the reader.


The conclusion of Lemma \ref{gradedbasicelements} can fail
if the hypothesis $\gr(R_nR)\geq2$ for all sufficiently large $n$ is
replaced by $\gr R_+\geq2$

\begin{ex}\label{basicelementcounterexample}
Let $R=k[x,y,z]$, where $|k|$ is infinite and
$\deg x=\deg y=\deg z=2$. Then $\operatorname{grade}R_+=3$,
but $R_n=0$ for every odd $n$. Let $J=(x,y)$ and
$M=J\oplus J(-1)$. This module is torsion-free of rank two
and locally free in depth one.

Every homogeneous element of $M$ lies in one summand, since
$J$ and $J(-1)$ are concentrated in even and odd degrees,
respectively. WLOG if $m=(f,0)\ne0$, then
$M/Rm\cong(J/fR)\oplus J(-1)$ has nonzero torsion. Thus no graded
exact sequence $0\to R(-d)\to M\to I(h)\to0$ exists, since
every ideal of $R$ is torsion-free.
\end{ex}

 \subsection{Orientable Modules} \begin{defi}
    Let $R$ be a Noetherian graded ring, and let $M$ be a finitely generated graded $R$-module that is locally free of constant rank $r\ge 1$ in depth 1. We define the determinant of $M$ as \[\det(M):=(\bigwedge^r M)^{**}.\] We call $M$ \emph{orientable} if there is a graded $R$-module isomorphism $\det(M)\cong R(d)$ for some $d\in\mathbb{Z}$.
\end{defi}
We also set $\det(0) = R$ so that the zero module, which occurs as a minimal syzygy of a free module, is
orientable. Note that this definition is equivalent to the classical definition of orientability found in \cite{bruns1986orientations}, which defines $M$ to be orientable if it is locally free of constant rank r in depth one and there exists an $R$-linear map $\mu:\bigwedge^r M\x R$ such that $\gr\img \mu\ge 2$. The following lemma is well known in the literature; for the convenience of the reader, we include a proof.

\begin{lem}\label{orientableexactsequences}
Let $R$ be a Noetherian graded ring, and let $
0\to A \to B \to C \to 0$ 
be a short exact sequence of finitely generated graded $R$-modules that are locally free of constant positive rank in depth 1.
Then there is a natural graded isomorphism
\[
\det(B)\ \cong\ \big(\det(A)\otimes_R \det(C)\big)^{**}.
\]
Consequently, if any two of $A,B,C$ are orientable, then so is the third.
\end{lem}

\begin{proof}
The second statement will follow from the first. Set \(a=\rank A\), \(c=\rank C\), and \(b=\rank B=a+c\). We use the following observation: if \(\varphi:M\to N\) is a graded map of finitely generated modules such that \(\varphi_\mathfrak p\) is an isomorphism for every \(\mathfrak p\) with \(\depth R_\mathfrak p\leq1\), then both $\varphi^*$ and \(\varphi^{**}\) are isomorphisms. 

For \(i\geq0\), let \(F^i\subseteq\bigwedge^bB\) be the image of the natural graded map \(\bigwedge^iA\otimes_R\bigwedge^{b-i}B\to\bigwedge^bB\). There is a natural graded map
$\theta:\bigwedge^aA\otimes_R\bigwedge^cC\longrightarrow F^a/F^{a+1}$,
where, for homogeneous \(u\in\bigwedge^aA\) and homogeneous \(v=c_1\wedge\cdots\wedge c_c\), we choose homogeneous lifts \(b_i\in B\) of \(c_i\) and send \(u\otimes v\) to the class of \(u\wedge b_1\wedge\cdots\wedge b_c\). This is well defined since, changing a lift changes the wedge by an element of \(F^{a+1}\).

Let \(\mathfrak p\) satisfy \(\depth R_p\leq1\). Then the localized sequence splits
as \(C_\mathfrak p\) is free. So, $(F^a)_\mathfrak p=\bigwedge^bB_\mathfrak p$, $(F^{a+1})_\mathfrak p=0$, 
the inclusion \(F^a\to\bigwedge^bB\), the quotient \(F^a\to F^a/F^{a+1}\), and \(\theta\) are all isomorphisms after localizing at \(\mathfrak p\). By the observation above, we obtain $
\left(\bigwedge^bB\right)^{**}\cong
\left(\bigwedge^aA\otimes_R\bigwedge^cC\right)^{**}$.

The canonical maps \(\bigwedge^aA\to\det(A)\) and \(\bigwedge^cC\to\det(C)\) are also isomorphisms in depth \(1\). Hence, their tensor product becomes an isomorphism after taking double duals, and $\left(\bigwedge^aA\otimes_R\bigwedge^cC\right)^{**}
\cong
\left(\det(A)\otimes_R\det(C)\right)^{**}$. Combining these gives $\det(B)\cong\left(\det(A)\otimes_R\det(C)\right)^{**}$.
\end{proof}
Note that $\det(\bigoplus_{i=1}^r R(d_i))\cong R(\sum_{i=1}^r d_i)$. By the lemma, we get that syzygies of an orientable module are orientable. It is important to note the following remarks. 

\begin{remark}
Let \(R\) be a Noetherian graded ring, and let \(M\) be a finitely generated graded \(R\)-module that is locally free of constant rank \(r\ge 1\) in depth one. If \(M\) is orientable, then \(M^*\) is orientable.
\end{remark}

\begin{proof}
Suppose \(\det(M)\cong R(d)\), and let \(N=\bigwedge^rM\), and \(D=N^{**}=\det(M)\). 

Let \(\eta:N\to D=N^{**}\) be the canonical graded map. If \(\depth R_\mathfrak p\leq1\), then \(M_\mathfrak p\) is free of rank \(r\), so \(N_\mathfrak p\cong R_\mathfrak p\), and hence \(\eta_\mathfrak p\) is an isomorphism. By the observation made in the previous proof, we conclude that \(\eta^*:D^*\to N^*\) is an isomorphism. Since \(D\cong R(d)\), we obtain $N^*\cong D^*\cong R(-d)$.

Now consider the natural graded determinant pairing
\[
\alpha:\bigwedge^r(M^*)\longrightarrow N^*,
\qquad
\alpha(\varphi_1\wedge\cdots\wedge\varphi_r)
(m_1\wedge\cdots\wedge m_r)
=\det(\varphi_i(m_j)).
\]
For every \(\mathfrak p\) with \(\depth R_\mathfrak p\leq1\), \(M_\mathfrak p\) is free of rank \(r\), so \(\alpha_\mathfrak p\) is the usual isomorphism
\(\bigwedge^r(M_\mathfrak p^*)\cong(\bigwedge^rM_\mathfrak p)^*\). Hence \(\alpha^{**}\) is an isomorphism, and therefore $
\det(M^*)=\left(\bigwedge^r(M^*)\right)^{**}
\cong (N^*)^{**}\cong R(-d)$.
\end{proof}

\begin{remark}\label{N+N*isorientable}
Let $R$ be a Noetherian graded ring and let $N$ be a finitely generated graded $R$-module, locally free of constant rank $r>0$ in depth 1. Then $N\oplus N^*$ is orientable.
\end{remark}
\begin{proof}[sketch of proof] Projecting \(\bigwedge^{2r}(N\oplus N^*)\) onto its summand \(\bigwedge^r N\otimes_R\bigwedge^r N^*\) and applying the determinant pairing gives a natural map to \(R\). In depth one, \(N\) is free of rank \(r\), so this is the only nonzero summand and the pairing is an isomorphism. Hence taking double dual we get \(\det(N\oplus N^*)\cong R\).
    
\end{proof}

\begin{remark}Let \(R\) be a graded-local Gorenstein ring, and let \(I\subset R\) be a homogeneous ideal of height at least \(2\). Then \(I\), viewed as an \(R\)-module of rank \(1\), is orientable. \end{remark} 

\subsection{N-type resolutions} 

We say an ideal $I$ of $R$ is Cohen-Macaulay if $R/I$ is a Cohen-Macaulay ring. Let $R$ be a graded-local ring. We call a graded exact sequence of the form $$0\x F\x M\x I(d)\x 0$$ an \emph{N-type resolution} of $I$ if $F$ is graded free, $M$ is graded, orientable, and an MCM $R$-module, and $I$ is a proper homogeneous Cohen-Macaulay ideal of height $2$. This exact definition first appeared in \cite[Chapter~2, Definition~4.5]{perrin1990classification}. Many existing works refer to this as Bourbaki sequences. 

\begin{prop}\cite[Proposition~1.8, Proposition~1.9]{herzog1987maximal}\label{existenceNtype} Let $R=\bigoplus_{n\ge0}R_n$ be a graded-local Gorenstein ring,
where $(R_0,\mathfrak m_0)$ is local with an infinite residue field.
Suppose that $\gr(R_nR)\ge2$ for all $n\gg 0$.

\begin{itemize}
    \item[$\rm{(a)}$] Let $I$ be a homogeneous Cohen-Macaulay ideal
of height two and type $r$. Then, for some
$d_1,\ldots,d_r,h\in\mathbb Z$, there is an $N$-type
resolution
\[
0\longrightarrow \bigoplus_{i=1}^r R(d_i)
\longrightarrow M\longrightarrow I(h)\longrightarrow0.
\]
Moreover, for fixed $h$, this minimal $N$-type resolution is unique
up to graded isomorphism, and every other $N$-type
resolution of $I(h)$ is a trivial extension of it.

    \item[$\rm{(b)}$] Conversely if $M$ be a $\mathbb Z$-graded orientable MCM $R$-module of rank $r\ge 2$. Then there is an N-type resolution $0\x F\x M\x I(h)\x0$. In particular $F$ can be chosen to be $\bigoplus_{i=1}^{r-1} R(-d)$ for some $d\gg0$. 
\end{itemize}
\end{prop}
\begin{proof}
$(a)$ The proof given in \cite[Proposition~1.9]{herzog1987maximal} adapts to graded extensions, with shifts chosen
from the degrees of the minimal generators of $\Ext_R^1(I(h),R)$. 

(b) Let $M$ be generated in degrees at most $a$. Since $M$ is MCM of constant rank over the Gorenstein ring $R$, it is torsion-free. Lemma~\ref{gradedbasicelements} gives a graded exact sequence
\[
0\longrightarrow R(-d)^{r-1}\longrightarrow M\longrightarrow L\longrightarrow0, \qquad\text{for some }d\gg a
\]
where $L$ is torsion-free of rank one and free in depth one. As $d\gg a$, we have $M_d\subseteq R_+M$, so the first map has the image in $\mathfrak mM$. Hence $\mu(L)=\mu(M)\geq r\geq2$, and $L$ is not cyclic.

By Lemma \ref{orientableexactsequences}, $L^{**}\cong R(h)$ for some $h\in\mathbb Z$. Since $L$ is torsion-free of rank one, its evaluation map is injective; identify $L$ with $I(h)$ for a homogeneous ideal $I$. Noncyclicity gives $I\ne R$. In depth one, the evaluation map is an isomorphism, so $I_{\mathfrak p}=R_{\mathfrak p}$; since $R$ is Cohen-Macaulay, $\operatorname{ht}I\geq2$. The depth lemma gives $\operatorname{depth}I\geq\dim R-1$, and therefore $\operatorname{depth}R/I\geq\dim R-2$. Together with $\dim R/I\leq\dim R-2$, this proves that $I$ has height two and $R/I$ is Cohen--Macaulay.
\end{proof}


We say a homogeneous ideal $I$ of $R$ is Gorenstein if $R/I$ is a Gorenstein ring. It is important to note that a Gorenstein ideal need not have a finite projective dimension. Now we apply the proposition for height 2 Gorenstein ideals to record the following.

\begin{cor}\label{rank2Gorenstein}
    Let $R$ be graded-local Gorenstein ring as above. Let $K$ be a proper homogeneous Cohen-Macaulay ideal of height two. Then $K$ is Gorenstein ideal if and only if there exists an N-type resolution of the form  $0\x R(d)\x M\x K(h)\x 0$, where $M$ is a rank 2 MCM orientable $\mathbb Z$-graded $R$-module. 

    Moreover, let \(M\) be an orientable graded MCM \(R\)-module of rank two,
and let \(\alpha\in \mathfrak mM\) be a homogeneous basic element of degree \(d\).
Then \(\alpha\) induces an \(N\)-type resolution
\[
0\longrightarrow R(-d)\xrightarrow{\,1\x\alpha\,}
M\longrightarrow K'(\ell)\longrightarrow 0,
\]
where \(K'\) is a homogeneous Gorenstein ideal of height two.
\end{cor}

The reverse implication in the first assertion follows from graded duality $\Ext_R^1(K(h),R)\cong \Ext_R^2(R/K,R)(-h)$, which is a graded quotient of \(R(-d)\); hence \(R/K\) has type one and is
Gorenstein. For the moreover statement, the usual Bourbaki argument shows
that \(M/R\alpha\) is, up to a shift, a height-two Cohen-Macaulay ideal
\(K'\), which is Gorenstein by the first assertion.

\subsection{Fundamentals of Linkage} Throughout this section, let $R$ be a Noetherian $\mathbb Z$-graded Gorenstein ring with
a unique homogeneous maximal ideal. An ideal $K$ is Gorenstein if $R/K$ is a Gorenstein ring. We refer to \cite{schenzel1982notes} for most of the proofs. Although this article works with local Gorenstein rings, all the proofs can be extended to graded-local Gorenstein rings.

\begin{defi}
Let \(R\) be a graded-local Gorenstein ring, and let \(I,J\subsetneq R\) be homogeneous ideals of height \(g\). We say that \(I\) and \(J\) are \emph{homogeneously Gorenstein linked} by a homogeneous Gorenstein ideal \(K\) if $K:I=J$ and $K:J=I$.
\end{defi}

\begin{prop}
Let \(R\) be a graded-local Gorenstein ring, and let \(I,J\subsetneq R\) be homogeneous ideals directly Gorenstein linked by a homogeneous Gorenstein ideal \(K\). Then
\[
\Ass_R(R/I)\cup\Ass_R(R/J)=\Ass_R(R/K).
\]
In particular, if \(\height K=g\), then \(I\) and \(J\) are unmixed ideals of height \(g\).
\end{prop}

\begin{prop}\cite[Proposition~2.2]{schenzel1982notes}\label{symetriclinkage}
Let \(R\) be a Gorenstein graded-local ring, let \(I\subsetneq R\) be a homogeneous unmixed ideal of height \(g\), and let \(K\subsetneq I\) be a homogeneous Gorenstein ideal of height \(g\). Set \(J=K:I\). Then \(J\) is unmixed of height \(g\) and $
K:J=I$. Thus \(I\) and \(J\) are homogeneously Gorenstein linked.
\end{prop}
\begin{prop}\cite[Corollary~3.3]{schenzel1982notes}
    Let \(R\) be a $\mathbb{Z}$-graded local Gorenstein ring, and let \(I,J\subsetneq R\) be homogeneous ideals that are homogeneously Gorenstein linked by a Gorenstein ideal K. Then $R/I$ is Cohen-Macaulay if and only if $R/J$ is Cohen-Macaulay.
    
    If \(R/I\) is Cohen-Macaulay and \(\omega_{R/K}\cong (R/K)(a)\), then $\omega_{R/I}\cong (J/K)(a)$, and $\omega_{R/J}\cong (I/K)(a)$ as graded modules.
\end{prop}

\begin{prop}\cite{peskine1974liaison}
    Let $R$ be a graded-local Gorenstein ring, and let $I$ and $J$ be homogeneous ideals that are homogeneously Gorenstein linked by a Gorenstein ideal $K$. If $K$ is perfect, that is, $\operatorname{pd}_R(R/K)=\gr K
    $. Then $I$ is perfect if and only if $J$ is perfect.
\end{prop}
This follows from the mapping cone construction due to Ferrand. In particular, licci ideals are perfect. 

\begin{defi}
Let \(R\) be a graded-local Gorenstein ring, and let \(I\subsetneq R\) be a homogeneous unmixed ideal.

\begin{enumerate}
\item[(a)] The \emph{homogeneous Gorenstein linkage class} of \(I\) is the set of all homogeneous ideals \(J\subsetneq R\) that can be obtained from \(I\) by a finite sequence of homogeneous Gorenstein links. The \emph{even homogeneous Gorenstein linkage class} of \(I\) is defined similarly, using an even number of homogeneous Gorenstein links.

\item[(b)] We say that \(I\) is \emph{homogeneously glicci} (in the Gorenstein linkage class of a complete intersection) if \(I\) and a complete intersection ideal are in the same homogeneous Gorenstein linkage class.
\end{enumerate}
\end{defi}

\begin{thm} \cite{apery1945courbes}, \cite{gaeta1952quelques},\cite[Theorem~3.2]{peskine1974liaison}
    Let $R$ be a positively graded-local Gorenstein ring, with an infinite residue field and $\gr R_nR\ge 2$ for all $n\gg0$. Let $I$ be a proper homogeneous ideal of grade 2. Then $I$ is homogeneously licci $\iff I$ is perfect.    
\end{thm}

The local version of the theorem holds over any Gorenstein local ring and can be proved by induction on the deviation. The same argument works in the standard graded setting. For more general grading, however, the local proof does not carry over directly; instead, we use Theorem~\ref{rao}.

\subsection{Technical remarks}

The following standard facts will be used throughout the paper. We assume $R$ is a Noetherian positively graded ring with a unique homogeneous maximal ideal. All modules
are finitely generated and $\mathbb Z$-graded, and all maps in exact sequences
and diagrams have degree zero. We give
$M^*=\operatorname{Hom}_R(M,R)$ its natural grading, so
$\operatorname{Hom}_R(M,R)_t=\operatorname{Hom}_{\mathrm{gr}R}(M,R(t))$.
The modules $\operatorname{Ext}^i_R(M,N)$ carry the corresponding grading.

\begin{lem}\label{split}
Let
\[
 (*)\qquad 0\longrightarrow A\oplus P\xrightarrow{i}X\xrightarrow{q}E
 \longrightarrow0
\]
be a graded exact sequence, where $P$ is graded projective and
$\operatorname{Ext}^1_R(E,P)_0=0$. 
Then $(*)$ is noncanonically isomorphic, as a graded extension,
to the direct sum of an extension $0\to A\to X'\to E\to0$ and the split
sequence $0\to P\xrightarrow{\mathrm{id}}P\to0\to0$.
\end{lem}

\begin{lem}[Graded horseshoe lemma]
Let $0\to A\xrightarrow{\iota}B\xrightarrow{\pi}C\to0$ be a graded exact
sequence, and let $F_\bullet\twoheadrightarrow A$ and
$G_\bullet\twoheadrightarrow C$ be graded free resolutions. Then there is
a graded free resolution $H_\bullet\twoheadrightarrow B$ with
$H_n=F_n\oplus G_n$ and differentials
\[
d_n^H=\begin{pmatrix}d_n^F&\sigma_n\\0&d_n^G\end{pmatrix}
\qquad(n\geq1),
\]
where $\sigma_n:G_n\to F_{n-1}$ is homogeneous of degree zero.
Consequently, for every $n\geq0$, there is a finite graded free module
$P_n$ and a graded exact sequence
\[
0\longrightarrow\Omega_R^n(A)\longrightarrow
\Omega_R^n(B)\oplus P_n\longrightarrow\Omega_R^n(C)
\longrightarrow0.
\]
\end{lem}

\begin{lem}
Suppose that $R$ is Gorenstein. Then $\operatorname{Hom}_R(-,R)$ is an
exact contravariant functor on the category of graded maximal
Cohen-Macaulay $R$-modules. Moreover, if $M$ is a graded MCM module, then
$M^*$ is MCM, and the canonical graded map $M\to M^{**}$ is an isomorphism.
\end{lem}

\begin{remark}
Let $M$ be a graded $R$-module of rank $r$. If $\mu_R(M)=r$, then $M$ is graded free.
\end{remark}

\begin{remark} \label{minimalsyzygydualresolution}
Suppose that $R$ is Gorenstein, and let $M$ be a graded module with no graded free direct summand. If $\operatorname{Ext}^1_R(M,R)=0$,
then $
\mu_R(M)=\mu_R\bigl((\Omega_R(M))^*\bigr)$.
\end{remark}

We note the following remark on pushouts and pullbacks along short exact sequences. 

\begin{remark}
Work in the category of graded $R$-modules with degree-zero maps.

\textit{(i) Pushout.}
The pushout $P$ of $0\to A\to B\to C\to0$ and
$0\to A\to B'\to C'\to0$ gives exact sequences
\[
0\to B'\to P\to C\to0,
\qquad
0\to B\to P\to C'\to0.
\]
\textit{(ii) Pullback.}
The pullback $P$ of $0\to A\to B\to C\to0$ and
$0\to A'\to B'\to C\to0$ gives exact sequences
\[
0\to A'\to P\to B\to0,
\qquad
0\to A\to P\to B'\to0.
\]
\end{remark}
 
 Let $u_\bullet:G_\bullet\x F_\bullet$ be a morphism of complexes of graded free modules with $G_i=0=F_i$ for all $i<0$. 
The \emph{mapping cone of u}, denoted by $C(U_\bullet)$, is a complex given by $$ \cdots F_{i+1}\oplus G_{i}\xrightarrow{}F_{i}\oplus G_{i-1}\xrightarrow{}F_{i-1}\oplus G_{i-2}\cdots$$

where differentials are $C(u)_{i}=F_i\oplus G_{i-1}\xrightarrow{\begin{bmatrix}
    \partial_i^F & u_{i-1}\\
    0    & -\partial_{i-1}^G
\end{bmatrix}}C(U)_{i-1}=F_{i-1}\oplus G_{i-2}$

\begin{remark}\label{mappingcone}
    \begin{itemize}
         \item[(a)] $0\x F_\bullet\hookrightarrow{}C(u_\bullet)\xrightarrow{proj}G[-1]\x 0$ is an exact sequence of complexes.
        \item[(b)] If $H_i(F_\bullet) = H_i(G_\bullet) = 0$ for all $i > 0$ and $H_0(u_\bullet) : H_0(G_\bullet) \rightarrow H_0(F_\bullet)$ is injective, then $C(u_\bullet)$ is acyclic, that is $H_i(C(u_\bullet)) = 0$ for all $i > 0$ and $H_0(C(u_\bullet)) \cong \text{coker } H_0(u_\bullet)$.
        
    \end{itemize}
\end{remark}

The following lemma is essentially a consequence of the horseshoe lemma, pullbacks, and the mapping-cone.

\begin{lem}\label{lem:syz}
The following statements hold for graded modules.
\begin{enumerate}
\item[$(a)$] If $0\to U\to V\to W\to0$ is exact, then for each $q\geq1$ there are graded free modules $H_q,G_q$ and exact sequences
\begin{align}
0&\x\Omega_R^qU\x\Omega_R^qV\oplus H_q\x\Omega_R^qW\x0, \tag{1}\\
0&\x\Omega_R^{q+1}W\x\Omega_R^qU\oplus G_q\x\Omega_R^qV\x0. \tag{2}
\end{align}
\item[$(b)$] If $u$ is homogeneous and $T$-regular, then for $q\geq2$, there is a graded free module $H$ and an exact sequence
\begin{equation}\tag{3}
0\x\Omega_R^qT\x\Omega_R^q(T/uT)\oplus H\x\Omega_R^{q-1}T(-\deg u)\x0.
\end{equation}
\item[$(c)$] If $0\to F\to E\to\mathfrak a\to0$ is exact with $F$ free, and $T=R/\mathfrak a$, then
\[
\Omega_R^qT\cong_{st}\Omega_R^{q-1}E\qquad(q\geq2).
\]
\end{enumerate}
\end{lem}

\begin{proof}[Sketch of proof]
The horseshoe resolution of $V$ gives (a)(1). For the second, choose a graded free cover $0\to\Omega_RW\to F\to W\to0$ and pull it back along $V\twoheadrightarrow W$. Since the resulting sequence $0\to U\to P\to F\to0$ splits, $P\cong U\oplus F$, and hence $0\to\Omega_RW\to U\oplus F\to V\to0$ is exact. Applying the first construction to this sequence gives $0\to\Omega_R^{q+1}W\to\Omega_R^qU\oplus G_q\to\Omega_R^qV\to0$ for a graded free module $G_q$.

For (b), take a free resolution of $T$ and the mapping cone of multiplication by $u$. It resolves $T/uT$, since $0\to T(-\deg u)\xrightarrow{u}T\to T/uT\to0$ is exact. Remark~\ref{mappingcone}(a), Remark~\ref{mappingcone}(b), and the horseshoe lemma give (3). Part (c) is standard.
\end{proof}

\section{Cohen–Macaulay Ideals That Are Not Glicci}
\subsection{Height 1 examples}
We show that Cohen-Macaulay does not imply glicci in height one by constructing a Gorenstein ideal that is not glicci. Let us first classify the complete intersection linkage classes of ideals of height one.

\begin{lem}\label{CIlinkageclassesofheight1}
Let \(R\) be a Gorenstein local ring, and let \(I,J\subset R\) be proper height-one unmixed ideals. Then \(I\) and \(J\) are in the same even complete intersection linkage class if and only if \(I\cong J\) as \(R\)-modules. They are in the same odd complete intersection linkage class if and only if \(I\cong J^*\).
\end{lem}

We give a sketch of the proof. 
\begin{proof}[Sketch of proof]
If \(A\) is an unmixed height-one ideal and \(x\in A\) is \(R\)-regular, then for \(B=(x):A\) multiplication gives \(B\cong A^*\); the double-colon property gives \(A=(x):B\), hence \(A\cong A^{**}\). Thus each direct CI-link replaces an ideal by its dual, so an even chain preserves the module up to isomorphism, while an odd chain gives the dual. Conversely, if \(I\cong J\), then $J=qI$ for some unit $q$ in the total quotient ring. Choose an $R$-regular element $x\in\mathfrak m I$ and set $y=qx\in\mathfrak m J$, then we have \((x):I=(y):J\), giving a common direct link and hence an even chain. If \(I\cong J^*\), link \(I\) once to an ideal $K$, but now \(I^*\cong K\cong J\), and apply the even case.
\end{proof}

In particular, two evenly CI-linked ideals can be linked in two steps, while two oddly CI-linked ideals can be linked in at most three steps. Now we are ready to present the example.

\begin{ex} The original proof for the following was substantially revised by ChatGPT;

Let $R=k[[x,y,z,w]]/(xw-yz)$, then $R$ is Gorenstein normal local domain. Take $I=(x,z)R$, now $I$ is a height one Gorenstein ideal of $R$. Now observe $(x,y)R=(x)R:(x,z)R$. 
Hence, by the lemma above $I^*\cong (x)R:I=(x,y)R$. Now we have the following remark. 

\begin{remark}
Let $R=k[[x,y,z,w]]/(xw-yz)$, and set $I=(x,z)R$. Let $L\subset R$ be a height-one ideal with $L\cong I$ or $L\cong I^*$. If $K\subsetneq L$ is a height-one Gorenstein ideal, then $K$ is principal. 
\end{remark}

\begin{proof}
We prove the case $L\cong I$; the case $L\cong I^*\cong (x,y)R$ is the same, after interchanging $y$ and $z$. Since $L\cong I$ and both are rank-one ideals, $L=qI$ for some $q\in \operatorname{Frac}(R)$. As $L\subset R$, we have $q\in R:I$. Since $(x)R:I=(x,y)R$, we get $R:I=x^{-1}(x,y)R=R+Rt$, where $t=y/x=w/z$. Hence $q=a+bt$ for some $a,b\in R$, and so $L=(u,v)R$, where $u=ax+by$ and $v=az+bw$. Let $S=k[[x,y,z,w]]$, let $g=xw-yz$, and let $A\subset S$ be the inverse image of $K$. Then $g\in A\subset (g,u,v)$. Since $R/K\cong S/A$ is Gorenstein and $S$ is regular, $A$ is a height-two Gorenstein ideal of $S$, hence a complete intersection. In particular, $\mu(A)=2$.

If $g\notin \mathfrak n A$, where $\mathfrak n=(x,y,z,w)$, then $g$ is part of a minimal generating set of $A$, so $A=(g,h)$ for some $h\in S$. Hence $K=(h)R$ is principal. Thus we may assume $g\in \mathfrak n A$. If $a,b\in \mathfrak n$, then $u,v\in \mathfrak n^2$, so $A\subset (g,u,v)\subset \mathfrak n^2$, hence $\mathfrak n A\subset \mathfrak n^3$, contradicting the fact that $g$ has order $2$. Therefore one of $a,b$ is a unit. Then the identities $wu-yv=ag$ and $zu-xv=-bg$ show that $g\in (u,v)$, so the inverse image of $L$ is $(u,v)$.

Now $u,v$ have independent linear parts. Since $A\subset (u,v)$ is two-generated and $g\in \mathfrak n A$, if the image of $A$ in $(u,v)/\mathfrak n(u,v)$ were one-dimensional, then the initial form of $g$ would be divisible by a nonzero linear form in the span of the initial forms of $u$ and $v$. This is impossible: after a linear change of variables, the initial form of $g=xw-yz$ has the form $UV'-U'V$, which is not divisible by any nonzero linear form in the span of $U,V$. Hence the image of $A$ in $(u,v)/\mathfrak n(u,v)$ is two-dimensional, so $A+\mathfrak n(u,v)=(u,v)$. By Nakayama's lemma, $A=(u,v)$, and therefore $K=L$, contradicting $K\subsetneq L$. Thus the only proper case is $A=(g,h)$, and then $K=(h)R$ is principal.
\end{proof}

Suppose, to the contrary, that \(I\) is glicci, and consider a sequence of Gorenstein links from \(I\) to a complete intersection. By induction on the links, Lemma 3.1 and Remark 3.3 show that every link in the sequence is a complete intersection link. Thus \(I\) is licci, hence perfect, a contradiction since height one perfect ideals are principal. Therefore \(I\) is not glicci.
\end{ex}
We next give an example in dimension one. If \(R\) is a one-dimensional normal local domain, then \(R\) is a DVR, so $\mathfrak m$-primary ideals are complete intersections. We therefore consider a nonnormal Gorenstein local domain.

\begin{ex}
Let $k$ be a field with $\operatorname{char} k\neq 2$ and there exist a $\delta\in k$ which is not a square. Set
\[
R=k[[u,v]]/(v^2-\delta u^2),\qquad \mathfrak m=(u,v)R.
\]
Then $R$ is a one-dimensional Gorenstein local domain, but it is not normal, also $\mathfrak m$ is a height-one Gorenstein ideal. We claim that $\mathfrak m$ is not glicci. First observe that $\mathfrak m^*\cong \mathfrak m$, this follows because $(u)R:\mathfrak m=\mathfrak m$. 

The proof of the remark is completed with the help of ChatGPT.
\begin{remark}
    Let $L\subset R$ be a height-one ideal with $L\cong \mathfrak m$. If $K\subsetneq L$ is a height-one Gorenstein ideal, then $K$ is principal.
\end{remark}
\begin{proof}
  Since $L\cong \mathfrak m$ and both are rank-one ideals, $L=q\mathfrak m$ for some $q\in \operatorname{Frac}(R)$. As $L\subset R$, we have $q\in R:\mathfrak m=u^{-1}\mathfrak m=R+Rt$, where $t=v/u$. Thus $q=a+bt$ for some $a,b\in R$, and hence
\[
L=(au+bv,\; av+\delta bu)R.
\]
Let $S=k[[u,v]]$, let $f=v^2-\delta u^2$, and let $A\subset S$ be the inverse image of $K$. Since $R/K\cong S/A$ is Gorenstein and $S$ is regular, $A$ is a height-two Gorenstein ideal of $S$, hence a complete intersection. Thus $\mu(A)=2$ and $f\in A$. If $f\notin (u,v)A$, then $A=(f,h)$ for some $h\in S$, and so $K=(h)R$ is principal. Suppose $f\in (u,v)A$. If $a,b\in (u,v)$, then $L$ has inverse image contained in $(u,v)^2$, so $A\subset (u,v)^2$, whence $(u,v)A\subset (u,v)^3$, contradicting the fact that $f$ has order $2$. Therefore one of $a,b$ is a unit. Since $\delta$ is not a square, $a^2-\delta b^2$ is a unit, and so $au+bv,\;av+\delta bu$ generate $(u,v)$. Thus the inverse image of $L$ is $(u,v)$. Now, if the image of $A$ in $(u,v)/(u,v)^2$ were one-dimensional, then the initial form $v^2-\delta u^2$ would be divisible by a linear form over $k$, impossible because $\delta$ is not a square. Hence $A+(u,v)^2=(u,v)$, and by Nakayama $A=(u,v)$. Thus $K=L$, a contradiction. Hence every proper height-one Gorenstein ideal contained in any $L\cong\mathfrak m$ is principal.
\end{proof}

Now again by similar arguments to the previous example we see that $\mathfrak m$ is not glicci. 
\end{ex}

\subsection{Height 2 
examples}
We now construct height-two Cohen-Macaulay ideals that are not glicci. In Section~7, we show that every height-two Cohen–Macaulay ideal in a two-dimensional Gorenstein local ring is glicci. The following examples occur in higher dimensions. As in the height-one case, we consider nonregular Gorenstein local rings over which every height-two Gorenstein ideal is a complete intersection. Over such rings, a height-two Cohen–Macaulay ideal is glicci if and only if it is licci, equivalently, if and only if it is perfect. Thus it suffices to construct nonperfect height-two Cohen–Macaulay ideals.

 
\begin{cor}\cite[3.4]{herzog1987maximal}\label{herzogrankisolatedsing} Let $R$ be a hypersurface local domain, i.e., $R\cong S/(f)$ where $S$ is a regular local ring and $(f)$ is a proper nonzero prime ideal of $S$. Suppose $R$ has an isolated singularity, i.e., $R_\mathfrak p$ is a regular local ring for all $\mathfrak p\in\Spec R\backslash \{\mathfrak m\}$, where $\mathfrak m$ is the maximal ideal of $R$. Let $M$ be a MCM $R$-module, which is not free, then $$2\rank(M)+1\ge \dim R.$$   
\end{cor}

This follows from a result of Bruns \cite[Corollary~2]{bruns1981eisenbud} and the 2 periodic resolutions of MCM modules over hypersurfaces. 

It follows that if \(R\) is an isolated hypersurface singularity of dimension at least six, then every rank-two MCM \(R\)-module is free. By the local version of Corollary~\ref{rank2Gorenstein}, \cite[Proposition~1.9]{herzog1987maximal}, every height-two Gorenstein ideal of \(R\) is a complete intersection. Consequently, a height-two Cohen–Macaulay ideal of \(R\) is glicci if and only if it is perfect.

\begin{thm}\cite[Corollary~1.3]{dutta1989syzygies}\label{dutta} Let $(R,\mathfrak m, k)$ be a local Noetherian ring. Suppose the $n^{th}$ syzygy module of $k$ has a non-zero projective module as its direct summand for some
$n\ge 0$. Then $R$ is regular.
\end{thm}

\begin{ex}\label{isolatedsingularityexample}
     Let $R$ be a non-regular, hypersurface local ring, with an isolated singularity, and $\dim R=d\ge 6$. Now $\Omega^d(k)$ is an orientable, MCM module, so consider an N-type resolution of it, $0\x F\x \Omega^d(k)\x I\x 0$, where $F$ is a free $R$-module. So we get $I$ is a height 2 Cohen-Macaulay ideal, that is not perfect: if not then we get $\Omega^d(k)$ has finite projective dimension, hence it is free contradicting $R$ not being regular ring. Hence $I$ is not glicci.

    Over this ring $R$,  let $M$ be MCM orientable module, that is not free over $R$. Then any  ideal corresponding to N-type resolution of $M$ is not glicci over $R$.
\end{ex}

\begin{ex}
    Let $R=\mathbb{R}[[x_1,x_2,x_3,x_4]]/(x_1^2+x_2^2+x_3^2+x_4^2)$, from \cite[Example~3.3]{herzog1987maximal} we get all rank 2 MCM orientable module of $R$ is free. So again by the local version of Corollary~\ref{rank2Gorenstein} every height 2 Gorenstein ideal is a complete intersection. So height 2 Cohen-Macaulay ideal is glicci iff licci iff perfect. So again one can consider the N-type ideals of a non free MCM orientable modules, then these ideals are not perfect hence they are not glicci. Note that $R$ is 3 dimensional non-regular, hypersurface local domain, so there exists non free MCM orientable modules; for example, $\Omega^3_R(k)$. 
\end{ex}

It is interesting to note that by Proposition~\ref{Quadratic hypersurfacesingularities} in $R'=\mathbb{C}[[x_1,x_2,x_3,x_4]]/(x_1^2+x_2^2+x_3^2+x_4^2)$, every Cohen-Macaulay ideal of height two is glicci in $R'$. Indeed, after a linear change of coordinates,
\(R'\cong\mathbb C[[u,v,z,w]]/(uv-zw)\).
Kn\"orrer's classification gives
\(M\cong(R')^a\oplus L_+^{\oplus b}\oplus L_-^{\oplus c}\)
for every MCM \(R'\)-module \(M\), where \(L_+=(u,z)\) and \(L_-=(u,w)\).
These ideals represent opposite generators of
\(\operatorname{Cl}(R')\cong\mathbb Z\), so orientability forces \(b=c\).
Thus, every orientable MCM module is stably a direct sum of the
type~(i) module \(L_+\oplus L_-\), and Theorem~\ref{DLstablyGlicci} gives the assertion.

This shows that glicciness need not descend along faithfully flat extensions of residue fields. But by the following remark, we get that licci ideals descend along the faithful flat extensions of residue field. 

Let $(R,\mathfrak m,k)$ be a Noetherian local ring, and let $k\subset l$ be a field extension. Then there exists a Noetherian local ring $(R',\mathfrak m')$ and a faithfully flat local homomorphism $R\to R'$ such that $\mathfrak m'=\mathfrak mR'$ and $R'/\mathfrak m'\cong l$. The proof of the existence of $R'$ is standard; see  \cite[\href{https://stacks.math.columbia.edu/tag/03C3}{Tag 03C3}]{stacks-project}. From the construction, one can see that $R'$ is Noetherian (Cohen-Macaulay, Gorenstein) when $R$ is Noetherian (Cohen-Macaulay, Gorenstein).

\begin{remark}
   Let $R$ be a Gorenstein local ring, with an infinite residue class field $k$, and let $R'$ be as above, for some field extension $k\seq l$. Let \(I\subset R\) be a Cohen-Macaulay ideal. If \(IR'\) is licci in \(R'\), then \(I\) is licci in \(R\). 
\end{remark}

\begin{proof}
By \cite[Proposition~2.3]{johnson2014ubiquity}, since \(IR'\) is licci, there exists \(n\) such that \(L_n(IR')=S_n'\).
Universal links commute with flat base change, hence \(L_n(IR')=L_n(I)S_n'\).
Thus \(L_n(I)S_n'=S_n'\), and faithful flatness gives \(L_n(I)=S_n\).
Again by \cite[Proposition~2.3]{johnson2014ubiquity}, \(I\) is licci.
\end{proof}

Here is another example over non normal Gorenstein local ring. 

\begin{ex}
Let \(A=k[[x_0,\ldots,x_6]]/(q)\), where \(q=x_0^2+\cdots+x_6^2\), \(k\) is an infinite field of characteristic different from \(2\), and let \(I_0\subset A\) be a non perfect height-two Cohen-Macaulay Bourbaki ideal as in Example~\ref{isolatedsingularityexample}, arising from an exact sequence \(0\x A^{r-1}\x M\x I_0\x0\) with \(M\) an orientable, non free MCM \(A\)-module. Set
\[
R=A[[u,v]]/(uv)\qquad\text{and}\qquad I=I_0R.
\]
Since \(A\) is a domain, \(u,v\neq0\) in \(R\), whereas \(uv=0\); hence the local ring \(R\) is not a domain and therefore is not normal. Moreover,
\[
R/I\cong A[[u,v]]/(I_0A[[u,v]],uv)\cong (A/I_0)[[u,v]]/(uv),
\]
so \(R/I\) is Cohen--Macaulay and \(\operatorname{ht}_R I=\operatorname{ht}_A I_0=2\).

Let \(\mathfrak p=\mathfrak m_AR+(u)\). Since \(v\notin\mathfrak p\), \(v\) becomes a unit in \(R_{\mathfrak p}\), and the relation \(uv=0\) gives \(u=0\). Hence
\[
B:=R_{\mathfrak p}\cong\bigl(A[[v]][v^{-1}]\bigr)_{\mathfrak m_AA[[v]][v^{-1}]}
 \cong\left(\frac{k[[x_0,\ldots,x_6,v]][v^{-1}]}{(q)}\right)_{(x_0,\ldots,x_6)}.
\]
Thus \(B\) is a six-dimensional hypersurface local ring with residue field \(k((v))\), and
\[
\widehat B\cong k((v))[[x_0,\ldots,x_6]]/(q).
\]

The Jacobian ideal of \(q\) in \(B\) is \((x_0,\ldots,x_6)B\), so \(B\) has an isolated singularity. The map \(A\to B\) is flat, being the composite of the flat power-series extension and localizations, and it is faithfully flat since it is local and \(\mathfrak m_AB\neq B\).

Now \(I_{\mathfrak p}=I_0B\). If \(I_0B\) were perfect, then the base change
\[
0\longrightarrow B^{r-1}\longrightarrow M\otimes_AB\longrightarrow I_0B\longrightarrow0
\]
would imply that \(M\otimes_AB\) has finite projective dimension. It is MCM over \(B\), hence by Auslander-Buchsbaum formula it is free. Faithfully flat descent would then imply that \(M\) is free over \(A\), a contradiction. Therefore \(I_{\mathfrak p}\) is not perfect. From the above examples we get $I_\mathfrak p$ is not glicci. Now, by Corollary~\ref{Gorensteinidealsareglicci} Gorenstein ideals of height atleast 2 are glicci, so we get that glicci-ness localizes. Hence, \(I\) is not glicci in \(R\).
\end{ex}

We also note that \(R\) admits height-two Gorenstein ideals that are not complete intersections. The examples in this section naturally lead to the following
question.

\begin{question}\label{whenisCMglicci}
Under what additional assumptions on a local Gorenstein
ring $R$ is every Cohen-Macaulay ideal of $R$ glicci?
\end{question}

The examples in this section suggest a common obstruction to glicciness: the lack of sufficiently many Gorenstein ideals of the appropriate height to construct nontrivial links. It would be interesting to understand the existence and abundance of such ideals and their effect on Gorenstein linkage classes.

\section{Linkage classes of Golod ideals}
A Cohen-Macaulay ideal over a Gorenstein local ring need not be glicci. This raises the question of whether there is a class of ideals to which every Cohen-Macaulay ideal is linked. A natural motivation for this question comes from the following theorem.

\begin{thm}\cite[Theorem~1.2,3.7]{christensen2020linkage} Let $I$ be a grade 3 perfect ideal in a regular local ring $Q$. Then $I$ is in the complete intersection linkage class of a grade 3 perfect ideal $J$, such that $Q/J$ is a complete intersection or a Golod ring.     
\end{thm}

The Poincaré series of a finitely generated module $M$ over a local ring $(R,\mathfrak m,k)$ is defined as $P_M^R(t)=\sum_{n=0}^\infty \dim_k\Tor_n^R(M,k)t^n$. Let $(S,\mathfrak n,k)$ be a regular local ring, and $(R,\mathfrak m,k)$ a local ring with a surjective local homomorphism $S\x R$. Then Serre in \cite{serre1965algebre} showed that \[P_k^R(t)\le \frac{P_k^S(t)}{1-t(P_R^S(t)-1)}\] A modern proof for this can be found in \cite[Proposition~3.3.2]{avramov1998infinite}. A local ring is called \emph{Golod} if equality holds for a minimal regular
presentation of its completion. 
\begin{defi}
A surjective homomorphism of local rings $\phi:S\twoheadrightarrow R$ is \emph{Golod} if $P_k^R(t)= {P_k^S(t)}/{(1-t(P_R^S(t)-1))}$.
The source S need not be regular.
\end{defi}

\begin{thm}\cite[Theorem~3.1]{gover1980class}\label{GoverSalmon}
Let \((R,\mathfrak m)\) be a local ring with ideals $I_1\subseteq \cdots \subseteq I_n\subseteq \mathfrak m$, 
and suppose that \(x_1,\ldots,x_n\in \mathfrak m\) satisfy: \(x_1\) is regular on \(R\), and, for each \(i\geq 2\), \(x_i\) is regular on $R_{i-1}=R/\sum_{j<i}x_jI_j$. 
Let \(M\) be a finitely generated \(R\)-module such that \(I_nM=0\). Then $P_M^{R_n}(t)=\frac{P_M^R(t)}{1-t\bigl(P_{R_n}^R(t)-1\bigr)}$. 

In particular, the natural map \(R\to R_n\) is a Golod homomorphism.
\end{thm}

\begin{prop}
    Let $(R,\mathfrak m)$ be a Gorenstein local ring, and let $I$ be an unmixed proper ideal of height one or two. Then there exists an unmixed proper ideal $J$ such that the natural projection map $R\x R/J$ is a Golod homomorphism and $I\sim I'\sim J$. Here $\sim$ denotes a complete intersection link.
\end{prop}
\begin{proof}\emph{Height 1}: Let $a\in I\seq \mathfrak m$ be an $R$-regular element. Set $J=aI$, then $I\cong J$. 
As $I$ and $(a)$ are unmixed, we conclude that $J$ is a proper unmixed ideal of $R$. Hence, by Lemma \ref{CIlinkageclassesofheight1}, we have $I$ and $J$ are complete intersection linked in two steps, i.e., $I\sim I'\sim J$. Apply Theorem~\ref{GoverSalmon} with $n=1$, $I_1=I$, and $x_1=a$. It follows that $R\x R/aI$ is Golod.\\

\emph{Height 2:} Let \(I\) be an unmixed ideal of height two, let \(x,y\in I\) be an \(R\)-regular sequence, and set \(z=x+y\). Let \(B=(xz)\) and \(J=B+yI\). By taking the filtration \((z)\subseteq I\), and as \(y\) is regular on \(R/(xz)\), by the theorem above, we get that \(R\to R/J\) is a Golod homomorphism. Set \(S=R/(xz)\) and \(\mathfrak a=I/(xz)\). Choose an \(S\)-regular element \(c\in\mathfrak a\), put \(\mathfrak a'=(c):_S\mathfrak a\), and let \(I'\subset R\) be its inverse image. Let \(\widetilde c\in I\) be a lift of \(c\); then \(G_1=(xz,\widetilde c)\) and \(G_2=(xz,y\widetilde c)\) are complete intersections since \(c\) and \(y\) are \(S\)-regular. Moreover, \((yc):_S\mathfrak a'=y\mathfrak a\), since \((c):_S\mathfrak a'=\mathfrak a\) and \(y,c\) are regular. Hence the double-colon property gives \(G_1:I=I'\), \(G_1:I'=I\), \(G_2:J=I'\), and \(G_2:I'=J\). Therefore \(I\sim I'\sim J\).
\end{proof}

\begin{lem}\cite[Proposition~4.1]{huneke2007liaison}\label{Gorensteindoublelink} Let $(R,\mathfrak m)$ be a Gorenstein local ring, $J$ and $K$ be proper ideals of $R$, and $x\in\mathfrak m$. Assume $R/J$ is a Cohen-Macaulay and generically Gorenstein ring, and that the ideal $J+xK$ is unmixed of height greater than $\height J$. Then $J+xK$ and $J+K$ are Gorenstein doubly linked (Gorenstein linked in 2 steps).   
\end{lem}

\begin{thm}\label{GolodGorensteinlinkage}
    Let $R$ be a Gorenstein local ring, and let $I$ be a proper unmixed ideal with $\height I=g\ge 2$. Then there exists an unmixed ideal $D$ of height g such that $I\sim_G I'\sim_G D$ and the natural projection $R\x R/D$ is a Golod homomorphism. Here $\sim_G$ stands for Gorenstein link. 
\end{thm}

\begin{proof}
Set \(d=\dim R\), choose an \(R\)-regular sequence \(z_1,\ldots,z_{g-1},x\) contained in \(I\). For \(1\leq i\leq g-1\), set \(y_i=z_i+x\), \(A_i=(z_1,\ldots,z_i)\), and \(B_i=\sum_{j=1}^{i}y_jA_j\). Let \(B=B_{g-1}\) and define \(D=B+xI\). We prove that \(D\) is unmixed ideal of $R$, the natural projection \(R\to R/D\) is a Golod homomorphism, and that \(I\) and \(D\) are connected by two Gorenstein links.

For \(0\leq r\leq i\), set \(C_{r,i}=(z_1,\ldots,z_r,y_{r+1},\ldots,y_i)\), then $C_{r,i}$ is a complete intersection ideal of height $i$. We claim, by induction on \(i\), that \(R/B_i\) is Cohen-Macaulay of dimension \(d-i\), and every minimal prime of \(B_i\) is minimal over one of the ideals \(C_{r,i}\), and \(B_i\) is generically a complete intersection. Set \(B_0=(0)\), suppose that the assertions hold for \(B_{i-1}\), and let \(\mathfrak p\in\Min(R/B_{i-1})\). By induction, \(\mathfrak p\) is minimal over some \(C_{r,i-1}\). Since \(C_{r,i-1},y_i\) is a regular sequence of length \(i\), whereas \(\height\mathfrak p=i-1\), one has \(y_i\notin\mathfrak p\). As \(R/B_{i-1}\) is Cohen-Macaulay, it has no embedded associated primes, and therefore \(y_i\) is regular on \(R/B_{i-1}\).

Since \(B_{i-1}\subseteq A_i\), we have \(B_i:y_i=(B_{i-1}+y_iA_i):y_i=(B_{i-1}:y_i)+A_i=B_{i-1}+A_i=A_i\). Hence multiplication by \(y_i\) gives an exact sequence \(0\to R/A_i\xrightarrow{\cdot y_i}R/B_i\to R/(B_{i-1},y_i)\to0\). The module \(R/A_i\) is Cohen-Macaulay of dimension \(d-i\), while \(R/(B_{i-1},y_i)\) is Cohen--Macaulay of the same dimension because \(y_i\) is regular on \(R/B_{i-1}\). The depth lemma therefore shows that \(R/B_i\) is Cohen-Macaulay of dimension \(d-i\).

Let \(\mathfrak p\in\Min(R/B_i)\). Then \(\height\mathfrak p=i\). If \(y_i\notin\mathfrak p\), the containment \(y_iA_i\subseteq B_i\subseteq\mathfrak p\) implies \(A_i\subseteq\mathfrak p\), so \(\mathfrak p\) is minimal over \(C_{i,i}=A_i\). If \(y_i\in\mathfrak p\), choose \(\mathfrak q\in\Min(R/B_{i-1})\) with \(\mathfrak q\subseteq\mathfrak p\). By induction, \(\mathfrak q\) is minimal over some \(C_{r,i-1}\), and hence \(C_{r,i}=C_{r,i-1}+(y_i)\subseteq\mathfrak p\). Since \(C_{r,i}\) and \(\mathfrak p\) both have height \(i\), the prime \(\mathfrak p\) is minimal over \(C_{r,i}\). Moreover, \(B_i\subseteq C_{r,i}\) for every \(r\): Indeed, every generator of \(B_i\) has the form \(z_\ell y_j\), with \(1\leq \ell\leq j\leq i\); if \(\ell\leq r\), then \(z_\ell\in C_{r,i}\), while if \(\ell>r\), then \(j\geq \ell>r\), hence \(y_j\in C_{r,i}\).

Fix \(\mathfrak p\in\Min(R/B_i)\), and suppose that \(\mathfrak p\) is minimal over \(C_{r,i}\). Since \(C_{r,i},x\) is a regular sequence of length \(i+1\), one has \(x\notin\mathfrak p\). For \(j\leq r\), one has \(z_j\in\mathfrak p\) and \(y_j=z_j+x\notin\mathfrak p\), whereas for \(j>r\), one has \(y_j\in\mathfrak p\) and \(z_j=y_j-x\notin\mathfrak p\). Since \(z_jy_j\in B_i\) for every \(j\), localization at \(\mathfrak p\) gives \((C_{r,i})_\mathfrak p\subseteq (B_i)_\mathfrak{p}\), while the reverse inclusion follows from \(B_i\subseteq C_{r,i}\). Hence \((B_i)_\mathfrak p=(C_{r,i})_\mathfrak p\). Thus \(B_i\) is generically a complete intersection.

Taking \(i=g-1\), we obtain \(\height B=g-1\), and \(S:=R/B\) is Cohen--Macaulay and generically Gorenstein. The preceding description of the minimal primes also shows that \(x\) belongs to no minimal prime of \(B\). Since \(R/B\) is Cohen-Macaulay, \(x\) is regular on \(R/B\).

The unmixedness of $D$ follows from the short exact sequence
$0\x R/I\xrightarrow{\cdot x}R/D\x R/(B,x)\x0$. Now, consider the ascending chain \(A_1\subseteq A_2\subseteq\cdots\subseteq A_{g-1}\subseteq I\subseteq \mathfrak m\), together with the elements \(y_1,\ldots,y_{g-1},x\). The element \(y_1\) is \(R\)-regular. For \(2\leq i\leq g-1\), the argument above shows that \(y_i\) is regular on \(R/B_{i-1}=R/\sum_{j<i}y_jA_j\), and \(x\) is regular on \(R/B=R/\sum_{j=1}^{g-1}y_jA_j\). 
Theorem~\ref{GoverSalmon} applies to this chain and gives \(R\to R/D\) is a Golod homomorphism.

Now by the Lemma \ref{Gorensteindoublelink}, we get $D=B+xI$ and $B+I=I$ are Gorenstein linked in 2 steps. 
\end{proof}


\begin{question}\label{questiongolod}
    Is every proper unmixed ideal of height greater than 2, over a Gorenstein local ring in the complete intersection linkage class of some ideal $J$, such that the natural projection map $R\x R/J$ is a Golod homomorphism.
\end{question}

\section{Linkage and N-type resolutions}

We stick to the same notations, but additionally $R$ is a graded-local Gorenstein ring, i.e., $(R=R_0\oplus\bigoplus_{i>0} R_i,\mathfrak m,k)$ is a positively graded-local Gorenstein ring, where $(R_0,\mathfrak m_0)$ is a local ring with $|R_0/\mathfrak m_0|=|k|=\infty$ and $\gr(R_nR)\ge 2$ for all $n\gg 0$, where $R_nR$ is the ideal generated by $R_n$. The reason we are interested in N-type resolution is because of the following remark:

\begin{remark}\label{rmk:heightcutglicci}
Let $R$ be a local Gorenstein ring with infinite residue field, or a graded-local Gorenstein ring as above with $R_0$ Artinian. Let $I\subsetneq R$ be a Cohen-Macaulay ideal of height $g\geq2$, and let $K\subseteq I$ be a Gorenstein ideal with $\operatorname{ht}K\leq g-2$. In the graded case, assume that $I,K$ are homogeneous and that $\operatorname{grade}_{R/K}((R/K)_n(R/K))\geq2$ for all $n\gg0$.

If $I/K$ is (homogeneously) glicci in $R/K$, then $I$ is (homogeneously) glicci in $R$. When $K$ is a complete intersection, the graded assertion holds without the Artinian hypothesis on $R_0$.
\end{remark}

\begin{proof}
Let us note a simple remark:  Let $I$, $J$, $K$ be $R$-ideals such that $K\seq I$, then $\frac{I}{K}:_{\frac{R}{K}} \frac{J}{K}=\frac{I:_R J}{K}$.

So by this observation, one can lift the links of $I/K$ to see that $I$ is in the homogeneous Gorenstein linkage class of a homogeneous Gorenstein ideal. 

Now, if $K$ is a complete intersection, the lifted terminal ideal $J$ is itself a complete intersection, proving the last assertion.

In the local case, choose a complete intersection $C\subseteq J$ of height $g-2$. Then $J/C$ is a height-two Gorenstein ideal, hence is glicci by Corollary~\ref{Gorensteinidealsareglicci}; lifting its links through $C$ shows that $J$, and therefore $I$, is glicci.

In the graded case, since $R_0$ is Artinian, we have $\operatorname{ht}J_+=g$. The construction in the proof of Theorem~\ref{thm:clean-gorenstein-graded}(b) gives a homogeneous complete intersection $C\subseteq J_+$ of height $g-2$ such that $R/C$ satisfies the required grading hypothesis. Thus $J/C$ is homogeneously glicci by Corollary~\ref{Gorensteinidealsareglicci}, and lifting its links through $C$ proves that $I$ is homogeneously glicci.

\end{proof}

Thus, to prove that a Cohen-Macaulay ideal \(I\) of height \(g\geq 2\) is homogeneously glicci, one may find a homogeneous Gorenstein ideal \(K\subseteq I\) of height \(g-2\) and show that \(I/K\) is homogeneously glicci in \(R/K\). This motivates the following definition.

\begin{defi}
Let $R$ be a graded-local Gorenstein ring as above, and let
$I\subsetneq R$ be a homogeneous Cohen--Macaulay ideal of
height $t\ge2$. We say that $I$ is \emph{homogeneously $2$-glicci}
if $I/K$ is homogeneously glicci in $A=R/K$ for every
homogeneous ideal $K\subseteq I$ such that
\[
\operatorname{ht}K\le t-2,\qquad
A\text{ is Gorenstein},\qquad
\operatorname{grade}_A(A_nA)\ge2\quad\text{for all }n\gg0.
\]
\end{defi}

\begin{remark}\label{existenceofgorensteinsubideal}
Let \(I\subseteq R\) be a proper homogeneous Cohen-Macaulay ideal of
height \(g\geq 3\), and let $I^{+}:=\bigoplus_{n>0}I_n$. If
$\height_R I^{+}=g$, then there exists a homogeneous \(R\)-regular sequence
\(f_1,\ldots,f_{g-2}\in I^{+}\) such that, for
\(K=(f_1,\ldots,f_{g-2})\),
\[
\gr_{R/K}\bigl((R/K)_n(R/K)\bigr)\geq 2
\qquad\text{for all }n\gg0.
\]The sequence can be inductively chosen using the homogeneous prime avoidance; see Theorem~\ref{thm:clean-gorenstein-graded}(b). The height condition on \(I^{+}\) is automatic when \(R_0\) is
Artinian local. 

These hypotheses guarantee the existence of a suitable \(K\), but they
do not imply that every homogeneous Gorenstein subideal of height
\(g-2\) has the required property. The latter conclusion does hold
under the stronger assumption
\[
\sqrt{R_nR}=\mathfrak m \qquad\text{for all }n\gg0.
\]
In particular, this applies when \(R\) is standard graded over an
Artinian local ring \(R_0\), since then
\(\sqrt{R_nR}=\sqrt{R_+}=\mathfrak m\) for every \(n\geq1\).

Standardness is necessary for this last assertion. For example, let $R=k[x,z,y]$,$\deg x=\deg z=1$,$\deg y=2$, and let \(I=(x,z,y)\). Although \(R_0=k\) and \(I^{+}=I\), the
Gorenstein subideal \(K=(x)\) does not preserve the grade condition:
for every odd \(n\), $\sqrt{R_nR}=(x,z)
\quad\text{and}\quad
\sqrt{(R/K)_n(R/K)}=(z)$, which has grade one in \(R/K\).
\end{remark}

So we record the following in a remark so that it becomes easy to reference it in the text.

The goal of this section is to prove Rao's theorem on the linkage class of height two Cohen-Macaulay ideals in a non normal, non domain setting. The proof for local rings is essentially the same as \cite[Theorem~2.1]{herzog1987maximal}; they prove this for local Gorenstein domains, but if one carefully reads the proof, one could see that the integral domain hypothesis can be removed by using prime avoidance. In fact, its introduction says that most results hold in non domain setting too. Considering the length of the paper, we provide the proof for it in the graded case only. \cite[Theorem~6.4]{nagel1998even} proves this, but again over integral domains; we closely follow their proof to reprove it in a more general setting. 

For $t\in\mathbb Z$, the graded shift $M(t)$ is defined by $M(t)_n=M_{n+t}$, with the same underlying $R$-module structure. We say that graded $R$-modules $M$ and $N$ are \emph{stably isomorphic up to shift} if there exist finitely generated graded free modules $F,G$ and an integer $t$ such that $M\oplus F\cong N(t)\oplus G$ by a degree-zero isomorphism. If $t=0$ can be chosen, we call them \emph{stably graded-isomorphic}.

\begin{thm}\label{rao}\cite[Theorem~2.1]{herzog1987maximal}\cite[Theorem~6.4]{nagel1998even}\cite{rao1981liaison}
Let $R=\bigoplus_{n\ge0}R_n$ be a positively graded-local Gorenstein ring, where $(R_0,\mathfrak m_0)$ is a local ring with an infinite residue field, and $\operatorname{grade}(R_nR)\ge2$ for all $n\gg0$.

Let
\[
0 \longrightarrow F' \longrightarrow M' \longrightarrow I'(h_1)\longrightarrow 0,\qquad 
0 \longrightarrow F''\longrightarrow M'' \longrightarrow I''(h_2) \longrightarrow 0
\]
be any two N-type resolutions. Then the following
statements are equivalent:
\begin{enumerate}
    \item There exist graded free modules $P,Q$ and an integer $t$ such that $M'\oplus P\cong M''(t)\oplus Q$ as graded modules.
    \item \(I'\) and \(I''\) are evenly homogeneously complete intersection linked.
\end{enumerate}
\end{thm} In particular, when $R$ is of the above format, by Proposition~\ref{existenceNtype} we get that every homogeneous Cohen-Macaulay ideal of height 2 has a N-type resolution. So this classifies complete intersection linkage classes of height 2 ideals. We break the proof into multiple small parts to keep it readable. 

\begin{prop}\cite[Proposition~2.13]{CasanellasDrozdHartshorne+2005+149+171}\label{cone}
Suppose that \(R\) is  graded-local Gorenstein ring (not necessarily positively graded), and that \(I\) has a resolution of the form
\[
0 \longrightarrow \mathcal A \longrightarrow \mathcal B \longrightarrow  I \longrightarrow 0
\]
with \(\mathcal A\) and \(\mathcal B\) graded MCM $R$-modules. Let \(K=(f,g)\) be a height two graded complete intersection in \(R\), contained in \(I\), with $\deg f=a$, and $\deg g=b$. Then \(I'=K:I\) has a resolution
\[
0 \longrightarrow \mathcal B^*(-a-b)
\longrightarrow
\mathcal A^*(-a-b)\oplus R(-a)\oplus R(-b)
\longrightarrow I'
\longrightarrow 0.
\]
\end{prop}
The proof given in \cite{CasanellasDrozdHartshorne+2005+149+171} works in the general graded case too. 
\begin{proof}
    (\ref{rao}) $(2)\implies(1)$ It is enough to prove for 2 links. Let $I'\sim_{\underline{\alpha}}K\sim_{\underline{\beta}} I''$ where $\underline{\alpha}$ and $\underline{\beta}$ are homogeneous linking complete intersection ideals. Now, by the mapping cone, we get the following resolutions for $K$
    $$0\x M'^*(h_1-a_1-a_2)\x F'^*(h_1-a_1-a_2)\oplus R(-a_1)\oplus R(-a_2)\x K\x 0,$$ $$
0\x M''^*(h_2-b_1-b_2)\x F''^*(h_2-b_1-b_2)\oplus R(-b_1)\oplus R(-b_2)\x K\x 0.$$
Hence, by Schanuel's lemma, we get that $M'^*$ and $M''^*$ are stably isomorphic up to shift. As $M'$ and $M''$ are reflexive, we conclude that $M'$ and $M''$ are stably isomorphic up to shift.
\end{proof}

\noindent\emph{Simultaneous choices.}
Let \((\mathfrak p_i,\alpha_i)\), \(1\leq i\leq q\), be finitely many pairs, where \(\operatorname{ht}\mathfrak p_i\leq 1\) and \(\alpha_i:M\to R(t_i)\) are homogeneous maps whose localization is nonzero modulo \(\mathfrak p_i\). For \(d\gg 0\), as \(M_d\) generates \(M_{\mathfrak p_i}\), we can choose \(z_i\in M_d\) with \(\alpha_i(z_i)\notin\mathfrak p_i\). Indeed, we get that 
\(W_i=\{v\in M_d:\alpha_i(v)\in\mathfrak p_i\}\) to be a proper \(R_0\)-submodule of \(M_d\). Since $k$ is infinite, we can find a $v\in M_d$ with $\alpha_i(v)\notin\mathfrak p_i$ for every $i$. 

\begin{lem}\label{lem:replace} Let $R$ be a graded-local Gorenstein ring, as in Theorem \ref{rao}. Let $0\to F\xrightarrow{\iota}M\xrightarrow{\pi}I(h)\to0$ be a N-type resolution. Write $F=Re\oplus F_0$, where $\deg e=a$, and identify $F$ with its image in $M$. Choose $d$ strictly larger than the generating degrees of $M$ and $u,v\in M_d$. Regard $x=\pi(u)$ and $y=\pi(v)$ as elements of $I$, of degree $d+h$, and write
\[
yu-xv=\lambda e+f_0,\qquad f_0\in F_0.
\]
Suppose that $x,y$ and $\lambda,x+cy$ are $R$-regular sequences for some $c\in R_0$. Then $I$ is homogeneously CI linked in two steps to a proper height-two Cohen-Macaulay ideal $J$, with an N-type resolution
\[
0\longrightarrow F_0\oplus R(-d)
\xrightarrow{\ (f,r)\mapsto f+r(u+cv)\ }
M\longrightarrow J(h+d-a)\longrightarrow0.
\]
\end{lem}
\begin{proof}
Set $K=(x,y):I$, and let $w=u+cv$ and $z=x+cy$. Then $yw-zv=\lambda e+f_0$,
and the hypotheses imply that $z,y$ and $z,\lambda$ are regular
sequences. For $m\in M$, let $\rho(m)$ be the coefficient of $e$
in $zm-\pi(m)w\in F=Re\oplus F_0$, and set $J=\rho(M)$.
This defines a homogeneous $R$-linear map $M\x J(h+d-a)$.
Moreover, $\rho(e)=z$, $\rho(v)=-\lambda$, and $\rho$ vanishes
on $F_0$ and $w$.

The identity
\[
y\bigl(zm-\pi(m)w\bigr)-z\bigl(ym-\pi(m)v\bigr)
=-\pi(m)(\lambda e+f_0)
\]
gives $\lambda\pi(m)+y\rho(m)\in zR$, upon taking the coefficient
of $e$. If $\rho(m)=0$, then $\lambda\pi(m)\in zR$. Since
$z,\lambda$ is a regular sequence, hence $\pi(m)=zr$ for some $r\in R$.
As $\pi(w)=z$, it follows that $m-rw\in\ker\pi=F$. Write
$m-rw=se+f$ with $f\in F_0$. Since $\rho(e)=z$ and
$\rho(F_0)=\rho(w)=0$, we obtain $sz=\rho(m-rw)=0$. The
regularity of $z$ gives $s=0$, so $m=rw+f$. Thus
$\ker\rho=Rw+F_0$. This sum is direct: if $rw\in F_0$, then
$rz=\pi(rw)=0$, and hence $r=0$. Therefore
$\ker\rho=Rw\oplus F_0$. Accounting for degrees, we obtain
\[
0\longrightarrow F_0\oplus R(-d)
\xrightarrow{\ (f,r)\mapsto f+rw\ }
M\xrightarrow{\rho}J(h+d-a)\longrightarrow0.
\]

Since $w\in\mathfrak mM$, we have
$\mu(J)\ge\mu(M)-\operatorname{rank}F_0\ge2$, so $J\ne R$.
Also $(z,\lambda)\seq J$. The depth lemma applied to the
sequence above shows that $J$ has height two and $R/J$ is
Cohen-Macaulay.

It remains to identify the second link. Write bars for reduction
modulo $z$. Since $\pi(M)=I$ and $\rho(M)=J$, the inclusion
$\lambda\pi(m)+y\rho(m)\in zR$ for all $m\in M$ gives
$\overline{\lambda}\overline{ I}=-\overline{y}\overline{J}$. Canceling the nonzerodivisors
$\overline{\lambda}$ and $\overline{y}$ shows that
$\overline{r}\overline{I}\seq(\overline{y})$ if and only if
$\overline{r}\overline{J}\seq(\overline{\lambda})$. Consequently,
\[
K=(x,y):I=(z,y):I=(z,\lambda):J.
\]
Now by the double-colon property \ref{symetriclinkage}, we get $I$ and $J$ are linked in two steps,
through $(x,y)$ and $(\lambda,x+cy)$, with the asserted
Bourbaki sequence.
\end{proof}

\begin{remark}\label{rem:ideals-up-to-shift} Let $R$ be as above, and let $I$ and $J$ be homogeneous ideals of grade at least two. If
$I\cong J(t)$ as graded $R$-modules for some $t\in\mathbb Z$, then
$t=0$ and $I=J$.
\end{remark}
It is interesting to note that this simple remark is false for arbitrary $\mathbb Z$ graded-local rings. 

\begin{lem}\label{lem:compare}
Suppose that
\[
0\longrightarrow F_1\longrightarrow M\xrightarrow{\pi_1}I_1(h_1)\longrightarrow0,
\qquad
0\longrightarrow F_2\longrightarrow M\xrightarrow{\pi_2}I_2(h_2)\longrightarrow0
\]
are N-type resolutions. Then $I'$ and $I''$ are evenly homogeneously CI-linked.
\end{lem}
\begin{proof}
Now, both $F_1$ and $F_2$ have rank $s=\rank M-1$. We show how to replace one chosen basis element
$e_i\in F_i$ with the same element of $M$. We will repeatedly use the observation made on simultaneous choices.

For $d\gg0$, the module $M_d$ generates $M_{\mathfrak p}$ for
every prime $\mathfrak p$ of height at most one. Now, as $(I_i)_{\mathfrak p}=R_{\mathfrak p}$ for every $\mathfrak p\in\Ass R$, we get $\pi_i(M)\not\seq \mathfrak p$. Applying the simultaneous choice to the finitely many pairs $(\mathfrak p,\pi_i)$ gives a
common $u\in M_d$ such that $
x_i=\pi_i(u)\notin\mathfrak p$ for all $\mathfrak p\in\Ass R$ and $i=1,2$.
Thus both $x_1$ and $x_2$ are nonzerodivisors.

Write $F_i=Re_i\oplus F_{i,0}$, and let
$\varepsilon_i:F_i\to R$ be the coefficient of $e_i$. Define
\[
\gamma_i(v)=
\varepsilon_i\bigl(\pi_i(v)u-x_iv\bigr).
\]
This is an $R$-linear homogeneous map, and
$\gamma_i(e_i)=-x_i$. In particular, for every
$\mathfrak p\in\Ass R$, the localized map
$(\gamma_i)_{\mathfrak p}:M_{\mathfrak p}\to R_{\mathfrak p}$
is surjective.

We next choose a common $v\in M_d;$ For
$\mathfrak q\in\Ass(R/(x_i))$, we have
$\operatorname{ht}\mathfrak q=1$ and
$(I_i)_{\mathfrak q}=R_{\mathfrak q}$, so
$\pi_i(M_d)\not\seq\mathfrak q$. Also,
$\gamma_i(M_d)\not\seq\mathfrak p$ for every
$\mathfrak p\in\Ass R$. Applying the simultaneous choice argument to
all these finitely many conditions gives $v\in M_d$ such that
\[
y_i=\pi_i(v)\notin\bigcup\Ass(R/(x_i))
\quad\text{and}\quad
\lambda_i=\gamma_i(v)\notin\bigcup\Ass R
\]
for $i=1,2$. Hence $x_i,y_i$ is a regular sequence and
$\lambda_i$ is a nonzerodivisor.

Finally, consider the finitely many primes
$\mathfrak q\in\Ass(R/(\lambda_i))$ for $i=1,2$. Each has height
one. Since $x_i,y_i$ is a regular sequence, $\mathfrak q$ cannot
contain both $x_i$ and $y_i$. Therefore, for each such
$\mathfrak q$, at most one $c\in D$ satisfies
$x_i+cy_i\in\mathfrak q$, where $D\seq R_0$ is a set containing exactly one representative of each residue class in the infinite field $R_0/\mathfrak m_0$. Choose a common $c\in D$ avoiding all
these finitely many values. Then
\[
\lambda_i,\ x_i+cy_i
\]
is a regular sequence for both $i=1,2$.

Lemma~\ref{lem:replace} now replaces $e_1$ and $e_2$ by the same
element $w=u+cv$, using two homogeneous CI links on each side,
while leaving all other basis elements unchanged. Repeating this
construction for the remaining basis elements produces two even
linkage chains whose terminal N-type resolutions have the same
kernel $H\seq M$:
\[
0\longrightarrow H\longrightarrow M\longrightarrow
L_1(a)\longrightarrow0,
\qquad
0\longrightarrow H\longrightarrow M\longrightarrow
L_2(b)\longrightarrow0.
\]
Thus $L_1(a)\cong M/H\cong L_2(b)$. By
Remark~\ref{rem:ideals-up-to-shift}, $L_1=L_2$. Joining the two
even linkage chains proves that $I'$ and $I''$ are evenly
homogeneously CI-linked.
\end{proof}

Now the rest of the proof for Rao's theorem is immediate.
\begin{proof} (\ref{rao}) ($1\implies 2$) Twist the second N-type resolution by $t$, then add $P$
to the kernel and middle module of the first sequence, and $Q$ to those of the second.
The isomorphism $M'\oplus P\cong M''(t)\oplus Q$ identifies their middle modules.
Lemma~\ref{lem:compare} therefore shows that $I'$ and $I''$ are evenly
homogeneously CI-linked.    
\end{proof}

The following corollary classifies the complete intersection linkage classes of homogeneous height two Cohen-Macaulay ideals by using their N-type modules and their syzygies.
\begin{cor}\cite[Corollary~2.4, 2.5]{herzog1987maximal}\label{oddCIlinkageclass}
Let \(R\) be the same as in Theorem~\ref{rao} and let \(I'\) and \(I''\) be proper homogeneous Cohen-Macaulay ideals of height \(2\), with graded N-type resolutions
\[
0\longrightarrow F'\longrightarrow M'\longrightarrow I'(h_1)\longrightarrow0
\qquad\text{and}\qquad
0\longrightarrow F''\longrightarrow M''\longrightarrow I''(h_2)\longrightarrow0.
\]
Then the following hold.
\begin{enumerate}
\item[(a)] \(I'\) and \(I''\) are evenly homogeneously complete intersection linked if and only if $
\Omega(I')$ and $\Omega(I'')$ 
are stably isomorphic up to shift as graded \(R\)-modules.

\item[(b)] \(I'\) and \(I''\) belong to the same odd homogeneous complete intersection linkage class if and only if $M'$ and $\Omega(M'')^*$ are stably isomorphic up to shift as graded \(R\)-modules.
\end{enumerate}
\end{cor}

The same proof given in \cite{herzog1987maximal} works for the graded setting, too.

\section{Doubly layered modules}

For the rest of this section, we stick to the same notations, i.e., $R=R_0\oplus\bigoplus_{i>0} R_i$ be a graded-local Gorenstein ring, with $R_0$ local ring with an infinite residue field and $\gr R_nR\ge 2$, for all $n\gg 0$.

The main result of this section is to prove the result of \cite[Corollary~5.3]{CasanellasDrozdHartshorne+2005+149+171} for Gorenstein rings of the form above. We will follow an algebraic approach to the proof, but for the most part, the proof is almost the same. When we use the corresponding local rings results, the same proof as in the graded case will prove the result in the local case. Firstly, we observe the existence of affine graded basic elements. The proof of this is very close to Lemma~\ref{lem:bourbaki-rank-reduction}.

\begin{lem}[Affine basic element lemma]\label{affinebasiclemma}
Let $R$ be a graded-local ring as in section 2 and $M$ be a torsion-free $\mathbb Z$-graded
$R$-module, locally free of constant rank $r>1$ in depth one. Fix $x\in M_d$,
and let $V\subseteq M_d$ be an $R_0$-submodule such that
$(RV)_{\mathfrak p}=M_{\mathfrak p}$ whenever $\depth R_{\mathfrak p}\leq 1$.
Then there exists $v\in V$ such that $x+v$ is basic in depth one.
\end{lem}
\begin{proof}[Sketch of proof]
For $R_0$-generators $u_1,\ldots,u_q$ of $V$, apply Remark~\ref{rem:homogeneous-choice} to $x+cu_1+\cdots+c^qu_q$ to obtain $z=x+v_0$ basic at $\Ass R$.
Apply the same argument to the images of $V$ in $M/Rz$ to choose $y\in V$, as in Lemma~\ref{lem:bourbaki-rank-reduction}.
The Fitting-ideal argument there then gives $c\in R_0$ such that $z+cy=x+(v_0+cy)$ is basic in depth one.
\end{proof}


\begin{prop}\cite[Proposition~3.4]{CasanellasDrozdHartshorne+2005+149+171}\label{lemma3.4CDH} Let $R$ be a graded-local Gorenstein ring, as above. Let $I$ be height 2 Cohen-Macaulay homogeneous proper ideal of $R$, with an N-type resolution $$ 0\x L\x N\x I(h)\x 0.$$ Suppose there is a graded exact sequence $$0\x E\xrightarrow{\iota} N\xrightarrow{\pi} N'\x0$$ where $E$ is rank 2 MCM orientable graded $R$-module, and $N'$ MCM orientable graded $R$-module with $\rank N'\ge 2$. Then there is an height 2 Cohen-Macaulay proper homogeneous ideal $J$ of $R$, in the same even homogeneous Gorenstein linkage class as $I$, with an N-type resolution $$0\x L'\x N'\x J(d)\x 0$$ for some $d\in\mathbb{Z}$.    
\end{prop} 

\begin{proof}
Set \(\rho=\rank N'\), by Remark~\ref{existenseofbasicelements}, choose \(b\gg0\), the degree \(b\) pieces of \(E,N,N'\) generate these modules locally in depth one.
By the Theorem \ref{gradedbasicelements}, we get an N-type resolution 
\[
0\longrightarrow L'\longrightarrow N'\longrightarrow J(d)\longrightarrow0,
\qquad L'=R(-b)^{\rho-1}.
\]
We claim that \(J\) and \(I\) are in the same even homogeneous Gorenstein linkage class. Choose homogeneous elements \(x_1,x_2\in J\) whose images in \(J(d)\) have degree \(b\), and which form an \(R\)-regular sequence and \(\mathfrak c=(x_1,x_2)\subsetneq J\). Let \(B=R(-b)^2\), 
\(A=R(-d-2b)\), and let
\[
0\longrightarrow A\longrightarrow B\longrightarrow \mathfrak c(d)\longrightarrow0
\]
be the Koszul resolution of \(\mathfrak c(d)\). Set \(J_1=\mathfrak c:J\); then \(J\) and \(J_1\) are directly homogeneously CI-linked by \(\mathfrak c\). Choose homogeneous lifts of the images of \(x_1,x_2\) in \(J(d)\) to \(N'\). Comparing the Koszul resolution of \(\mathfrak c(d)\) with the given N-type resolution of \(J\), and using the mapping cone, we get an exact sequence
\[
0\longrightarrow A\longrightarrow L'\oplus B\longrightarrow N'\longrightarrow (J/\mathfrak c)(d)\longrightarrow0.\tag{1}
\]
Set \(P=L'\oplus B=R(-b)^{\rho+1}\). Now lift the map \(P\to N'\) in (1) to a homogeneous map \(\alpha:P\to N\), using the exact sequence \(0\to E\to N\to N'\to0\). Since the composite \(A\to P\to N'\) is zero, the restricted map \(A\to P\xrightarrow{\alpha}N\) factors through \(E\). Denote this map by \(\beta:A\to E\), i.e., we have the following commutative diagram
\[
\begin{tikzcd}[row sep=small, column sep=small,
              font=\small, cramped]
  & 0 & 0 \\
0 \arrow[r]
  & \mathcal F \arrow[u] \arrow[r]
  & N' \arrow[u] \arrow[r]
  & (J/c)(d) \arrow[r] & 0 \\
  & P \arrow[u] \arrow[r,"\alpha"]
  & N \arrow[u] \\
  & A \arrow[u] \arrow[r,"\beta"]
  & E \arrow[u] \\
  & 0 \arrow[u]
  & 0 \arrow[u]
\end{tikzcd}
\]
Here \(\mathcal F\) is just the kernel of the map \(N'\to (J/\mathfrak c)(d)\) the top row and left column are from (1), and the middle column is the one from the hypothesis.

We next perturb \(\alpha\) and \(\beta\); Since the cokernel \(\mathcal{F}=\operatorname{coker}(A\to P)\) is locally free at every homogeneous prime \(\mathfrak p\) with \(\depth R_\mathfrak p\leq1\); the map $
\operatorname{Hom}_R(P,E)\longrightarrow \operatorname{Hom}_R(A,E)$ 
is surjective after localizing at every such prime. 

Set $H=\h_R(A,E)\cong E(d+2b)$ and
$W=\img\big(\h_R(P,E)_0\to H_0\big)$. By the choice of $b$, we have
$(R\h_R(P,E)_0)_{\mathfrak p}=\h_R(P,E)_{\mathfrak p}$; hence, we get $(RW)_{\mathfrak p}=H_{\mathfrak p}$ whenever $\dep R_{\mathfrak p}\leq1$.
Now $H$ is torsion-free and locally free of rank two in depth one, and $\beta\in H_0$; the affine basic element lemma (Lemma~\ref{affinebasiclemma})

gives $s\in W$ such that $\beta+s$ is basic in depth one. Hence by Corollary~\ref{rank2Gorenstein}, we can find a height 2 homogeneous Gorenstein ideal $K$, such that $0\x A\xrightarrow{\beta+s}E\x K(c)\x 0$ is exact for some $c\in\mathbb Z$. 

Since $s\in W$, there exists a homogeneous map $t:P\to E$ of degree zero such that
$t|_A=s$. Set $\beta'=\beta+s$ and $\alpha'=\alpha+\iota t$. Then $\alpha'|_A=\iota\beta'$ and
$\pi\alpha'=\pi\alpha$, so the new diagram remains commutative. By the snakes lemma, we get the following diagram

\[
\begin{tikzcd}[row sep=1.1em, column sep=2em, cramped]
  & 0 & 0 & 0 \\
0 \arrow[r]
  & \mathcal F \arrow[u] \arrow[r]
  & N' \arrow[u] \arrow[r]
  & (J/c)(d) \arrow[u] \arrow[r] & 0 \\
0 \arrow[r]
  & R(-b)^{\rho+1} \arrow[u] \arrow[r,"\alpha+\iota t"]
  & N \arrow[u] \arrow[r]
  & G \arrow[u] \arrow[r] & 0 \\
0 \arrow[r]
  & R(-d-2b) \arrow[u] \arrow[r,"\beta+s"]
  & E \arrow[u] \arrow[r]
  & K(c) \arrow[u] \arrow[r] & 0 \\
  & 0 \arrow[u]
  & 0 \arrow[u]
  & 0 \arrow[u]
\end{tikzcd}
\]
Here \(G=\operatorname{coker}(\alpha+\iota t)=\operatorname{coker}(\alpha')\). We claim that \(G\), up to a twist, is a proper homogeneous height two Cohen-Macaulay ideal;
Let $n=\dim R$, the right-hand column shows that $G_{\mathfrak p}\cong
R_{\mathfrak p}(c)$ whenever $\depth R_{\mathfrak p}\leq1$, and using depth lemma on
$0\to P\to N\to G\to0$ we get $G$ is
torsionless of rank one and $\depth G\geq n-1$. Since $K(c)\to G$ is an
isomorphism in depth one, $G^{**}\cong K(c)^{**}\cong R(c)$; hence
$G\cong I_1(c)$ for a homogeneous ideal $I_1$ with $\height I_1\geq2$.
By the choice of $b$, $\alpha'(P)\subseteq\mathfrak mN$, so
$\mu(I_1)=\mu(G)=\mu(N)>1$ and therefore $I_1$ is proper. Finally,
$\depth I_1\geq n-1$ gives $\depth R/I_1\geq n-2$; since
$\height I_1\geq2$, it follows that $\dim R/I_1=n-2$, so $I_1$ has height
two and is Cohen-Macaulay.

Now $I_1(c)$ and $I(h)$, have the same N-type modules, by Theorem \ref{rao}, we get $I$ and $I_1$ to be in the same even homogeneous complete intersection linkage class. Now as $K\subsetneq I_1$, set $I_2=K:I_1$, then $\omega_{R/I_2}\cong (I_1/K)(c)\cong (J/\mathfrak c)(d)$. We get $\Ann_R(I_1/K)=\Ann_R(J/\mathfrak c)$, so $K:I_1=\mathfrak{c}:J$
\end{proof}

\begin{remark}\label{linksin3.4CDH}
    It is important to note the links here $I\sim\cdots\sim I_1\sim_G J_1\sim J$, where $\sim_G$ is a Gorenstein link that is not necessarily a complete intersection, while the rest of the links are complete intersections. Also, note that the linking Gorenstein ideal has an N-type resolution $0\x R\x E\x K\x 0$. So the linking Gorenstein ideal $\sim_G$ is a CI if and only if $E$ is free.
\end{remark}
This is the only lemma that needs to be noted to prove Corollary 5.3 of \cite{CasanellasDrozdHartshorne+2005+149+171}, for graded-local Gorenstein rings in our setting, or local Gorenstein rings with an infinite residue field. The same proofs from their paper hold in the general setting too; for the other statements, there is no mention of normality or the domain being used in their proofs. 

\begin{defi}\cite[Definition~4.4]{CasanellasDrozdHartshorne+2005+149+171}\label{doublylayereddefinition} Let $R$ be a Gorenstein local (graded) ring. A $R$-module $M$ is called \emph{doubly layered} if there exists a finite filtration of submodules $$0=M_0\seq M_1\seq\cdots\seq M_n=M$$ with the following properties $(a),(b)$
\begin{itemize}
    \item[$(a)$] Each factor $M_{i}/M_{i-1}$ is either:
    \begin{itemize}
        \item[(i)] a rank 2 (graded) MCM orientable module that is not free, or
        \item[(ii)] a module isomorphic to $(\Omega(F))^*$, where $F$ is of type~(i).
    \end{itemize}
    \item[(b)] There exists an integer $r$, with $0\le r\le n$, such that for $i\le r$, $M_{i}/M_{i-1}$ is of type~(i) and for all $i>r$, $M_{i}/M_{i-1}$ is of type~(ii). 
\end{itemize}    
\end{defi}

The rank 2 terms in the definition are motivated by Corollary~\ref{rank2Gorenstein}, and the syzygy dual of rank 2 modules is due to Corollary~\ref{oddCIlinkageclass}. We call rank 2 MCM (graded), orientable modules that are not free as \emph{type~(i)} modules. And Modules of the form $\Omega(F)^*$, where $F$ is of type~(i) as \emph{type~(ii)} modules. 

\begin{remark}

    Let $0\xrightarrow{} M'\xrightarrow{\alpha} M\xrightarrow{\beta} M''\x 0$ be an exact sequence, where $M'$ and $M''$ are doubly layered modules such that all the layers (successive quotient's in its filtration) of $M'$ are of type (i) or all layers of $M''$ are of type (ii). Then $M$ is Doubly layered.

\end{remark} 
\begin{proof} Let $0=M'_0\seq\cdots\seq M'_k=M'$ and $0=M''_0\seq\cdots\seq M''_l=M''$ be the doubly layered filtration of $M'$ and $M''$, which satisfy the assumptions of the remark. Then $$0=\alpha(M'_0)\seq \alpha(M'_1)\seq\cdots\seq \alpha(M'_k)=\beta^{-1}(M''_0)\seq \beta^{-1}(M''_1)\seq\cdots\seq\beta^{-1}(M''_l)=M$$ is a doubly layered filtration of $M$.
\end{proof}

\begin{prop}\cite[Proposition~4.5]{CasanellasDrozdHartshorne+2005+149+171}\label{Omega(M)^*isDL}
    Let $R$ be a local (graded-local) Gorenstein ring. Let $M$ be MCM (graded) orientable module, if $M$ is doubly layered, then there exist a (graded) free module $F$, such that $\Omega(M)^*\oplus F$ is doubly layered.

In particular, if $M$ admits a doubly layered filtration whose successive quotients are all of type~(i) (resp. type~(ii)), then there exists a (graded) free
module $F$ such that $\Omega_R(M)^*\oplus F$ admits
a filtration whose successive quotients are all
of type~(ii) (resp. type~(i)).
\end{prop}

\begin{thm}\label{DLstablyGlicci}\cite[Corollary~5.3]{CasanellasDrozdHartshorne+2005+149+171}
    Let $R$ be a Gorenstein local ring with an infinite residue field, and $\dim R\ge 2$ or $R=R_0\oplus\bigoplus_{i>0} R_i$ graded-local Gorenstein ring, with $R_0$ is a local ring with an infinite residue field, and $\gr R_nR\ge 2$ for all $n\gg 0$. Let $I$ be a proper (homogeneous) height two Cohen-Macaulay ideal of $R$. Then $I$ is Glicci if and only if for some (equivalently for any) N-type resolution $0\x F\x M\x I(d)\x 0$, $M$ is stably isomorphic (up to shift) to a doubly layered module.  
\end{thm}

\begin{cor}\label{Gorensteinidealsareglicci} Let $R$ be as in Theorem~\ref{DLstablyGlicci}. Let $K$ be a homogeneous Gorenstein ideal of $\height 2$. Then $K$ is homogeneously glicci.   

 If $K$ is homogeneous Gorenstein ideal of $\height g\ge 3$. If $\height K_+=\height K$, then $K$ is glicci. 
\end{cor}
\begin{proof}
For height two Gorenstein ideals, it follows from Theorem~\ref{DLstablyGlicci}. This also gives 2-glicci. If $\height K_+=\height K$, then by Remark~\ref{existenceofgorensteinsubideal} there exists a homogeneous complete intersection $(f_1,...,f_{g-2})\seq K$ such that $K/(f_1,...,f_{g-2})$ is glicci in $R/(f_1,...,f_{g-2})$. Hence $K$ is glicci in $R$. 
\end{proof}

Using Theorem~\ref{DLstablyGlicci}, and the following structure theorem given by Buchweitz, Eisenbud, Herzog, and Knorrer, we get the following. 

\begin{thm}\cite{10.1007/BFb0078838}\label{thm:graded-mcm-quadrics}
Let \(k\) be an algebraically closed field of arbitrary characteristic and let 
\[
A=k[z_1,\ldots,z_r]/(q),
\qquad 3\leq r\leq6,
\]
where \(q\) is a quadratic form such that \(\operatorname{Proj}A\) is
nonsingular. Then, after a linear change of coordinates,
$q$ has one of the standard forms below, and up to graded isomorphism and degree shift, the
indecomposable graded MCM \(A\)-modules are precisely
\[
\begin{array}{c|c|l|c}
r & \text{Standard form of }q
  & \text{Indecomposable graded MCM modules}
  & \text{Ranks}\\ \hline
3 & xy+z^2
  & A,\ L=(x,z) & 1,1\\
4 & xy+zw
  & A,\ L_+=(x,z),\ L_-=(x,w) & 1,1,1\\
5 & xy+zw+t^2
  & A,\ E=\operatorname{coker}\Phi_5 & 1,2\\
6 & xy+zw+uv
  & A,\ E_+=\operatorname{coker}\Phi_6,\
       E_-=\operatorname{coker}\Psi_6 & 1,2,2
\end{array}
\]

\noindent The modules $E,E_+,E_-$ are  orientable MCM modules, defined by the
following matrix factorizations:
\[
\Phi_5=
\begin{pmatrix}
t&x&z&0\\
y&-t&0&z\\
w&0&-t&-x\\
0&w&-y&t
\end{pmatrix}.
\]
And for $r=6$,
\[
P=\begin{pmatrix}x&z\\-w&y\end{pmatrix},
\qquad
Q=\begin{pmatrix}y&-z\\w&x\end{pmatrix},\text{ and}\qquad
\Phi_6=
\begin{pmatrix}P&uI_2\\vI_2&-Q\end{pmatrix},
\qquad
\Psi_6=
\begin{pmatrix}Q&uI_2\\vI_2&-P\end{pmatrix}.
\]
\end{thm}

Consequently, every finitely generated graded MCM \(A\)-module is isomorphic, as a graded module, to a finite direct sum, up to degree shifts, of the modules appearing in the corresponding row. 

\begin{prop}\label{Quadratic hypersurfacesingularities}
    Let $R=k[x_1,..,x_n]/(f)$, where $3\le n\le 6$, $k$ is an algebraically closed field, $(f)$ is a degree 2 homogeneous irreducible polynomial, such that $\Proj R$ is nonsingular. Then every proper homogeneous Cohen-Macaulay ideal of height two is glicci. 
\end{prop}
\begin{proof}
    By Theorem~\ref{DLstablyGlicci} and using the theorem above, we obtain the required result.
\end{proof}

Note that, though $L$, $L_+$, and $L_-$ are not orientable, we still have $L\oplus L$, and $L_+\oplus L_-$ is orientable, and this is enough for the proof. 
 We observe the following lemmas on doubly layered modules, which will be used in later sections.

\begin{lem}\label{shortexactDL} Let $R$ be a graded-local Gorenstein ring.
Let $0\x N\x M\x E\x 0$ be a graded exact sequence of MCM $R$-modules, such that $N\oplus F_1\cong D_1(t)\oplus G_1$, and $E\oplus F_2\cong D_2(e)\oplus G_2$, where $F_1,F_2,G_1$ and $G_2$ are graded free modules, $D_1$ and $D_2$ are doubly layered modules. If either 
\begin{itemize}
    \item[$(a)$] $D_1$ has a finite filtration, in which all the layers are of type (i) or,
    \item[$(b)$] $D_1$ has a finite filtration, in which all the layers are type (ii), and $D_2$ has a finite filtration, in which all the layers are type (ii).
\end{itemize}

Then $M$ is stably isomorphic to a doubly layered module up to shift, i.e., there is graded free modules, $F,G$ and a graded doubly layered module $D$ such that $M\oplus F\cong D(c)\oplus G$, where $c\in\mathbb Z$.
\end{lem}
\begin{proof}
Take the direct sum of the given exact sequence and $0\x F_1\x F_1\x0\x0$ to get an exact sequence $0\x N\oplus F_1\x M\oplus F_1\x E\x 0$. Replace $N\oplus F_1$ with $D_1(t)\oplus G_1$ to get $0\x D_1(t)\oplus G_1\x M\oplus F_1\x E\x0$. Now by Lemma~\ref{split} we get an exact sequence $0\x D_1(t)\x M'\x E\x 0$, where $M'\oplus G_1\cong M\oplus F_1$.

Now take the direct sum of $0\x D_1(t)\x M'\x E\x 0$ and $0\x 0\x F_2\x F_2\x0$ to get an exact sequence $0\x D_1(t)\x M'\oplus F_2\x E\oplus F_2\x0$. Replace $E\oplus F_2$ with $D_2(e)\oplus G_2$ to get $0\x D_1(t)\x M'\oplus F_2\x D_2(e)\oplus G_2\x0$. Now consider the commutative diagram 
\[\begin{tikzcd}[row sep=small, column sep=small,
              font=\small,cramped]
	& 0 & 0 & 0 \\
	0 & D_1(t) & {M''} & D_2(e) & 0 \\
	0 & D_1(t) & {M'\oplus F_2} & {D_2(e)\oplus G_2} & 0 \\
	0 & 0 & {G_2} & {G_2} & 0 \\
	& 0 & 0 & 0 
	\arrow[from=1-2, to=2-2]
	\arrow[from=1-3, to=2-3]
	\arrow[from=1-4, to=2-4]
	\arrow[from=2-1, to=2-2]
	\arrow[from=2-2, to=2-3]
	\arrow[from=2-2, to=3-2]
	\arrow[from=2-3, to=2-4]
	\arrow[from=2-3, to=3-3]
	\arrow[from=2-4, to=2-5]
	\arrow[from=2-4, to=3-4]
	\arrow[from=3-1, to=3-2]
	\arrow[from=3-2, to=3-3]
	\arrow[from=3-2, to=4-2]
	\arrow["\psi", from=3-3, to=3-4]
	\arrow["{\pi\circ\psi}", from=3-3, to=4-3]
	\arrow[from=3-4, to=3-5]
	\arrow["\pi", from=3-4, to=4-4]
	\arrow[from=4-1, to=4-2]
	\arrow[from=4-2, to=4-3]
	\arrow[from=4-2, to=5-2]
	\arrow["id", from=4-3, to=4-4]
	\arrow[from=4-3, to=5-3]
	\arrow[from=4-4, to=4-5]
	\arrow[from=4-4, to=5-4]
\end{tikzcd}\] here, the exactness of the top row is due to the snake lemma. Note that $M''\oplus G_2\cong M'\oplus F_2$. So consider $0\x D_1(t)\x M''\x D_2(e)\x 0$, and lifting the filtrations of $D_1$ and $D_2$, we get a filtration for $M''$. The assumptions $(a)$ or $(b)$ take care of the condition $(b)$ in  Definition~\ref{doublylayereddefinition} for this filtration. Hence, $M''$ is doubly layered; this proves $M$ to be stably isomorphic to a doubly layered module. 
\end{proof}

\section{Ideals with clean Gorenstein filtration}

Let $R$ be a commutative graded ring. We say a graded $R$-module $M$ is \emph{clean} if $M$ has a graded finite filtration $0=M_0\seq M_1\seq\cdots \seq M_{n-1}\seq M_n=M$, such that $M_i/M_{i-1}\cong R/P_i(d_i)$, with $P_i\in\Ass(M)$ and $d_i\in\mathbb Z$. We say a homogeneous ideal $I$ of $R$ is clean if $R/I$ is a clean $R$-module. 

\begin{defi}\label{deficleangorenstein}
    Let $R$ be a graded ring, and let $I\subsetneq R$ be a homogeneous ideal. We say $I$ has a \emph{clean Gorenstein filtration} if all the associated primes of $R/I$ are Gorenstein ideals, and there exists a finite filtration $$0=M_0\seq M_1\seq\cdots\seq M_n=R/I,\qquad M_{i}/M_{i-1}\cong R/P_i(d_i)$$ where $P_i\in\Ass(R/I)$ and $d_i\in\mathbb Z$.
\end{defi}

If $(R,\mathfrak m)$ is a graded-local ring, then $\mathfrak m$-primary homogeneous ideals have clean Gorenstein filtration. Let $S$ be a polynomial ring, then clean Monomial ideals has a clean Gorenstein filtration, as all the associated primes of monomial ideals are generated by some collection of variables.

\begin{thm}\label{cleanGorensteinfiltrations}
    Let $(R=R_0\oplus\bigoplus_{i>0}R_i,\mathfrak m=\mathfrak m_0\oplus\bigoplus_{i>0} R_i,k)$ be a graded-local Gorenstein ring, with $\gr R_nR\ge 2$ for all $n\gg 0$, $(R_0,\mathfrak m_0)$ a local ring and $|k|=\infty$. Let $I$ be a proper homogeneous Cohen-Macaulay ideal of height two. If $I$ has a clean Gorenstein filtration, then $I$ is homogeneously glicci.   
\end{thm}

First, we note the following key lemma which is used to prove the theorem. 

\begin{lem}\label{IglicciiffOmega(I)^*}
    Let $R$ be as in Theorem~\ref{cleanGorensteinfiltrations}, and let $I\subsetneq R$ be a homogeneous Cohen-Macaulay ideal of height two. Then $I$ is homogeneously glicci if and only if $\Omega_R(I)^*=\Omega^2_R(R/I)^*$ is stably isomorphic to a doubly layered module.  
\end{lem}

\begin{proof}
    Consider an N-type resolution of $I$ [Proposition~\ref{existenceNtype}] $0\x F\x M\x I\x 0$. Notice that $\Omega_R(I)^*$ is an orientable MCM $R$-module as the height of $I$ is two. So, consider an $N$-type resolution of it  $0\x G\x\Omega(I)^*\x J\x 0$. From the first exact sequence, notice that $\Omega(M)^*$ is stably isomorphic to $\Omega(I)^*$. Hence, by Corollary~\ref{oddCIlinkageclass}, we get that $I$ and $J$ are complete intersection linked by an odd number of links. Hence, $I$ is glicci $\iff$ $J$ is glicci $\iff$ (by Theorem~\ref{DLstablyGlicci}) $\Omega_R(I)^*$ is stably isomorphic up to shift to a doubly layered module. 
\end{proof}

\begin{proof}[Proof of Theorem~\ref{cleanGorensteinfiltrations}] By the previous lemma, it is enough to show $\Omega^2(R/I)^*$ is stably isomorphic to a doubly layered module. As $I$ is Cohen-Macaulay, all its associated primes are height 2 Gorenstein ideals. Consider the clean Gorenstein filtration $$0=M_0\seq M_1\seq\cdots\seq M_n=M$$ such that $M_{i}/M_{i-1}\cong R/P_i(d_i)$ and $P_i\in\Ass(R/I)$. Consider the short exact sequence $$0\x M_{i}\x M_{i+1}\x R/P_i(d_i)\x 0$$ and take the dual of the second syzygies to the exact sequence above. Hence we get the following $$0\x\Omega^2(R/P_i(d))^*\x \Omega^2(M_{i+1})^*\oplus F\x\Omega^2(M_i)^*\x 0$$ where $F$ is a graded free module. As $P_i$'s are Gorenstein ideals of height two, consider an N-type resolution $0\x R(\beta_i)\x E_i\x P_i\x 0$ where $\rank E_i=2$[Corollary~\ref{rank2Gorenstein}]. Then $\Omega_R(P_i)^*=\Omega^2_R(R/P_i)^*$ is stably isomorphic up to shift to $\Omega(E)^*$. They $\Omega^2(R/\overline{P_i}(d))^*$ is stably isomorphic up to shift to a type~(ii) module.  Hence by induction on $i$, and using Lemma~\ref{shortexactDL} we get that $\Omega^2(M_n)^*=\Omega^2(R/I)^*$ is stably isomorphic up to shift to a doubly layered module.     
\end{proof}

\begin{defi}
    Let $\mathcal P_{\mathrm{i}}(R)$ (respectively,
$\mathcal P_{\mathrm{ii}}(R)$) denote the class of (graded) orientable
MCM $R$-modules that are stably isomorphic, up to grading
shift, to a (graded) orientable doubly layered module $D$ admitting a finite filtration
\[
0=D_0\subseteq D_1\subseteq\cdots\subseteq D_s=D,
\]
whose successive quotients are all of type (i)
(respectively, type (ii)). We include the zero module
and all (graded) free modules in both classes.
\end{defi}

\begin{thm}\label{2dimDL}
Let $(R,\mathfrak m, k)$ be a two-dimensional local Gorenstein ring with an infinite residue field, or $(R=R_0\oplus\bigoplus_{i>0}R_i,\mathfrak m=\mathfrak m_0\oplus\bigoplus_{i>0}R_i,k)$ be a two-dimensional graded-local Gorenstein ring with an infinite residue field, $(R,\mathfrak m_0)$ a local ring, and $\gr R_nR\ge 2$ for all $n\gg 0$. Let $N$ be a ($\mathbb Z$ graded) orientable MCM $R$-module.
Then $N\in\mathcal{P}_i(R)\cap \mathcal{P}_{ii}(R)$, that is; N  is stably isomorphic (up to grading shift)
to (graded) doubly layered modules $D_1$ and $D_2$,
where $D_1$ admits a filtration whose successive
quotients are all of type~(i), and $D_2$ admits
a filtration whose successive quotients are all
of type~(ii). 
\end{thm}

\begin{proof}
We may assume that $N$ is not free. Set
$T(M)=\Omega_R(M)^*$ for any MCM $R$-module $M$.
We first show that $T(N)\in\mathcal P_{\mathrm{ii}}(R)$. 

Since $N$ is orientable, it admits an $N$-type resolution
\[
0\longrightarrow F\longrightarrow N
\longrightarrow I(a)\longrightarrow0,
\]
where $I$ is an $\mathfrak m$-primary ideal. Taking syzygies and dualizing, we obtain $
T(N)\cong_{\mathrm{st}}
\bigl(\Omega_R^2(R/I)\bigr)^*(-a)$. Since $R/I$ has finite length, the proof of
Theorem~\ref{cleanGorensteinfiltrations} shows that $\bigl(\Omega_R^2(R/I)\bigr)^*
\in\mathcal P_{\mathrm{ii}}(R)$.
Hence, $T(N)\in\mathcal P_{\mathrm{ii}}(R)$.

Now, note that for every MCM $R$-module $M$, $T^2(M)\cong_{\mathrm{st}}M$. 
So, applying the first part of the proof
to $T(N)$, we get that
\[
N\cong_{\mathrm{st}}T^2(N)\in\mathcal P_{\mathrm{ii}}(R).
\]

Finally, Proposition~\ref{Omega(M)^*isDL} and the horseshoe lemma
show that $T\bigl(\mathcal P_{\mathrm{ii}}(R)\bigr)
\subseteq\mathcal P_{\mathrm{i}}(R)$. 
Since $T(N)\in\mathcal P_{\mathrm{ii}}(R)$, we get
$T^2(N)\in\mathcal P_{\mathrm{i}}(R)$, and hence
$N\in\mathcal P_{\mathrm{i}}(R)$.
\end{proof}

In the local case, we use Definition~\ref{cleanGorensteinfiltrations} with the grading shifts omitted.

\begin{thm}\label{thm:clean-gorenstein-local}
Let $(R,\mathfrak m,k)$ be a local Gorenstein ring with an infinite residue
field. Let $I\subsetneq R$ be a Cohen-Macaulay ideal of height $g\geq 2$
admitting a clean Gorenstein filtration. Then $I$ is $2$-glicci. More
precisely, for every Gorenstein ideal $K\subseteq I$ of height at most
$g-2$, the ideal $I/K$ is glicci in $R/K$. In particular, $I$ is glicci.
\end{thm}

\begin{proof}
\emph{Claim.} It suffices to consider Gorenstein subideals of height $g-2$.
Indeed, let $h=\operatorname{ht}_R K<g-2$, and write $A=R/K$ and $J=I/K$.
Since $A$ is Cohen--Macaulay and $\operatorname{ht}_A J=g-h$, choose a
$A$-regular sequence of length $g-h-2$ in $J$. The inverse image $K'$
of the ideal it generates is a Gorenstein subideal of $I$ of height $g-2$.
If $I/K'$ is glicci in $R/K'$, Remark~\ref{rmk:heightcutglicci} applied to $K'/K\subseteq J$
shows that $J$ is glicci in $A$. This proves the claim.

Now suppose that $\operatorname{ht}_R K=g-2$. Then $J$ has height two
in $A$, and $A/J\cong R/I$ is Cohen-Macaulay. Since $K$ annihilates
$R/I$, its clean Gorenstein filtration is also a clean Gorenstein filtration of $A/J$ over $A$. Thus, $J$ has a clean
Gorenstein $A$-filtration. 

Hence, by the local version of Theorem~\ref{cleanGorensteinfiltrations}, we conclude that $J$ is glicci; Hence, $I$ is $2$-glicci.
Taking $K=(0)$ gives the last assertion.
\end{proof}

\begin{remark}\label{rem:clean-gorenstein-grading}
In the graded case, a reduction to height two must also preserve the
grade conditions of Theorem~\ref{cleanGorensteinfiltrations}. The example of weighted polynomial ring in
Remark~\ref{existenceofgorensteinsubideal}
shows that this condition can fail after quotienting by a
homogeneous Gorenstein subideal, even when $R_0$ is a field. If $R$ is
standard graded over an Artinian local ring, every homogeneous Gorenstein
quotient $A$ of dimension at least two satisfies the grade conditions,
since $\sqrt{A_nA}=\mathfrak m_A$ for every $n\geq 1$.
\end{remark}

\begin{thm}\label{thm:clean-gorenstein-graded}
Let \(R=\bigoplus_{n\geq 0}R_n\) be a graded-local Gorenstein ring such
that \(R_0\) is local with infinite residue field and
\(\gr_R(R_nR)\geq 2\) for all \(n\gg 0\). Let \(I\subsetneq R\) be a
homogeneous Cohen-Macaulay ideal of height \(g\geq 2\) admitting a clean
Gorenstein filtration, and set \(I_+=\bigoplus_{n>0}I_n\).

\begin{enumerate}
\item[\textup{(a)}] $I$ is homogeneously 2-glicci i.e., if there exists a homogeneous Gorenstein subideal $K\seq I$, such that $\height I=\height K+2$ and $\gr_{R/K}(R/K)_nR/K\ge 2$ for all $n\gg 0$; then $I/K$ is glicci in $R/K$.

\item[\textup{(b)}] If \(\height_R I_+=g\). Then there exists a
homogeneous complete intersection \(C\subseteq I_+\) of height \(g-2\)
such that \(R/C\) satisfies the same grade condition and \(I/C\) is
homogeneously glicci in \(R/C\). In particular, \(I\) is homogeneously
glicci in \(R\).

\item[\textup{(c)}] If \(R_0\) is Artinian and $R$ is standard graded over $R_0$. Then, for every homogeneous Gorenstein
ideal \(K\subseteq I\) of height at most \(g-2\), the quotient
\(A=R/K\) satisfies $  \gr_A(A_nA)\geq 2$ for all $n\gg0$. 
Hence, by part~\textup{(a)}, we get 
that \(I/K\) is homogeneously glicci in \(R/K\).
Thus \(I\) is homogeneously \(2\)-glicci and in particular $I$ is homogeneously glicci.
\end{enumerate}
\end{thm}

\begin{proof}
For \textup{(a)}, set \(A=R/K\) and \(J=I/K\) where $K$ satisfies the height and grade assumptions. The clean Gorenstein
filtration of \(R/I\) is also a clean Gorenstein filtration of \(A/J\) as $A$-module. Thus \(J\) is homogeneously glicci by
Theorem~\ref{cleanGorensteinfiltrations}.

For \textup{(b)}, we first record the finiteness needed for prime
avoidance. For a homogeneous ideal \(D\subseteq R\), set
\[
   \mathcal M(D)=\bigcup_{n\geq 1}\Min_R\bigl(R/(D+R_nR)\bigr).
\]
We claim that \(\mathcal M(D)\) is finite.

Choose homogeneous \(R_0\)-algebra generators \(x_1,\ldots,x_s\) of
positive degree, and let \(d\) be a common multiple of their degrees.
Every monomial of sufficiently large degree is divisible by some
\(x_i^{d/\deg x_i}\). Thus, writing \(J_n=R_nR\), we have
\(J_{n+d}=J_dJ_n\) for \(n\gg0\). Since \(J_dJ_n\) and \(J_d^2J_n\)
have the same radical, \(D+J_{n+d}\) and \(D+J_{n+2d}\) have the same
minimal primes. Choose \(N\) so that \(D+J_n\) and \(D+J_{n+d}\) have the same minimal
primes for every \(n\geq N\). Every such \(n\) is obtained from one of
\(N,\ldots,N+d-1\) by repeatedly adding \(d\), so no new minimal primes
occur after degree \(N+d-1\). Hence \(\mathcal M(D)\) is finite.

We choose the regular sequence one element at a time. Suppose that \(C_j=(f_1,\ldots,f_j)\) has been chosen and that \(A=R/C_j\) satisfies the grade condition, starting with \(C_0=(0)\). Choose \(N\) so that \(A_nA\) has height at least two for every \(n\geq N\). If a prime of height two contains one of these ideals, it must be minimal over that ideal. By the claim with \(D=C_j\), only finitely many such primes occur. When \(j<g-2\), the image of \(I_+\) in \(A\) has height \(g-j>2\). Homogeneous prime avoidance therefore gives \(f_{j+1}\in I_+\) whose image avoids both these primes and \(\Ass A\). Thus \(\bar f_{j+1}\) is \(A\)-regular, and every prime containing both \(A_nA\) and \(\bar f_{j+1}\) has height at least three for \(n\geq N\). Passing to \(A/(\bar f_{j+1})\) lowers the height of each such prime by one, so \(R/C_{j+1}\) still satisfies the grade condition.

So, \(C=C_{g-2}\) is the required complete intersection. Part~\textup{(a)} gives that \(I/C\) is homogeneously glicci in \(R/C\), and by the Remark~\ref{rmk:heightcutglicci} gives that \(I\) is homogeneously glicci in \(R\).

For \textup{(c)}, by the reduction in the proof of
Theorem~\ref{thm:clean-gorenstein-local}, it suffices to consider
\(\height K=g-2\). For \(A=R/K\), standard gradedness and the Artinian
hypothesis give \(\sqrt{A_nA}=\mathfrak m_A\) for every \(n\geq1\); hence
\(\gr_A(A_nA)=\dim A\geq2\). Part~\textup{(a)} gives that \(I/K\) is
homogeneously glicci, so \(I\) is homogeneously \(2\)-glicci and therefore
homogeneously glicci by Remark~\ref{rmk:heightcutglicci}.
\end{proof}

\subsection{Applications}

Now, we are ready to prove that $\mathfrak m$-primary ideals are glicci. 

\begin{thm}\label{thm:clean-primary-applications}\label{mprimaryglicci}
Let $(R,\mathfrak m,k)$ be a local or graded-local Gorenstein ring of
dimension $d\geq 2$ with infinite residue field, and let $I$ be an
$\mathfrak m$-primary ideal. In the graded case, assume that $I$ is
homogeneous and that $\operatorname{grade}_R(R_nR)\geq 2$ for all $n\gg 0$.
\begin{enumerate}
\item[\textup{(a)}] If $R$ is local, then $I$ is $2$-glicci and hence glicci.
\item[\textup{(b)}] If $R$ is graded-local and $d=2$, then $I$ is homogeneously
glicci.
\item[\textup{(c)}]  Let $R$ is graded local and $d\ge 3$. Suppose if $I_+$ is also $\mathfrak m$-primary. Then there exists a homogeneous complete intersection $C\seq I_+\seq I$ such that $\height I=\height C+2$ and $I/C$ is homogeneously glicci in $R/C$. Hence $I$ is homogeneously glicci in $R$.

\item[\textup{(d)}] Suppose that $R_0$ is Artinian, and $R$ is standard graded over $R_0$. Then $I$ is
homogeneously $2$-glicci; more precisely for any homogeneous Gorenstein subideal $K$ of $I$ with height at most $d-2$, $I/K$ is homogeneously glicci in $R/K$. In particular $I$ homogeneously glicci. 
\end{enumerate}
\end{thm}

\begin{proof}
Since \(I\) is \(\mathfrak m\)-primary, \(\height_R I=d\) and \(R/I\) has finite length. Choose a composition series of \(R/I\), and a graded composition series in the graded case. Its factors are copies of \(k=R/\mathfrak m\), up to shift. Since \(\mathfrak m\) is the unique prime over \(I\) and \(k\) is Gorenstein, this is a clean Gorenstein filtration.

For \textup{(a)}, Theorem~\ref{thm:clean-gorenstein-local} shows that \(I/K\) is glicci in \(R/K\)
for every Gorenstein subideal \(K\subseteq I\) of height at most \(d-2\).
Thus \(I\) is \(2\)-glicci, and taking \(K=(0)\) shows that \(I\) is
glicci. For \textup{(b)}, \(I\) itself has height two. The graded
filtration above and the assumed grade condition allow us to apply
Theorem~\ref{cleanGorensteinfiltrations} directly, so \(I\) is homogeneously glicci.

In \textup{(c)}, \(I_+\) is \(\mathfrak m\)-primary, so
\(\height_R I_+=d=\height_R I\). Theorem~\ref{thm:clean-gorenstein-graded}\textup{(b)} therefore
gives the required complete intersection \(C\) and both glicci
assertions. Part~\textup{(d)} follows from Theorem~\ref{thm:clean-gorenstein-graded}\textup{(c)}.
\end{proof}

Now, the counterexamples in Section~3 and the preceding result show that the glicci property need not, in general, ascend under
deformation over Gorenstein rings. The situation is more subtle when
the ambient rings are regular. Indeed, let $S$ be a regular local ring
with infinite residue field and let $J\subseteq S$ be Cohen-Macaulay.
By successively taking sufficiently general elements that are regular
on $S/J$ and belong to a regular system of parameters of $S$, one may
specialize $J$ to an $\mathfrak m$-primary ideal in a regular local
ring. Since $\mathfrak m$-primary ideals in regular local rings are
glicci, the glicci conjecture is therefore reduced to the following
question.

\begin{question}\label{deformationquestion}
    Let $(S,J)$ be a deformation of $(R,I)$; that is, there exists
a sequence $\alpha_1,\ldots,\alpha_n$ which is regular on both $S$ and
$S/J$, together with an isomorphism
\[
R \cong S/(\alpha_1,\ldots,\alpha_n)\quad\text{under which}\quad 
I \cong
\frac{J+(\alpha_1,\ldots,\alpha_n)}
     {(\alpha_1,\ldots,\alpha_n)}.
\]
Assume that $S$ and $R$ are regular local rings and that $I$ is
$\mathfrak m_R$-primary. If $I$ is glicci, then is $J$ be glicci?
\end{question}

\begin{thm}\label{thm:clean-monomial-application}
Let $S=k[x_1,\ldots,x_n]$ be a standard graded polynomial ring over an
infinite field. Let $I\subsetneq S$ is a clean monomial ideal and $S/I$
is Cohen-Macaulay. If
$\height_S I\geq 2$, then $I$ is homogeneously $2$-glicci i.e., for any homogeneous Gorenstein ideal $K\seq I$ with $\height I\ge \height K+2$, $I/K$ is homogeneously glicci in $S/K$. In particular, $I$ is homogeneously glicci.
\end{thm}

\begin{proof}
Every associated prime $P$ of $S/I$ is generated by a subset of the
variables. Thus $S/P$ is a polynomial ring, and the given clean filtration
of $S/I$ is a clean Gorenstein filtration. Now, Theorem~\ref{thm:clean-gorenstein-graded}\textup{(c)} applies since $S_0=k$. Hence we conclude that $I$ is homogeneously $2$-glicci and glicci. 
\end{proof}

\begin{thm}(A.Dress 1993, Proof in English can be found in \cite[Theorem~2.1.1]{simon1994combinatorial})  Let $I$ be a squarefree monomial ideal of $S=k[x_1,..,x_n]$, and let $I=I_\Delta$, where $\Delta$ is the corresponding simplicial complex.
$S/I$ is a clean $S$-module if and only if $\Delta$ is a shellable simplicial complex. 
\end{thm}

So this extends the work of \cite{nagel2008glicci}, where they show 
weakly vertex decomposable square free monomial ideals are glicci. Now we have extended this result in a stronger setting: Any shellable simplicial complex is 2-Glicci.

\begin{thm}\cite[Corollary~3.4]{popescu2009stanley}\label{Popsescu} Let $S=k[x_1,..,x_n]$ and $I$ be Cohen-Macaulay monomial ideal such that $\dim S/I=2$. Then $I$ is clean.    
\end{thm}

\begin{lem}
Let \(S=K[x_1,\ldots,x_n]\), and let \(I\subset S\) be a monomial ideal such that \(S/I\) is Cohen-Macaulay. If \(\operatorname{ht}I\ge n-2\), equivalently \(\dim S/I\le 2\), then \(S/I\) is clean.
\end{lem}

For $\dim S/I=1$, apply Theorem~\ref{Popsescu} to $(S/I)[y]$; its clean
filtration descends modulo $y$, since $y$ is regular on
every successive quotient. We leave the standard details
to the reader.


\begin{cor}
    Let $S=k[x_1,..,x_n]$, and let $I$ be Cohen-Macaulay monomial ideals, with $\dim S/I\in\{0,1,2\}$. Then $I$ is a 2-glicci ideal of $S$.
\end{cor}

It was shown that Cohen-Macaulay monomial ideals with at most 5 generators are clean; see \cite[Theorem~1.2]{enomoto2019cleanness}. In \cite[Example~2.8]{enomoto2019cleanness} they show a 6 generated height 3 Cohen-Macaulay monomial ideal which is not clean. This example has 60 associated primes.

There exist height-three Cohen--Macaulay monomial ideals
with ten associated primes that are not clean. On the
other hand, the author has shown that every
Cohen-Macaulay monomial ideal with at most seven
associated primes is clean. Since the proof lies
outside the scope of the present paper, we omit it
here. This leads to the following question.

\begin{question}\label{minimalassquestion} Let $S$ be a polynomial ring over a field. Does there exist a Cohen-Macaulay monomial ideal of $S$ with exactly eight or nine associated primes that is not clean? \end{question}

\begin{thm}\label{dimoneclean}
    Let \(k\) be an algebraically closed field and let \(S=k[x_1,\ldots,x_n]\) be standard graded. If \(I\subset S\) is a homogeneous ideal such that \(S/I\) is Cohen-Macaulay of dimension one, then \(S/I\) admits a clean Gorenstein filtration.
\end{thm}
\begin{proof}
We first consider a homogeneous prime \(P\) with \(\dim S/P=1\). Since \(k\) is algebraically closed, there is a point \([a_1:\cdots:a_n]\in\mathbb P^{n-1}_k\) at which every element of \(P\) vanishes. The map \(x_i\mapsto a_it\) induces a surjection \(S/P\to k[t]\), since some \(a_i\ne 0\). This is an isomorphism because both rings are one-dimensional domains. In particular, \(S/P\) is Gorenstein, and the image of any linear form outside \(P\) is a nonzero scalar multiple of \(t\).

Now let \(A=S/I\) and choose a linear nonzerodivisor \(\ell\) on \(A\). Let \(P\) be a minimal prime of \(I\). Since \(P\) is associated to \(A\), we may choose a nonzero homogeneous element \(u\in A\) of least degree among those with annihilator \(P\). Then \(u\notin\ell A\): if \(u=\ell v\), the regularity of \(\ell\) gives \(v\) the same annihilator as \(u\), but \(v\) has a smaller degree, hence a contradiction. Since \(\ell\notin P\), its image in \(S/P\cong k[t]\) is a nonzero multiple of \(t\), so every element of \(S/P\) is a polynomial in \(\bar\ell\). Since \(P\) is the annihilator of \(u\), sending \(1\) to \(u\) identifies \((S/P)(-\deg u)\) with the graded module \(Au=k[\ell]u\).

We next show that \(\ell\) remains a nonzerodivisor on \(A/Au\). Let \(\ell w\in Au\); since \(Au=k[\ell]u\), we may write \(\ell w=(c+\ell q(\ell))u\), where \(c\in k\) is the constant term and \(q(t)\in k[t]\). It follows that \(cu=\ell(w-q(\ell)u)\in\ell A\). If \(c\ne0\), dividing by \(c\) gives \(u\in\ell A\), contrary to our choice of \(u\). Thus \(c=0\), and we have \(\ell w=\ell q(\ell)u\). Since \(\ell\) is a nonzerodivisor on \(A\), we obtain \(w=q(\ell)u\in Au\). Hence, $\ell$ is a nonzero divisor on $A/Au$. If \(A/Au\ne0\), it therefore has depth at least one and, as a quotient of \(A\), dimension at most one. Hence it is again one-dimensional and Cohen-Macaulay.

Begin the filtration of \(A=S/I\) with \(0\subset Au\). If \(Au\ne A\), then \(A/Au\) is again one-dimensional and Cohen-Macaulay. Applying the same construction to \(A/Au\) gives a nonzero cyclic ideal, whose inverse image in \(A\) is the next term of the filtration. At each step, we obtain an ideal \(M_{i+1}\) strictly containing the preceding ideal \(M_i\), with \(M_{i+1}/M_i\) a graded shift of \(S/P_i\) for some homogeneous prime \(P_i\) with \(\dim S/P_i=1\). Since the construction continues whenever \(M_i\ne A\) and \(A\) is Noetherian, the chain reaches \(A\) after finitely many steps. Each \(P_i\) contains \(I\) and is minimal over \(I\), since \(\dim S/P_i=\dim S/I=1\). Finally, \(S/P_i\cong k[t]\) is Gorenstein by the first paragraph, so this is a clean Gorenstein filtration of \(S/I\).
\end{proof}

In fact, the proof tells us the following; for a polynomial ring over an algebraically closed field, Cohen-Macaulay ideals of one dimension are clean. The assumption of the algebraically closed field hypothesis is to ensure that the associated primes are Gorenstein ideals of one dimension. This need not to hold for arbitrary fields.

\begin{cor}\label{algcloseddimone}
    Let $S$ be a polynomial ring, over an algebraically closed field. Let $I$ be a Cohen-Macaulay homogeneous ideal with $\dim S/I=1$. Then $I$ is homogeneously 2-glicci and in particular homogeneously glicci.  
\end{cor}
\begin{proof}
    This follows from Theorem~\ref{dimoneclean} and Theorem~\ref{thm:clean-gorenstein-graded}(c).
\end{proof}

In Section~9 we use a different, longer argument to show that a
homogeneous Cohen--Macaulay ideal \(I\) in a polynomial ring \(S\) over
an arbitrary infinite field is glicci whenever \(\dim(S/I)=1\). Over an
algebraically closed field \(k\), the argument above gives the stronger
conclusion that \(I\) is homogeneously \(2\)-glicci when
\(\operatorname{ht} I\ge 2\). Geometrically, every zero-dimensional
closed subscheme of \(\mathbb{P}^n_k\) is 2-glicci.

\section{Going up one codimension} We retain the standing assumptions on R. We have seen that, when $\dim R\geq 2$, every $\mathfrak m$-primary ideal is glicci, whereas Cohen-Macaulay ideals need not be glicci. The following proposition, generalizing \cite[Lemma 2.1]{migliore2013glicci}, describes the situation after passing to codimension one higher.
\begin{prop}\label{oneheightup}
Let $R$ be a local Gorenstein
ring with an infinite residue field. Let $I\subsetneq R$ be a 
Cohen-Macaulay ideal of height $g\geq 1$ such that $R/I$ is generically
Gorenstein, i.e., $R_{\mathfrak p}/I_{\mathfrak p}$ is Gorenstein for every
$\mathfrak p\in\Ass_R(R/I)$. If $f\in\mathfrak m$ is a  nonzerodivisor on $R/I$, then
$(I,f)$ is glicci.
\end{prop}

\begin{proof}
Since $R/I$ is a Cohen-Macaulay and generically Gorenstein ring, its
canonical module can be identified, up to a shift in the graded case, with a Gorenstein ideal $J/I$ of $R/I$, and $\height_R J=g+1$ \cite[Proposition~3.3.18]{bruns1998cohen}. Note that $I+fJ$ is Gorenstein of height
$g+1$; since $f$ is
$R/I$-regular, $f(J/I)\cong J/I$ is again a canonical ideal of $R/I$.
Hence, $(R/I)/f(J/I)\cong R/(I+fJ)$ is Gorenstein. We claim that $(I,f)$ and $J$ are directly G-linked by $I+fJ$.
Indeed, since $I:f=I$,
\[
(I+fJ):(I,f)=(I+fJ):f=(I:f)+J=J.
\]
Hence, $(I,f)$ is directly G-linked to the Gorenstein ideal $J$. Since
$\operatorname{ht}J=g+1\geq2$, by Corollary~\ref{Gorensteinidealsareglicci}, we conclude that $J$ is glicci.
\end{proof}

This raises the following question. 

\begin{question}\label{questiononeheightup}
    Let $R$ be a graded-local Gorenstein ring. Let $I$ be a homogeneous Cohen-Macaulay ideal of $R$. If $f\in\mathfrak m$ is a homogeneous non zero divisor of $R/I$. Is $(I,f)$ homogeneously glicci in $R$?
\end{question}

Let $R$ be a local Gorenstein ring of dimension at least two. Let $I$ be a Cohen-Macaulay ideal of dimension 1, i.e., $\dim R/I=1$. Then the question has a positive answer. Let $R=k[x_1,..,x_n]$ or $k[[x_1,..,x_n]]$. By Theorem~\ref{mainsection9}, we get that every Cohen-Macaulay ideal of height $n-1$ or of height $n-2$ is (homogeneously) glicci. Hence, over these rings, the above question has a positive answer for ideals of height $n-3$.

\begin{remark}
Let $R$ be a local Gorenstein ring, let $I\subset R$ be a Cohen-Macaulay ideal that is glicci, and let $f$ be $R/I$-regular. Then $(I,f)$ is glicci.
\end{remark}

\begin{proof}
Choose a G-linkage chain $I=I_0\sim_{K_1}I_1\sim\cdots\sim_{K_s}I_s$, with $I_s$ a complete intersection, and let $\Sigma$ be the union of the associated primes of $R/I_i$ and $R/K_i$. Write $\Sigma=\Sigma_0\sqcup\Sigma_1$, where $\Sigma_0=\{P\in\Sigma:f\in P\}$ and $\Sigma_1=\{P\in\Sigma:f\notin P\}$ . If $P\in\Sigma_0$, then $I\nsubseteq P$; hence, $L:=I\cap\bigcap_{Q\in\Sigma_1}Q$ is not contained in any $P\in\Sigma_0$, so prime avoidance gives $a\in L\setminus\bigcup_{P\in\Sigma_0}P$. Thus, $h=f+a$ avoids every prime in $\Sigma$, while $(I,h)=(I,f)$.

Fix $i$ and set $A=R/K_i$, $U=I_i/K_i$, and $V=I_{i+1}/K_i=0:_AU$. Since $A/U$ is MCM over the Gorenstein ring $A$, $\Ext_A^1(A/U,A)=0$; so, $A\to\h_A(U,A)$ is surjective with kernel $V$. If $xU\subseteq hA$, regularity of $h$ on $A$ defines $\varphi\in\h_A(U,A)$ by $xu=h\varphi(u)$; writing $\varphi$ as multiplication by $y\in A$ gives $x-hy\in V$. Hence, $(h):_AU=V+hA$, and similarly $(h):_AV=U+hA$. Since $(K_i,h)$ is Gorenstein, $(I_i,h)\sim_{(K_i,h)}(I_{i+1},h)$ for every $i$; therefore, $(I,f)=(I,h)$ is glicci.
\end{proof}
The preceding argument is more delicate in the graded setting, since the
homogeneous element \(a\) obtained by prime avoidance need not satisfy
\(\deg a=\deg f\); thus \(h=f+a\) need not be homogeneous. Nevertheless,
if \(f\) is general of sufficiently large degree, then \(a\) may be chosen
with \(\deg a=\deg f\), and the argument shows that \((I,f)\) is
homogeneously glicci.

\begin{prop}\label{genericallylicci}
Let $(R,\mathfrak m,k)$ be a Gorenstein local ring with $k$ infinite, and let
$I\subsetneq R$ be a Cohen--Macaulay ideal of height $g\geq1$. Assume that
$I$ is generically licci, i.e.,
\[
I_{\mathfrak p}\ \text{is licci in }R_{\mathfrak p}
\qquad\text{for every }\mathfrak p\in\Ass_R(R/I).
\]
If $f\in\mathfrak m$ is a nonzerodivisor on $R/I$, then $(I,f)$ is glicci.
\end{prop}

\begin{proof}
For a licci $\mathfrak pR_{\mathfrak p}$-primary ideal $Q\subset R_{\mathfrak p}$,
let $\ell(Q)$ denote the minimum number of complete intersection links from $Q$
to a complete intersection, and set $\lambda(I)=\sum_{\mathfrak p\in\Ass(R/I)}\ell(I_{\mathfrak p})$. We induct on $\lambda(I)$; if $\lambda(I)=0$, then $I$ is generically complete intersection, so we get the claim by Proposition~\ref{oneheightup}.

We first show that, if $\lambda(I)>0$, there is a complete-intersection link
$I\sim_C I_1$ such that $I_1$ is generically licci and
$\lambda(I_1)<\lambda(I)$. Write
$\Ass(R/I)=\{\mathfrak p_1,\ldots,\mathfrak p_s\}$. For each $i$ with
$\ell(I_{\mathfrak p_i})>0$, choose the first link
$I_{\mathfrak p_i}\sim_{C_i}Q_i$ in a shortest licci chain; if
$\ell(I_{\mathfrak p_i})=0$, set $C_i=I_{\mathfrak p_i}$. Now by prime avoidance and Nakayama's lemma, we can choose a $C\seq I$, a complete intersection, such that $C_{\mathfrak p_i}=C_i$.

Set $I_1=C:I$. Then
$I\sim_C I_1$, and $I_1$ is Cohen-Macaulay of height $g$. At every associated prime of $I$, we have $(I_1)_{\mathfrak p_i}
=C_i:I_{\mathfrak p_i}=Q_i$, and if 
$\mathfrak q\in\Ass(R/I_1)\setminus\Ass(R/I)$, then
$I_{\mathfrak q}=R_{\mathfrak q}$; hence, $I_1$ is generically licci and
$\lambda(I_1)<\lambda(I)$. Iterating, we obtain complete-intersection links $
I=I_0\sim_{C_0}I_1\sim_{C_1}\cdots\sim_{C_{r-1}}I_r$, such that $\lambda(I_r)=0$. Hence, $I_r$ is a generically complete intersection ideal
and therefore generically Gorenstein.

It remains to carry $f$ through this finite linkage chain. Let $\Sigma$ be the
union of the associated primes of all $R/I_j$ and $R/C_j$. By prime avoidance arguments, similar to those in the proof above, choose a $a\in L\setminus\bigcup_{P\in\Sigma_0}P$. Then $h=f+a$ avoids every prime of $\Sigma$, while $(I,h)=(I,f)$. Hence $h$ is
regular on every $R/I_j$ and $R/C_j$. Again by the same arguments in the previous proof, we get $(I_j,h)\sim_{(C_j,h)}(I_{j+1},h)$. Hence, $(I,f)$ is glicci. 
\end{proof}

Suppose that if $R$ is a normal Gorenstein local ring and $I$ is a height-one
Cohen-Macaulay ideal, or if $R$ is a local Gorenstein ring that satisfies Serre's
condition $(R_2)$ and $I$ is a height-two Cohen-Macaulay ideal, then
$(I,f)$ is glicci for every $R/I$-regular element $f$, since in either
case $I$ is generically licci.

If $\height I\geq 2$, then $I$ is generically glicci. Unlike
the case of complete intersection linkage, however, local Gorenstein
links, in general, cannot be glued to obtain a global Gorenstein link.
This obstruction makes the corresponding global statement considerably
more difficult to establish.

\section{dimension 1 and dimension 2 ideals in polynomial ring}

A classical motivation for the study of complete intersection linkage
is its use in the classification of algebraic curves, particularly in $\mathbb P^3_k$. In this section, we show that any two curves, as well as any two finite set of points in $\mathbb P_k^n$ belong to the same
Gorenstein linkage class. Throughout this section, $k$ denotes an
infinite field. We treat the polynomial
and power series cases simultaneously:
\[
S=k[x_1,\ldots,x_n]
\qquad\text{and}\qquad
S=k[[x_1,\ldots,x_n]].
\]
In the polynomial case, all ideals, modules, and maps are assumed to
be graded, and all links are homogeneous. The dimension of an ideal
$I$ always means the dimension of the quotient ring $S/I$.

\begin{thm}\label{mainsection9}
Let $I\subsetneq S$ and suppose that $S/I$ is Cohen-Macaulay.
Then $I$ is (homogeneously) glicci in either of the following cases:
\begin{enumerate}
\item[$(a)$] $\dim S/I=1$.
\item[$(b)$] $\dim S/I=2$. 
\end{enumerate}
\end{thm}

We break the proof into parts to keep it readable. We prove this by constructing a complete intersection subideal with height two less than that of $I$; such that when going modulo it, $\overline{I}$ is glicci. The glicci-ness in the quotient ring is proved by doubly layered techniques. Then by
Remark~\ref{rmk:heightcutglicci} we get $I$ to be glicci. The proof is rather long and technical because, unlike in the case of a polynomial ring over an algebraically closed field (Corollary~\ref{algcloseddimone}), the associated primes of a one-dimensional ideal need not be Gorenstein.

Let $B$ be a standard graded ring or complete local $k$-algebra with residue field $k$, denote $B\langle t\rangle$ for $B[t]$ in the graded setting and $B[[t]]$ in the complete local setting. Let $\PiC(R)$ and $\PiiC(R)$ denote the stable classes of modules filtered
entirely by factors of type~(i) and type~(ii), respectively. We let both zero module in both the classes. 

\begin{lem}\label{lem:closure}
Let $A$ be a local Gorenstein ring with infinite residue field or $A$ is a graded-local Gorenstein ring, with an infinite residue field:
\begin{enumerate}
\item[$(a)$] Each of $\mathcal{P}_{i}(A)$ and $\mathcal{P}_{ii}(A)$ is closed under degree shifts, finite direct sums, stable equivalence, and extensions.
\item[(b)] $\Omega_A(\Omega_A(L)^*)^*\cong_{st} L$ where $\cong_{st}$ denotes stable isomorphism and $L$ a graded MCM $A$-module.
\item[(c)] $
(\Omega_A\mathcal{P}_{i}(A))^*\subseteq\mathcal{P}_{ii}(A)$ and
$(\Omega_A\mathcal{P}_{ii}(A))^*\subseteq\mathcal{P}_{i}(A)$.
\end{enumerate}
\end{lem}
\begin{proof}[sketch of proof] $(a)$ follows immediately from Lemma \ref{shortexactDL}. $(b)$ is standard; left for readers. $(c)$ follows from Proposition~\ref{Omega(M)^*isDL}. 
\end{proof}

\begin{remark}

Let $B$ be a Cohen-Macaulay ring, either local or positively
graded, and set $R=B[[t]]$ in the local case and $R=B[t]$
in the graded case. If $M$ is an orientable MCM $B$-module, then
$M_R:=M\otimes_B R$ is an orientable MCM $R$-module. 
\end{remark}
\begin{proof}
Since $R$ is flat over $B$, $M\otimes_B R$ is MCM, and by the going down theorem we get it to be locally free in depth 1.
If $\rank M=r$, then $
\det(M\otimes_B R)=\big(\textstyle\bigwedge^r_R (M\otimes_B R)\big)^{**}
 \cong \det_B(M)\otimes_B R$.
\end{proof}

In the preceding section, we showed that every orientable
MCM module over a two-dimensional Gorenstein local or graded-local ring, under the stated hypotheses, is stably isomorphic to a doubly layered module. However, the counterexamples in Section~3 show that this result does not extend to arbitrary three-dimensional Gorenstein local rings. For example, we can take $R=\mathbb R[[x_1,...,x_4]]/(x_1^2+...+x_4^2)$, then any orientable MCM module that is not free is not stably doubly layered; see Section 3.2. In the following theorem, we prove that every orientable MCM module over $B[t]$ or $B[[t]]$ is stably isomorphic to a doubly layered module, where $B$ is a two dimensional Gorenstein local ring as in the previous section. Lets collect some essential observations needed for the proof in the following lemmas.  

\begin{lem}\label{lem:9.4}
Let $B$ be a two-dimensional Gorenstein standard graded $k$-algebra, or a
two-dimensional complete local Gorenstein $k$-algebra with residue field $k$.
Every orientable MCM $B$-module belongs to both $\PiC(B)$ and $\PiiC(B)$.
Set $R=B\langle t\rangle$ and $\mathfrak a=\m_B R$, where $\m_B=B_+$ in the graded
case. Then,
\begin{enumerate}
\item[$(i)$] if $Z$ is a
one-dimensional Cohen-Macaulay $R$-module and $\Supp Z\subseteq V(\mathfrak a)$, then $(\Omega^2_R(Z))^*$, and $ \Omega_R((\Omega_R^2(Z))^*)\in\PiiC(R)$;
\item[$(ii)$] if $L$ is finite length $R$-module, then $(\Omega_R^3L)^*\in\PiiC(R)$.
\end{enumerate}
\end{lem}

\begin{proof}
The first part is Theorem~\ref{2dimDL}. Now, set $Q=R/\mathfrak a=k\langle t\rangle$.

For (i), since $\Supp_R Z\subseteq V(\mathfrak a)$, we have
$\mathfrak a^sZ=0$ for some $s$, and hence $Z$ is finite over
$Q=k\langle t\rangle$. Moreover, $\Ass_R Z=\{\mathfrak a\}$. 
To obtain a filtration by $R$-submodules for $Z$, set
\[
W_j=\bigcup_{n\ge 0}(\mathfrak a^jZ:_Z t^n),
\qquad 0\leq j\leq s.
\]
Since $t$ is $Z$-regular and $\mathfrak a^sZ=0$,
we have
\[
Z=W_0\supseteq W_1\supseteq\cdots\supseteq W_s=0.
\]
Moreover, $\mathfrak a W_j\subseteq W_{j+1}$,
and $W_{j+1}$ is $t$-saturated in $Z$.
Thus, each $W_j/W_{j+1}$ is a finite torsion-free
$Q$-module; hence it is a (graded) free $Q$-module.

For each $j$, choose a (homogeneous) $Q$-basis $\overline e_1,\ldots,\overline e_r$ of $W_j/W_{j+1}$. Since $\mathfrak a$ annihilates this quotient, the inverse images in $W_j$ of the submodules $Q\overline e_1+\cdots+Q\overline e_i$ are $R$-submodules. These give a filtration from $W_{j+1}$ to $W_j$ whose successive quotients are isomorphic to $Q$ up to degree shifts. Combining these filtrations yields \[ 0=Z_0\subseteq Z_1\subseteq\cdots\subseteq Z_h=Z, \qquad Z_i/Z_{i-1}\cong Q(-d_i). \]
In the complete local case, the grading shifts are omitted.

Since $Q\cong k\otimes_B R$ and $R$ is flat over $B$, we get $(\Omega_R^j(Q))^*
\cong(\Omega_B^j(k))^*\otimes_B R$, and $\Omega_R((\Omega_R^2Q)^*)\cong\Omega_B((\Omega_B^2k)^*)\otimes_B R$.
The two underlying $B$-modules are orientable MCM; hence they belong to
$\mathcal P_{ii}(B)$, and by the preceding remark, their extensions are
orientable MCM $R$-modules. Since, flat extension of scalars preserves
$\mathcal P_{ii}$, both displayed $R$-modules lie in $\mathcal P_{ii}(R)$.

Now apply the dual horseshoe lemma successively to
$0\to Z_{j-1}\to Z_j\to Q(-d_j)\to0$. So, up to free summands, this gives $0\to(\Omega_R^2Q)^*(d_j)\to(\Omega_R^2Z_j)^*\to(\Omega_R^2Z_{j-1})^*\to0$.
Applying the horseshoe lemma once more gives the corresponding exact
sequences for their first syzygies. By induction on $j$, using Lemma~\ref{shortexactDL}, we get $(\Omega_R^2Z)^*$, and $
\Omega_R((\Omega_R^2Z)^*)$ in $\mathcal P_{ii}(R)$. 

For (ii), first let $L=k$. The mapping cone of $0\to Q(-1)\xrightarrow{t}Q\to k\to0$, so Lemma~\ref{lem:syz} gives an exact sequence $0\to\Omega_R^3Q\to\Omega_R^3k\oplus F\to\Omega_R^2Q(-1)\to0$, for some graded free module $F$. This is an exact sequence of orientable MCM module; hence,  dualizing gives $$0\x(\Omega_R^2Q)^*(1)\x(\Omega_R^3k)^*\oplus F^*\x(\Omega_R^3Q)^*\x0.$$ 

Again, by the flatness, we have $(\Omega_R^j(Q))^* \cong(\Omega_B^j(k))^*\otimes_B R$. So, by Theorem~\ref{2dimDL} both $(\Omega_B^2(k))^*$ and $(\Omega^3_B(k))^*$ are in $\mathcal{P}_{ii}$. Since flat extension preserves the two types of layers and their filtrations, both $(\Omega_R^2(Q))^*$ and $(\Omega_R^3(Q))^*$ belong to $\mathcal P_{ii}(R)$. Hence, by Lemma~\ref{shortexactDL}, we get $(\Omega^3_R(k))^*\in\mathcal{P}_{ii}$. For $\ell(L)\geq 2$, the result follows by induction on $\ell(L)$, using a composition series of $L$ and Lemma~\ref{shortexactDL}.
\end{proof}

The next lemma allows us to replace a Cohen--Macaulay
ideal by a larger one, provided the two ideals agree
away from $V(\mathfrak a)$. The point is that their
difference is supported on $V(\mathfrak a)$, where
Lemma~\ref{lem:9.4}(i) applies.

\begin{lem}\label{lem:passing-to-larger-ideal}
Let $R=B\langle t\rangle$ and $\mathfrak a=\mathfrak m_B R$ be as in Lemma~9.4. Let $J\subseteq I$ be height-two Cohen-Macaulay ideals of $R$, homogeneous in the graded case. Suppose that $J$ and $I$ agree outside $V(\mathfrak a)$ i.e., $I_\mathfrak p=J_\mathfrak p$ for all $\mathfrak p\notin V(a)$. Then $(\Omega_R^2(R/J))^*\in\mathcal P_{ii}(R)\implies(\Omega_R^2(R/I))^*\in\mathcal P_{ii}(R)$.
\end{lem}

\begin{proof}
Set $Z=I/J$. The exact sequence
\[
0\longrightarrow Z\longrightarrow R/J
\longrightarrow R/I\longrightarrow0
\]
shows, by the depth lemma, that $Z$ is one-dimensional Cohen--Macaulay, unless $Z=0$. By hypothesis, $\operatorname{Supp}_R Z \subseteq V(\mathfrak a)$, so Lemma~\ref{lem:9.4}(i) gives \[ (\Omega_R^2Z)^*,\quad  \Omega_R((\Omega_R^2Z)^*)\in\mathcal P_{ii}(R). \]
Taking second syzygies and dualizing the preceding sequence gives, up to free summands,
\[
0\longrightarrow(\Omega_R^2(R/I))^*
\longrightarrow(\Omega_R^2(R/J))^*\oplus F
\longrightarrow(\Omega_R^2Z)^*
\longrightarrow0.
\]

Let $C=(\Omega_R^2Z)^*$, and choose a free presentation
\[
0\longrightarrow\Omega_R C
\longrightarrow F'\longrightarrow C
\longrightarrow0.
\]
Pulling back the preceding surjection onto $C$ gives a module $P$ fitting into exact sequences
\[0\longrightarrow\Omega_R C \longrightarrow P \longrightarrow(\Omega_R^2(R/J))^*\oplus F \longrightarrow0 \]
and \[ 0\longrightarrow(\Omega_R^2(R/I))^* \longrightarrow P\longrightarrow F' \longrightarrow0. \]
The latter splits since $F'$ is free. Hence $P\cong(\Omega_R^2(R/I))^*\oplus F'$, giving
\[ 0\longrightarrow\Omega_R C \longrightarrow(\Omega_R^2(R/I))^*\oplus F' \longrightarrow(\Omega_R^2(R/J))^*\oplus F \longrightarrow0.
\]
 The outer terms belong to $\mathcal P_{ii}(R)$, so the assertion follows from Lemma~\ref{lem:closure}.
\end{proof}

From the proof of Artinian ideals being glicci, one can see that it was enough to prove $\Omega^2_R(k)^*$ is stably isomorphic to a doubly layered module, where $R$ is a 2 dimensional local ring. Infact, this is a type~(ii) module. Similar to this, to show the result for $R$ lets first reduce it to showing for height 2 prime ideals of $R$. 

\begin{lem}\label{lem:reduction-to-primes}
Let $R=B\langle t\rangle$ be as in Lemma~\ref{lem:9.4}. Suppose that \[ (\Omega_R^2(R/\mathfrak p))^* \in\mathcal P_{ii}(R) \] for every height-two prime ideal $\mathfrak p$ of $R$, homogeneous in the graded case. Then every orientable MCM $R$-module belongs to $\mathcal P_i(R)\cap\mathcal P_{ii}(R)$. \end{lem}

\begin{proof}
We prove the assertion in four steps. We suppress the grading shifts to keep the proof readable. 

\smallskip
\noindent
\textbf{1. Reduction to one-dimensional
Cohen-Macaulay modules.}

Let $M$ be an orientable MCM $R$-module, and write $T(M)=(\Omega_R M)^*$. We may assume that $M$ is not free. Then there is an $N$-type resolution \[ 0\longrightarrow F\longrightarrow M \longrightarrow I(h)\longrightarrow 0, \]
where $F$ is free and $I$ is a height-two Cohen-Macaulay ideal of $R$. Hence, $\Omega_R M\cong_{\mathrm{st}} \Omega_R^2(R/I)(h)$. Dualizing, we obtain $\Omega(M)^*\cong_{\mathrm{st}} (\Omega_R^2(R/I))^*(-h)$. Hence, it is enough to establish
\[ (\Omega_R^2 U)^*\in\mathcal P_{ii}(R) \]
for every one-dimensional Cohen-Macaulay $R$-module $U$. Indeed, applying $(*)$ to $U=R/I$ shows that $T(M)\in\mathcal P_{ii}(R)$. By Lemma~\ref{lem:closure}, we have $T^2(M)\cong_{\mathrm{st}}M$ and $T(\mathcal P_{ii}(R))\subseteq\mathcal P_i(R)$, so $M\in\mathcal P_i(R)$. Moreover, $T(M)$ is again an orientable MCM module. Applying the same argument to $T(M)$ gives $T^2(M)\in\mathcal P_{ii}(R)$, and hence $M\in\mathcal P_{ii}(R)$ as well.\\

\noindent\textbf{2. A filtration by rank-one modules.}

Let $U$ be a one-dimensional Cohen--Macaulay module.
Choose a $U$-regular parameter $\ell$ such that $U$
is finite over $k\langle\ell\rangle$, with $\ell$
linear in the graded case. In the complete local case,
$U_\ell$ is finite-dimensional over $k((\ell))$,
and hence admits a composition series over $R_\ell$.

In the graded case, set $A=(R_\ell)_0$, $E=(U_\ell)_0$. Since $U$ is finite over $k[\ell]$, the $k$-vector space $E$ is finite-dimensional. Choose a composition series of $A$-modules \[ 0=E_0\subsetneq E_1\subsetneq\cdots \subsetneq E_s=E, \qquad E_i/E_{i-1}\cong A/\mathfrak n_i, \]
where each $\mathfrak n_i$ is a maximal ideal of $A$. Since $\ell$ has degree one and is invertible in $R_\ell$, we have $R_\ell=A[\ell,\ell^{-1}]$, $U_\ell=\bigoplus_{n\in\mathbb Z}\ell^n E$. For each $i$, define $W_i=\bigoplus_{n\in\mathbb Z}\ell^n E_i$. These are graded $R_\ell$-submodules: indeed, for homogeneous $r\in(R_\ell)_d$ and $e\in E_i$, we have $r\ell^n e=\ell^{n+d}(\ell^{-d}r)e \in W_i$, since $\ell^{-d}r\in A$.
Thus we obtain a filtration
\[
0=W_0\subsetneq W_1\subsetneq\cdots
\subsetneq W_s=U_\ell,
\]
with $W_i/W_{i-1}\cong (A/\mathfrak n_i)[\ell,\ell^{-1}]$. Each quotient is graded-simple, since $A/\mathfrak n_i$ is a field and every nonzero homogeneous element of this quotient is invertible.

In either case, write the resulting filtration as
\[
0=W_0\subsetneq W_1\subsetneq\cdots
\subsetneq W_s=U_\ell.
\]
Each quotient is isomorphic, up to a grading
shift, to $(R/\mathfrak p_i)_\ell$ for some
prime $\mathfrak p_i$ not containing $\ell$.
Since $\mathfrak p_i\in\operatorname{Supp}U$
and $\dim U=1$, we have
$\dim R/\mathfrak p_i=1$, so
$\operatorname{ht}\mathfrak p_i=2$.

Since $\ell$ is $U$-regular, we may regard
$U\subseteq U_\ell$. Setting $U_i=U\cap W_i$
gives a filtration
\[
0=U_0\subsetneq U_1\subsetneq\cdots
\subsetneq U_s=U.
\]

Moreover, $(U_i)_\ell=W_i$: indeed, every $w\in W_i$ can be written as $u/\ell^n$ for some $u\in U$, and then $u=\ell^nw\in U\cap W_i=U_i$. The natural map $V_i=U_i/U_{i-1}\x W_i/ W_{i-1}$ is injective, since $U_i\cap W_{i-1}=U_{i-1}$, and becomes an isomorphism after inverting $\ell$. Thus $V_i$ embeds into $(R/\mathfrak p_i)_\ell$ and has rank one over $R/\mathfrak p_i$. Since $R/\mathfrak p_i$ is a one-dimensional domain, $V_i$ is torsion-free and hence Cohen-Macaulay.

\smallskip
\noindent
\textbf{3. Reduction to the prime quotients.}

Fix $i$, and put $D=R/\mathfrak p_i$ and
$V=V_i$. Recall that $V$ embeds into $D_\ell$
and that $V_\ell\cong D_\ell$, up to a grading
shift. Since $V$ is finitely generated, we may
multiply this embedding by a sufficiently large
power of $\ell$ to obtain an injection
$V\hookrightarrow D$. This map remains an
isomorphism after inverting $\ell$. Hence its
cokernel $L$ satisfies $L_\ell=0$, so
$\ell^nL=0$ for some $n\geq 1$. Since $D$ is
a one-dimensional domain and $\ell$ is a
parameter on $D$, the ring $D/\ell^nD$
has finite length. Thus, suppressing grading
shifts, we obtain
\[
0\longrightarrow V\longrightarrow D
\longrightarrow L\longrightarrow 0,
\qquad \operatorname{length}_R L<\infty.
\]

Now, applying Lemma~\ref{lem:syz}(a-2), and then taking dual we get 
\[
0\longrightarrow(\Omega_R^2D)^*
\longrightarrow(\Omega_R^2V)^*\oplus F
\longrightarrow(\Omega_R^3L)^*
\longrightarrow0.
\]
The outer terms belong to
$\mathcal P_{ii}(R)$ by hypothesis and previous lemma
Lemma~\ref{lem:9.4}(ii), respectively. Hence
Lemma~\ref{lem:closure} gives
\[
(\Omega_R^2V_i)^*\in\mathcal P_{ii}(R)
\qquad\text{for every }i.
\]
\smallskip
\noindent\textbf{4. Passing from the filtration to $U$.}

For each $i$, the exact sequence \[ 0\longrightarrow U_{i-1} \longrightarrow U_i\longrightarrow V_i \longrightarrow0 \] and after taking second syzygies and dualizing it, we get \[ 0\longrightarrow(\Omega_R^2V_i)^*\longrightarrow(\Omega_R^2U_i)^*\oplus F_i \longrightarrow(\Omega_R^2U_{i-1})^* \longrightarrow0, \]
where $F_i$ is free. Now, by induction on $i$, we get that the outer terms belong to $\mathcal P_{ii}(R)$, (starting with $U_0=0$). So, by Lemma~\ref{lem:closure} we get $(\Omega^2(U_i))^*\in\mathcal{P}_{ii}$. Hence, the required follows from part~1.
\end{proof}

Now, we note the following one generator technique to work with prime ideals. This is a key part in proofs for both one and two dimensional ideals being glicci.

\begin{lem}\label{lem:one-equation}
Let $R$ be a Noetherian graded ring and $A=R[v]$ with
$\deg v=1$, or let $R$ be complete local and $A=R[[v]]$.
Let $P\subset A$ be a prime ideal, homogeneous in the
graded case, and set $\mathfrak q=P\cap R$.
Suppose $P\ne\mathfrak qA$. In the complete case, assume
also that $\dim R/\mathfrak q\leq1$ or $A/P$ is finite
over $R/\mathfrak q$.
Then there exist $f\in R\setminus\mathfrak q$ and
$\mu\in P$ (homogeneous in the graded case), such that
\[
P_f=(\mathfrak qA+(\mu))_f.
\]
\end{lem}

\begin{proof}
Put $D=R/\mathfrak q$ and $Q=P/\mathfrak qA$.
Then $D$ is a domain and $Q\cap D=0$.
In the graded case, choose a homogeneous polynomial
$\bar\mu\in Q$ of least $v$-degree, and let $\bar f\in D$
be its nonzero leading coefficient. After inverting
$\bar f$, polynomial division by the resulting monic
polynomial shows that $Q_{\bar f}=(\bar\mu)_{\bar f}$.

In the complete case, we first find a distinguished polynomial \[ h=v^d+a_{d-1}v^{d-1}+ cdots+a_0\in Q, \qquad a_i\in\mathfrak m_D. \] If $\dim D=0$, take $h=v$. If $\dim D=1$, the nonzero prime $Q$ is not contained in $\mathfrak m_DD[[v]]$; Thus, Weierstrass preparation (see \cite[Chapter~9,  Section~9, Theorem~9.2]{lang2012algebra}) gives such an $h$. If $A/P\cong D[[v]]/Q$ is finite over $D$, the quotient $k[[v]]/\overline Q$ is finite-dimensional over $k=D/\mathfrak m_D$. As $k[[v]]$ is a DVR, $\overline Q=(v^d)$ for some $d\geq1$, so multiplication by $v$ is nilpotent modulo $\mathfrak m_D$. The determinant trick therefore gives a distinguished polynomial in $Q$.

Weierstrass division states that every power series
has a unique remainder of $v$-degree less than $d$
upon division by $h$ (see \cite[Chapter~9, Section~9,Theorem~9.1]{lang2012algebra}). 
So, we have $D[[v]]/(h)\cong D[v]/(h)$, which is a finite free $D$-module. Hence $Q/(h)$ is finitely generated over $D$. Choose polynomial representatives $g_1,\ldots,g_s\in Q$ of its generators, so that $Q=(h,g_1,\ldots,g_s)$. 

Let $K$ be the fraction field of $D$. Since $Q\cap D=0$, the extension of $Q/(h)$ to $K[v]/(h)$ is a prime ideal, generated by a monic irreducible polynomial $p\in K[v]$ dividing $h$ and each $g_i$ in $K[v]$.

Since $p$ belongs to this extended ideal, there is a nonzero $c\in D$ such that $\bar\mu=cp\in Q$. Choose a nonzero $f\in D$, divisible by $c$, which clears all denominators in the coefficients of $p$, $h/p$, and $g_i/p$ for $1\leq i\leq s$. Then $\bar\mu\in Q_f$, and each of $h,g_1,\ldots,g_s$ is a multiple of $\bar\mu$ in $D_f[v]$. Thus $Q_f=(\bar\mu)_f$. 

Lifting $f$ and $\bar\mu$ to $R$ and $P$, respectively, gives $P_f=(\mathfrak qA+(\mu))_f$, as required.

In either case, lifting $\bar f$ and $\bar\mu$ to
$f\in R$ and $\mu\in P$ gives the assertion.
\end{proof}

\begin{thm}\label{prop:one-variable}
For $B$ as in Lemma~\ref{lem:9.4}, every orientable MCM module over
$R=B\langle t\rangle$ belongs to both $\PiC(R)$ and $\PiiC(R)$.
\end{thm}

\begin{proof}
We divide the proof into four steps. Again we suppress the grading shifts to keep it readable.

\smallskip \noindent \textbf{1. Reduction to an ideal generated by one element modulo $\mathfrak qR$.}

By Lemma~\ref{lem:reduction-to-primes}, it suffices to prove $(\Omega_R^2(R/\mathfrak p))^* \in\mathcal P_{ii}(R)$ for every height-two prime ideal $\mathfrak p$ of $R$, homogeneous primes in the graded case. Set $\mathfrak a=\mathfrak m_B R$, $\mathfrak q=\mathfrak p \cap B$.

If $\mathfrak q=\mathfrak m_B$, then $\mathfrak p=\mathfrak a$, and the assertion follows from Lemma~9.4(i). We may therefore assume that $\mathfrak q\ne\mathfrak m_B$. Then $D=B/\mathfrak q$ is a one-dimensional Cohen-Macaulay domain.

Since $D=B/\mathfrak q$ is a one-dimensional Cohen--Macaulay
domain, Lemma~\ref{lem:one-equation} gives
$f\in B\setminus\mathfrak q$ and $g\in\mathfrak p$,
homogeneous in the graded case, such that, on setting
$J=\mathfrak qR+(g)$, we have $J_f=\mathfrak p_f$. 

Since $g\notin\mathfrak qR$, it is regular on the domain
$R/\mathfrak qR$. As $R/\mathfrak qR$ is Cohen--Macaulay
of dimension two, $J$ is a height-two Cohen-Macaulay ideal.

\smallskip \noindent \textbf{2. $J$ is directly Gorenstein linked to a Gorenstein ideal.}

Set $S=R/\mathfrak qR=D\langle t\rangle$, and write $\overline g$ for the image of $g$ in $S$. Since $D$ is a one-dimensional Cohen-Macaulay domain, its canonical module is a torsion-free module of rank one. Thus, by clearing denominators, we may identify it with a canonical ideal $\mathfrak c\subseteq D$, this can be chosen to be homogeneous in the graded case.

Let $H_0$ be the inverse image of $\mathfrak c$ in $B$, and set
\[
H=H_0R,\qquad G=\mathfrak qR+gH.
\]
We first show that $H$ and $G$ are Gorenstein ideals. Since $S$ is Cohen-Macaulay with
canonical module $\mathfrak cS$, and $\overline g$ is a nonzerodivisor on $S$; both $\mathfrak cS$ and $\overline g\mathfrak cS$ are canonical ideals of $S$. Their quotients are Gorenstein, and hence $R/H\cong S/\mathfrak cS$, $R/G\cong S/\overline g\mathfrak cS$ are one-dimensional Gorenstein rings. Thus $H$ and $G$ are height-two
Gorenstein ideals.

It remains to verify the linkage. Since $\mathfrak cS$ is a canonical ideal, we have $(\mathfrak cS:\mathfrak cS)=S$. As $\overline g$ is a nonzerodivisor, $(\overline g\mathfrak cS:\overline g) =\mathfrak cS$, $ (\overline g\mathfrak cS:\mathfrak cS) =\overline gS$. Taking inverse images in $R$ gives \[ G:J=H,\qquad G:H=J. \] Therefore $J$ and the Gorenstein ideal $H$ are directly Gorenstein linked by $G$.

\smallskip\noindent\textbf{3. Proving the assertion for $J$}

Choose an $N$-type resolution of $H_0$ over $B$ and extend it to $R$, obtaining \[0\longrightarrow L\longrightarrow N\longrightarrow  H\longrightarrow0,\] where $L$ is free and $N$ is extended from a rank-two orientable MCM $B$-module. Also, since $G$ is Gorenstein, choose \[ 0\longrightarrow R \longrightarrow E_G \longrightarrow G\longrightarrow0, \] where $E_G$ is a rank-two orientable MCM $R$-module. Grading shifts are suppressed throughout.

Since $\operatorname{Ext}_R^1(E_G,L)=0$,
the inclusion $G\subseteq H$ lifts to a
map $E_G\to N$. Together with the map
$L\to N$, this gives an exact sequence
\[
0\longrightarrow R
\longrightarrow E_G\oplus L
\longrightarrow N
\longrightarrow H/G\longrightarrow0.
\]
Now $H/G$ is the canonical module of $R/J$. Dualizing the sequence, and using $\operatorname{Ext}_R^2(H/G,R)\cong R/J$, we obtain \[ 0\longrightarrow N^* \longrightarrow E_G^*\oplus L^* \longrightarrow J\longrightarrow0. \]

Choose a free presentation $0\longrightarrow\Omega_R N \longrightarrow F\longrightarrow N \longrightarrow0$. Dualizing and pushing out along $N^*\to E_G^*\oplus L^*$ gives a module $N_J$ fitting into exact sequences
\[
0\longrightarrow E_G^*\oplus L^*
\longrightarrow N_J
\longrightarrow(\Omega_R N)^*
\longrightarrow0
\]
and
\[
0\longrightarrow F^*
\longrightarrow N_J\longrightarrow J
\longrightarrow0.
\]

Since $E_G$ is rank-two orientable MCM,
$E_G^*\oplus L^*\in\mathcal P_i(R)$. Moreover, $N$ is extended from an orientable MCM $B$-module, so Lemma~\ref{lem:9.4} and Lemma~\ref{lem:closure} give $N\in\mathcal P_{ii}(R)$, $(\Omega_R N)^*\in\mathcal P_i(R)$.  

The first sequence therefore shows that
$N_J\in\mathcal P_i(R)$.
The second is an $N$-type resolution of
$J$, and hence
\[
(\Omega_R^2(R/J))^*
\cong_{\mathrm{st}}(\Omega_R N_J)^*
\in\mathcal P_{ii}(R),
\]
again by Lemma~\ref{lem:closure}.

\smallskip \noindent\textbf{3. Passing from $J$ to $\mathfrak p$.}

We now use the fact that $J$ and $\mathfrak p$ agree after inverting $z$. 
Set $Z=\mathfrak p/J$. Since $J_f=\mathfrak p_f$,
some power of $f$ annihilates $Z$. Moreover,
$\mathfrak qZ=0$, so
\[
\operatorname{Supp}_R Z
\subseteq V((\mathfrak q+(f))R).
\]
Since $D=B/\mathfrak q$ is a one-dimensional domain
and $f\notin\mathfrak q$, the quotient
$B/(\mathfrak q+(f))$ has dimension at most zero.
Hence $\operatorname{Supp}_R Z\subseteq V(\mathfrak m_BR)$.

Thus $J$ and $\mathfrak p$ agree outside $V(\mathfrak a)$. By Step~2 and Lemma~\ref{lem:passing-to-larger-ideal}, we obtain \[ (\Omega_R^2(R/\mathfrak p))^* \in\mathcal P_{ii}(R). \] The assertion for every orientable MCM $R$-module now follows from Lemma~\ref{lem:reduction-to-primes}.
\end{proof}

Now we are ready to prove \ref{mainsection9}(a). The previous result and Noether normalizations are the tools that we need for the proof.

\begin{proof}[Proof of Theorem~\ref{mainsection9}\textup{(a)}]
Set $C=S/I$. Since $C$ is one-dimensional and
Cohen-Macaulay; so, we may choose a linear form $\ell$
that is $C$-regular. In the graded case, Noether
normalization gives that $C$ is finite over $k[\ell]$.
In the complete local case, $\ell$ is a parameter
of $C$, so $C$ is finite over $k[[\ell]]$ by
complete local Noether normalization.

After a linear change of coordinates, we may assume
$\ell=x_{n-1}$, and set $t=x_n$, and write
\[
S=T\langle t\rangle,\qquad
T=k\langle x_1,\ldots,x_{n-1}\rangle,
\qquad J=I\cap T.
\]
Since $T/J$ embeds into $C$; so, $T/J$ is finite over
$k\langle\ell\rangle$. Hence, $\dim T/J=1$ and
$\operatorname{ht}_T J=n-2$.

Choose a (homogeneous) $T$-regular sequence
$f_1,\ldots,f_{n-3}\in J$, and set
\[
K=(f_1,\ldots,f_{n-3}),\qquad
B=T/K,\qquad A=S/KS=B\langle t\rangle.
\]
Then $B$ is a two-dimensional Gorenstein ring,
and $I/KS$ is a height-two Cohen-Macaulay
ideal of $A$. By Theorem~\ref{prop:one-variable} and
Theorem~\ref{DLstablyGlicci}, $I/KS$ is glicci.
Taking inverse images of the Gorenstein links
under $S\to A$ shows that $I$ is glicci.
\end{proof}

\begin{question}\label{DLextensionvariablequestion}
    Let $B$ be a two dimensional local Gorenstein ring, with an infinite residue field. Let $R=B[[t_1,..,t_n]]$ or $R=B[t_1,..,t_n]$. Is every orientable MCM $R$-module stably isomorphic to a doubly layered $R$-module?
\end{question}

We next prove that every two-dimensional Cohen-Macaulay ideal in $S$ is (homogeneously) glicci. As in the one-dimensional case, we divide the argument into several smaller statements to make the proof easier to follow. Although the overall strategy is similar, the technical details are substantially different. In particular, extending Theorem~\ref{prop:one-variable} to $B\langle u,v\rangle$ appears to be considerably more difficult than the one-variable case, and we do not know whether such an extension holds. Instead, we use a different construction, guided by the same ideas, to prove Theorem~\ref{mainsection9}(b).

Let $k$ be an infinite field, and let $S=k[x_1,\ldots,x_n]$ with its standard grading, or $S=k[[x_1,\ldots,x_n]]$. Let $I\subsetneq S$ be an ideal such that $S/I$ is Cohen--Macaulay of dimension two, and assume that $I$ is homogeneous in the graded case. The cases $n\leq4$ are already licci, so assume $n\geq5$.

Choose a complete intersection $K\subseteq I$ of height $n-4$, homogeneous in the graded case, and set $A=S/K$ and $J=I/K$. Thus $\operatorname{ht}I=\operatorname{ht}K+2$, the ring $A$ is four-dimensional and Gorenstein, and $J$ is a height two Cohen-Macaulay ideal. It suffices to prove that $J$ is glicci in $A$. Choose an N-type resolution
\[
0\longrightarrow F\longrightarrow M\longrightarrow J\longrightarrow0,
\]
where $F$ is free, $M$ is orientable, and MCM. By Theorem~\ref{DLstablyGlicci}, it is enough to show that $M\in\mathcal P_i(A)$.

Unlike in dimension one, the prime quotients occurring in a prime filtration of $S/I$ need not be Cohen--Macaulay. We therefore need a way to apply our MCM module classes to arbitrary finite modules. Since $\dim A=4$, the fourth syzygy $\Omega_A^4U$ is MCM for every finitely generated $A$-module $U$. Taking four syzygies, however, also changes the stable class of a module that is already MCM. Over a hypersurface, Eisenbud's periodicity theorem controls this change, this is a key tool in Section~10.1. For a general Gorenstein ring, we compensate by taking four cosyzygies.

For an MCM module $L$ over a local or graded-local Gorenstein ring $R$, define its \emph{first cosyzygy} by $\Omega_R^{-1}L=(\Omega_R(L^*))^*$, where $(-)^*=\operatorname{Hom}_R(-,R)$, using graded Hom in the graded case. Dualizing a free presentation of $L^*$ gives an exact sequence $0\to L\to G\to\Omega_R^{-1}L\to0$ with $G$ free. Thus cosyzygy is inverse to syzygy up to stable isomorphism; higher cosyzygies are defined by iteration. For a finitely generated $A$-module $U$, set
\[
[U]:=\Omega_A^{-4}\Omega_A^4U.
\]
This is an MCM module defined up to stable isomorphism, and $[U]\simeq_{\mathrm{st}}U$ when $U$ is already MCM. The N-type resolution above gives
\[
\Omega_A^4(A/J)\simeq_{\mathrm{st}}\Omega_A^3M,\qquad [A/J]\simeq_{\mathrm{st}}\Omega_A^{-1}M.
\]

\begin{remark}\label{rem:cosyzygy-types}
Let $R$ be a Noetherian local Gorenstein ring of dimension at least two with infinite residue field, or a positively graded Gorenstein ring $R=\bigoplus_{n\geq0}R_n$ such that $(R_0,\mathfrak m_0)$ is local with infinite residue field and $\operatorname{grade}_R(R_nR)\geq2$ for all $n\gg0$.  Then
\[
\Omega_R^{-1}\mathcal P_i(R)\subseteq\mathcal P_{ii}(R),\qquad \Omega_R\mathcal P_{ii}(R)\subseteq\mathcal P_i(R).
\]
\end{remark}

\begin{proof}
Duality preserves $\mathcal P_i(R)$; it preserves rank two orientable MCM modules and reverses exact MCM filtrations. By Proposition~\ref{Omega(M)^*isDL}, the operation $L\mapsto(\Omega_RL)^*$ sends $\mathcal P_i(R)$ into $\mathcal P_{ii}(R)$ and $\mathcal P_{ii}(R)$ into $\mathcal P_i(R)$. Applying the first assertion to $L^*$ proves the cosyzygy inclusion, and dualizing the second proves the syzygy inclusion. Graded duals and graded free presentations give the same argument in the graded case, with all shifts retained.
\end{proof}

The preceding remark shows that $M\in\mathcal P_i(A)$ if and only if $[A/J]\in\mathcal P_{ii}(A)$. We therefore call a finitely generated $A$-module $U$ \emph{of type~(ii)} if $[U]\in\mathcal P_{ii}(A)$. For MCM modules, this agrees with membership in the stable filtration class $\mathcal P_{ii}(A)$. We call $U$ \emph{strongly type~(i)} if
\[
\Omega_A^j[U]\in\mathcal P_i(A)\qquad\text{for every }j\in\mathbb Z,
\]
where negative powers denote cosyzygies. In particular, a strongly type~(i) module is of type~(ii), since $[U]\simeq_{\mathrm{st}}\Omega_A^{-1}\Omega_A[U]$ and $\Omega_A[U]\in\mathcal P_i(A)$. These definitions let us work directly with the factors of a prime filtration: to prove that $J$ is glicci, it suffices to show that $A/J$ is of type~(ii).

The following lemma records the closure properties needed for this argument. 

\begin{lem}\label{lem:type-calculus}
Let $k$ be an infinite field, and let $B$ be a two-dimensional standard graded Gorenstein $k$-algebra, or a two-dimensional complete local Gorenstein $k$-algebra with residue field $k$. In the graded case, set $R=B[u]$ and $A=R[v]$, with $\deg u=\deg v=1$; in the local case, set $R=B[[u]]$ and $A=R[[v]]$. Thus $\dim R=3$ and $\dim A=4$. 

\par\smallskip\noindent\textup{(a)} In an exact sequence $0\to E\to V\to U\to0$ of $A$-modules, if two terms are strongly type~(i), then so is the third. Modules of type~(ii) are closed under extensions. If $V$ is of type~(ii) and $[E]\in\mathcal P_i(A)$, then $U$ is of type~(ii).

\par\smallskip\noindent\textup{(b)} If $N$ is an $R$-module with $\dim_R N\leq1$, then $N\otimes_R A$ is strongly type~(i).

\par\smallskip\noindent\textup{(c)} If $G\subsetneq A$ is a Gorenstein ideal of height two, then $A/G$ is of type~(ii). If $J\subsetneq A$ is a Cohen--Macaulay ideal of height two and $A/J$ is of type~(ii), then $J$ is (homogeneously) glicci.
\end{lem}

\begin{proof}
\textup{(a)} The horseshoe lemma, followed by four cosyzygies, gives an exact sequence $$0\to[E]\to[V]\oplus F\to[U]\to0$$ after choosing stable representatives, where $F$ is free. This proves the extension closure of type~(ii). Push this sequence out along $0\to[E]\to F_2\to\Omega_A^{-1}[E]\to0$, where $F_2$ is free. Since $[U]$ is MCM and $A$ is Gorenstein, $\operatorname{Ext}_A^1([U],F_2)=0$. The resulting extension of $[U]$ by $F_2$ therefore splits, giving
\[
0\longrightarrow[V]\oplus F_1\longrightarrow[U]\oplus F_2\longrightarrow\Omega_A^{-1}[E]\longrightarrow0.
\]
By Remark~\ref{rem:cosyzygy-types}, the right-hand term belongs to $\mathcal P_{ii}(A)$, and the left-hand term does so by hypothesis. Extension closure and stability under adding free summands imply $[U]\in\mathcal P_{ii}(A)$. 

The first exact sequence also yields, for some free module $F_3$,
\[
0\longrightarrow \Omega_A[U]\longrightarrow [E]\oplus F_3
\longrightarrow [V]\oplus F\longrightarrow 0.
\]
The three exact sequences, together with the extension closure of $\mathcal P_i(A)$,
show that any two strongly type~$(i)$ terms force the third to be strongly
type~$(i)$, by applying arbitrary syzygies and cosyzygies, up to free summands.
This completes the proof of (a).

\par\smallskip
\textup{(b)} Since $\dim R=3$ and $\dim_R N\leq1$, the module $N$ vanishes at every prime of height at most one. Now note that its first syzygy is therefore orientable; in codimension at most one, it agrees with the free module in a free presentation of $N$, so its determinant is free, up to a grading shift. Thus $L=\Omega_R^3N$ is MCM, and all its stable syzygies and cosyzygies are orientable. By the one-variable theorem, Theorem~\ref{prop:one-variable}, they belong to $\mathcal P_i(R)$, with the stably free case included. Flat base change to $A$ commutes with syzygies and cosyzygies up to stable isomorphism and preserves $\mathcal P_i$. Precisely, for every $j\in\mathbb Z$, we have \[\Omega_A^j[N\otimes_R A]\simeq_{\mathrm{st}}\bigl(\Omega_R^{j-3}L\bigr)\otimes_R A\in\mathcal P_i(A).\] Hence $N\otimes_R A$ is strongly type~(i).

\par\smallskip
\textup{(c)} For a height-two Gorenstein ideal $G$, an N-type resolution $0\to A(a)\to M\to G\to0$ with $M$ orientable MCM of rank two shows that $A/G$ is of type~$(ii)$ by Remark~\ref{lem:type-calculus}. The second assertion follows from the preceding discussion.
\end{proof}

\begin{lem}\label{lem:projection-types}
Let $A$, $B$, and $R$ be as in the previous lemma. Let $P\subseteq A$
be a prime ideal, and set $\mathfrak q=P\cap R$ and
$e=\dim R/\mathfrak q$.

\par\smallskip\noindent\textup{(a)} If $e\leq1$, then $A/P$ is strongly type~(i). Consequently, if $\mathfrak b\subseteq R$ satisfies $\dim R/\mathfrak b\leq1$, every finitely generated $A$-module supported on $V(\mathfrak bA)$ is strongly type~(i).

\par\smallskip\noindent\textup{(b)} Suppose that $e=2$ and that $A/P$ is finitely generated as an $R$-module. If there is an ideal $\mathfrak a\subseteq\mathfrak q$ such that $R/\mathfrak a$ is Gorenstein of dimension two and $\mathfrak aR_{\mathfrak q}=\mathfrak qR_{\mathfrak q}$, then $A/P$ is of type~(ii).
\end{lem}

\begin{proof}
(a) If $P=\mathfrak qA$, the assertion follows from
Lemma~\ref{lem:type-calculus}(b). Otherwise,
Lemma~\ref{lem:one-equation} gives $f\in R\setminus\mathfrak q$
and $\mu\in P$ such that
$P_f=(\mathfrak qA+(\mu))_f$. Set $H=\mathfrak qA+(\mu)$.
Since $A/\mathfrak qA$ is a domain, multiplication by $\mu$
gives an exact sequence
\[
0\longrightarrow A/\mathfrak qA(-\deg\mu)
\xrightarrow{\ \mu\ } A/\mathfrak qA
\longrightarrow A/H\longrightarrow0,
\]
with the shift omitted in the local case. Thus $A/H$ is
strongly type~$(i)$ by Lemma~\ref{lem:type-calculus}(a),(b).
Moreover, $(P/H)_f=0$, so
\[
\operatorname{Supp}_A(P/H)
\subseteq V((\mathfrak q+(f))A),
\qquad
\dim R/(\mathfrak q+(f))<e.
\]

We now prove the assertion by induction on $e$.
If $e=0$, then $R/\mathfrak q$ is a field, so
$f\notin\mathfrak q$ acts invertibly on $P/H$, since
$\mathfrak q(P/H)=0$. As $(P/H)_f=0$, we obtain
$P/H=0$, and hence $A/P=A/H$ is strongly type~$(i)$.

Suppose $e=1$. Take a prime filtration of $P/H$,
with factors $A/Q_i$ (up to grading shifts).
Since $\operatorname{Supp}_A(P/H)\subseteq
V((\mathfrak q+(f))A)$, each $Q_i$ contains
$(\mathfrak q+(f))A$. Thus
$\dim R/(Q_i\cap R)=0$, so each $A/Q_i$ is strongly
type~$(i)$ by the case $e=0$. Extension closure
then shows that $P/H$ is strongly type~$(i)$.
Finally, the exact sequence
\[
0\longrightarrow P/H\longrightarrow A/H
\longrightarrow A/P\longrightarrow0
\]
and Lemma~\ref{lem:type-calculus}(a) show that
$A/P$ is strongly type~$(i)$.
The assertion for arbitrary modules follows by
prime filtrations and extension closure.

(b) By Lemma~\ref{lem:one-equation}, choose
$f\in R\setminus\mathfrak q$ and $\mu\in P$ such that
$P_f=(\mathfrak qA+(\mu))_f$.
Since $R/\mathfrak a$ is Gorenstein of dimension two,
$\mathfrak a$ is unmixed of height one. Moreover,
$\mathfrak aR_{\mathfrak q}=\mathfrak qR_{\mathfrak q}$,
so we may write
\[
\operatorname{Ass}_R(R/\mathfrak a)
=\{\mathfrak q,\mathfrak q_1,\ldots,\mathfrak q_r\},
\]
where the primes are distinct and all have height one.

Choose
$z\in(\bigcap_{i=1}^r\mathfrak q_i)\setminus\mathfrak q$
and $w\in\mathfrak q\setminus\bigcup_{i=1}^r\mathfrak q_i$.
If $r=0$, take $z=1$ and $w=0$.
In the graded case, these elements can be chosen
homogeneous and their degrees adjusted using a linear
form outside all these primes so that
$\deg(z\mu)=\deg w$.

Set $g=z\mu+w$ and $H=\mathfrak aA+(g)$.
Note that $g\in P$, so $H\subseteq P$.
Modulo $\mathfrak qA$, the element $g$ equals
$z\mu\ne0$, while modulo each $\mathfrak q_iA$,
it equals $w\ne0$.
Thus $g$ is regular on the Cohen-Macaulay ring
$A/\mathfrak aA$. Since this ring is Gorenstein,
$H$ is a height-two Gorenstein ideal.
Hence $A/H$ is of type~$(ii)$ by
Lemma~\ref{lem:type-calculus}(c).

We claim that $P/H$ is supported only over primes of
$R$ of dimension at most one. Indeed, let
$Q\in\operatorname{Supp}_A(P/H)$ and put
$\mathfrak p=Q\cap R$. Since $\mathfrak a\subseteq\mathfrak p$,
if $\dim R/\mathfrak p=2$, then
$\mathfrak p=\mathfrak q$ or $\mathfrak q_i$ for some $i$.

If $\mathfrak p=\mathfrak q$, then $f,z$ are units in
$A_Q$, and $w\in\mathfrak q$. Hence
\[
H_Q=(\mathfrak qA+(\mu))_Q=P_Q,
\]
contradicting $Q\in\operatorname{Supp}_A(P/H)$.
If $\mathfrak p=\mathfrak q_i$, then
$z\in Q$ and $w\notin Q$, so $g=z\mu+w\notin Q$,
again a contradiction.

Therefore every prime in $\operatorname{Supp}_A(P/H)$
contracts to a prime $\mathfrak p$ with
$\dim R/\mathfrak p\leq1$. By part~(a) and a prime
filtration, $P/H$ is strongly type~$(i)$.
Finally, applying Lemma~\ref{lem:type-calculus}(a) to
\[
0\longrightarrow P/H\longrightarrow A/H
\longrightarrow A/P\longrightarrow0
\]
shows that $A/P$ is of type~$(ii)$.

\end{proof}

\begin{proof}[Proof of Theorem~\ref{mainsection9}(b)]
The height-at-most-two case is licci, so assume $n\geq5$ and put $C=S/I$. All choices are homogeneous in the graded case. We break the proof into two steps to keep it readable.

\textbf{ 1. a common ambient ring.}

Choose $S=R_0\langle v\rangle$, where $R_0=k\langle x_1,\ldots,x_{n-1}\rangle$ contains two linear parameters of $C$; graded or complete local Noether normalization then makes $C$ finite over $R_0$. For each $\mathfrak P\in\operatorname{Min}(I)$, set $\mathfrak q_{\mathfrak P}=\mathfrak P\cap R_0$. Since $S/\mathfrak P$ is finite over $R_0/\mathfrak q_{\mathfrak P}$ and has dimension two, so $\operatorname{ht}\mathfrak q_{\mathfrak P}=n-3$. Choose a complete intersection $\mathfrak c_{\mathfrak P}\subseteq\mathfrak q_{\mathfrak P}$ of height n-3, satisfying $\mathfrak c_{\mathfrak P}(R_0)_{\mathfrak q_{\mathfrak P}}=\mathfrak q_{\mathfrak P}(R_0)_{\mathfrak q_{\mathfrak P}}$. Indeed, prime avoidance lets us choose a regular sequence whose images form a basis of the maximal ideal modulo its square in the regular local ring $(R_0)_{\mathfrak q_{\mathfrak P}}$.

Choose common linear parameters $\ell_1,\ell_2$ on $C$ and all $R_0/\mathfrak c_{\mathfrak P}$; this is possible as the family is finite and $k$ is infinite. Extend them to coordinates so that $R_0=T\langle u\rangle$ with $\ell_1,\ell_2\in T=k\langle y_1,\ldots,y_{n-2}\rangle$. These quotients are finite over $k\langle\ell_1,\ell_2\rangle$, hence over $T$. Hence, the finite inclusions $T/(I\cap T)\hookrightarrow C$ and $T/(\mathfrak c_{\mathfrak P}\cap T)\hookrightarrow R_0/\mathfrak c_{\mathfrak P}$ shows that all the contractions have height $n-4$. Their intersection has the same height, so choose a complete intersection $K\subseteq(I\cap T)\cap\bigcap_{\mathfrak P\in\Ass(S/I)}(\mathfrak c_{\mathfrak P}\cap T)$ of height $n-4$, and set
\[
B=T/K,\qquad R=R_0/KR_0=B\langle u\rangle,\qquad A=S/KS=R\langle v\rangle.
\]
These rings are Gorenstein of dimensions two, three, and four, respectively. Moreover, $J=I/KS$ is Cohen--Macaulay of height two, and $C=A/J$ is finite over $R$.

\textbf{2: a prime filtration.}

Take a prime filtration of $C$, allowing degree shifts in the graded case. Each factor $A/P$ is finite over $R$, so $\dim R/(P\cap R)=\dim A/P$. If $P$ is not minimal over $J$, this dimension is at most one, and Lemma~\ref{lem:projection-types}(a) makes $A/P$ strongly type~(i). For a minimal factor $P=\mathfrak P/KS$, put $\mathfrak a=\mathfrak c_{\mathfrak P}/KR_0$. Then $R/\mathfrak a\cong R_0/\mathfrak c_{\mathfrak P}$ is Gorenstein of dimension two, and the localization equality above gives $\mathfrak aR_{P\cap R}=(P\cap R)R_{P\cap R}$. Lemma~\ref{lem:projection-types}(b) therefore makes $A/P$ of type~(ii).

Thus every factor is of type~(ii), since strongly type~(i) implies type~(ii). Extension closure makes $C$ of type~(ii), so $J$ is glicci by Lemma~\ref{lem:type-calculus}. Lifting its links through the complete intersection $KS$ proves that $I$ is glicci, homogeneously in the graded case.
\end{proof}

\section{Gorenstein linkage class of monomial ideals}

Throughout this section, let $S=k[x_1,\ldots,x_n]$,
where $k$ is an infinite field, and let $I\subseteq S$
be a Cohen-Macaulay monomial ideal of height at least
three. In the first subsection, we show that height-three Cohen-Macaulay monomial ideals are 2-glicci. In the second subsection, we prove that every Cohen-Macaulay monomial ideal is glicci.
These results naturally lead to the following question.

\begin{question}\label{Questionmonomial2glicci}
    Let $S$ be a polynomial ring, with an infinite residue field, and $I$ be a Cohen-Macaulay monomial ideal. Is $I$ homogeneously 2-glicci?
\end{question}

\subsection{Cohen-Macaulay monomial ideals of height 3 are 2-glicci}

Again, let $S=k[x_1,..,x_n]$ and $I$ be height 3 Cohen-Macaulay proper monomial ideal of $S$. 
\begin{thm}\label{height3monomialis2glicci}
    Let $S$ be a polynomial ring, with an infinite residue field. Let $I$ be a height three Cohen-Macaulay monomial ideal. Then $I$ is 2-glicci. More precisely, for any homogeneous $f\in I$, $I/(f)$ is homogeneously glicci in $S/(f)$. In particular, $I$ is glicci.
\end{thm}
We break the proof into steps. As $S$ is a regular ring, we get that any height 1 Gorenstein ideal is just complete intersection of height 1. So for any homogeneous $f\in I$, we show that $I/(f)$ is glicci in $S/(f)$. If $\deg f=1$, then $I/(f)$ is perfect ideal of height 2 in $S/(f)$ hence licci; So we assume $\deg f\ge 2$. As $S/(f)$ is a hypersurface ring, we have the structure theorem due to Eisenbud for MCM modules over hypersurfaces. 

\begin{thm}\label{eisenbud}\cite{eisenbud1980homological},\cite[Lemma~1.7]{herzog1987maximal}
Let $M$ be an graded MCM-module over
an graded hypersurface ring R. Let $p=\text{free-rank}$ $M$, and write: $M\cong U\oplus F$, where $F$ is graded free module and $\rank F=p$, 
then $\Omega^2(M)\cong U(e)$ for some $e\in\mathbb Z$.In particular, if $p=0$, then $\Omega^2(M)\cong M(e)$ for some $e\in\mathbb Z$.
\end{thm}

\begin{remark}\label{primefiltrationformonomialideals}
    Let $I$ be a non zero monomial ideal in $S$, there is a finite filtration $0=M_0\seq M_{1}\seq\cdots\seq M_{r-1}\seq M_r=S/I$ such that each successive quotients $M_i/M_{i+1}$ is $S/P_i(-a_i)$, where $P_i$ an ideal generated by some collection of variables.
\end{remark}
\begin{proof}
We argue by Noetherian induction on monomial ideals, ordered by reverse inclusion. If \(I=S\), there is nothing to prove; assume every monomial ideal properly containing \(I\) has the desired filtration. If \(I\neq S\), choose \(P\in\Ass_S(S/I)\); since \(I\) is monomial, \(P=I:u\) for some monomial \(u\notin I\), and \(P\) is generated by variables. Then
\[
0\longrightarrow S/P(-\deg u)\xrightarrow{\cdot u} S/I\longrightarrow S/(I+(u))\longrightarrow0
\]
is exact. Since \(I+(u)\supsetneq I\), the induction hypothesis gives the required filtration of \(S/(I+(u))\), and lifting this filtrations using the  above exact sequence gives the desired filtration of \(S/I\).
\end{proof}

Every $P_j$ contains $I$; hence, contains a minimal prime $Q_j$ of $I$.
Since $S/I$ is Cohen-Macaulay, these minimal primes have height three. Thus, after ordering the variables generating $P_j$,
\[
 I\subseteq Q_j=(y_1,y_2,y_3)\subseteq
 P_j=(y_1,\ldots,y_r),\qquad r\ge3.
\]
In particular, for any $f\in I$, the filtration is also an $R=S/(f)$-module filtration, with factors
$(R/\overline P_j)(-a_j)$, where $\overline P_j=P_j/(f)$.

So for the rest of this subsection, fix; $0\neq f\in I$ is a homogeneous element; set $R=S/(f)$. Also, note that $\dim R=n-1$ and $\overline{I}$ is a height 2 Cohen-Macaulay ideal of $R$. Again, let $\mathcal{P}_{i}(R)$ denote the stable isomorphism class of doubly layered modules that have a finite filtration, with successive quotients of type~(i).

\begin{lem}\label{lem:filtration}We stick to the same notations and let $q\ge n-1$. If every $\Omega_R^q(R/\overline P_j)$ is stably
isomorphic to a module filtered by type~(i) factors,  where $P_i$'s are the primes appearing in the filtration above. Then so is
$\Omega_R^q(R/\overline{I})$.
\end{lem}
\begin{proof}
The graded horseshoe lemma gives an exact sequence
\[
0\longrightarrow\Omega_R^qM_{j-1}
\longrightarrow\Omega_R^qM_j\oplus F_j
\longrightarrow\Omega_R^q(R/\overline P_j)(-a_j)\longrightarrow0
\] where $F_j$ is a graded free $R$-module. All these syzygies are MCM and orientable because $q\ge\dim R$. Now, the rest follows from induction on $j$, using Lemma~\ref{lem:closure} or Lemma~\ref{shortexactDL}.
\end{proof}

\begin{thm} Let $R=k[x_1,..,x_n]/(f)$, and if $f\in(x_{i_1},x_{i_2},x_{i_3})\seq(x_{i_1},..,x_{i_r})$. Then for any $r\ge 3$, if $n$ is even, then  $\Omega^{n-1}_R(R/(x_{i_1},..,x_{i_r}))$ $($if $n$ is odd then $\Omega^n_R(R/(x_{i_1},..,x_{i_r})))$ is stably isomorphic to doubly layered module, that has a finite filtration such that every successive quotient is of type~(i) i.e., $\Omega_R^{n-1}(R/(x_{i_1},..,x_{i_r}))\in \mathcal{P}_{i}(R)$ if $n$ is even and $\Omega_R^{n}(R/(x_{i_1},..,x_{i_r})\in\mathcal{P}_i(R)$ when $n$ is odd. 
\end{thm}

We suppress the grading shifts and the bar notations $\overline{x_{i_j}}$ notations to keep the proof readable. 
   
\begin{proof} We prove the statement for $\Omega^{n-1}(R/(x_{i_1},..,x_{i_r}))$, the other case is similar. We prove this by induction on $r$: if $r=3$, from the Shamash resolution; we get $R$ betti numbers of $R/(x_{i_1},x_{i_2},x_{i_r})$ as $1,3,4,4,\cdots$. Hence, $\Omega_R^q(R/(x_{i_1},x_{i_2},x_{i_3}))$ is rank 2 orientable MCM $R$-module for any $q\ge 2$.

Set $T_r=k[x_1,..,x_n]/(x_{i_1},\cdots,x_{i_r})$. Now for $r>3$, multiplication by $x_{i_r}$ gives the graded exact sequence
\[
0\longrightarrow T_{r-1}(-1)\xrightarrow{\ \cdot x_{i_r}\ }T_{r-1}
\longrightarrow T_r\longrightarrow0.
\]
For each $q\ge r-1$, by Lemma~\ref{lem:syz}
\[
0\longrightarrow\Omega_R^qT_{r-1}
\longrightarrow\Omega_R^qT_r\oplus F_q
\longrightarrow\Omega_R^{q-1}T_{r-1}(-1)\longrightarrow0,
\]
where $F_q$ is graded free. Since $q-1\ge r-2$, both outer terms
are stably filtered by type~(i) factors by induction.
The extension results Lemma~\ref{lem:closure} or by Lemma~\ref{shortexactDL} we get same conclusion for $\Omega_R^qT_r$, completing the induction.
\end{proof}

\begin{proof}[Proof of Theorem~\ref{height3monomialis2glicci}] Set $R=S/(f)$ for some non zero homogeneous element $f\in I$. Consider an N-type resolution of $\overline{I}$ $$0\x F\x N\x \overline{I}\x 0.$$ Then $\Omega_R^3(R/\overline{I})$ is stably isomorphic up to shift to $\Omega^2_R(N)$, and this is stably isomorphic up to shift to $N$. Now $\Omega^3_R(R/\overline{I})$ is stably isomorphic to $\Omega_R^{n-1}(R/\overline{I})$ if $n$ is even, and $\Omega^n_R(R/\overline{I})$ if $n$ is odd. Now by the previous results, we get that $\Omega_R^{c}(R/\overline{I})$, $c\in\{n-1,n\}$ is stably isomorphic up to shift to a doubly layered module; hence, so is $N$. So, by Theorem~\ref{DLstablyGlicci}, we get $\overline{I}$ to be glicci.    
\end{proof}

\subsection{Cohen-Macaulay monomial ideals are glicci}

We stick to the same notations as in the previous subsections. The goal of this section is to prove that Cohen-Macaulay monomial ideals are glicci.

Recall that $S=k[x_1,\ldots,x_n]$ has a natural
$\mathbb Z^n$-grading, defined by
$\deg x_i=\mathbf e_i$, where $\mathbf e_i$ is
the $i$-th standard basis vector of $\mathbb Z^n$.
An $S$-module $M$ is called \emph{multigraded} if $M=\bigoplus_{\mathbf a\in\mathbb Z^n}M_{\mathbf a}$ such that $x_iM_{\mathbf a}\subseteq M_{\mathbf a+\mathbf e_i}$.
An ideal $I\subseteq S$ is monomial if and only
if $S/I$ is multigraded. A homomorphism
$\varphi:M\to N$ of multigraded modules is called
\emph{multigraded} if it preserves multidegrees,
that is,
$\varphi(M_{\mathbf a})\subseteq N_{\mathbf a}$
for every $\mathbf a\in\mathbb Z^n$.
In particular, every multigraded module and
homomorphism is graded with respect to the
usual total degree.

\begin{thm}\label{thm:main}
Let $S=k[x_1,\ldots,x_n]$ be a standard graded polynomial ring, where $k$ is an infinite field,
and let $I\subset S$ be a proper monomial ideal such that $S/I$ is
Cohen-Macaulay and $\height I=g\geq3$. There is a homogeneous
$S$-regular sequence $f_1,\ldots,f_{g-2}\in I$ such that, on setting
\[
 K=(f_1,\ldots,f_{g-2}),\qquad R=S/K,\qquad J=I/K,
\]
the ideal $J$ is homogeneously glicci in $R$. In particular, $I$ is
homogeneously glicci in $S$.
\end{thm}

We break the proof into multiple steps to keep it readable. We stick to the same notations; a \emph{type-(i)} module is a
nonfree rank-two orientable MCM module. A \emph{type-(ii)} module has
the form \((\Omega_A(E))^*\)
for a type~(i) module $E$. 
Let $\mathcal{P}_{i}(A)$, and $\mathcal{P}_{ii}(A)$ consist of modules that are up to shift stably isomorphic to modules with finite filtrations entirely of type (i) and entirely of type (ii), respectively. The zero module, with its empty
filtration, is allowed. Thus, free modules belong to both stable classes. Denote $M\cong_{st} N$ if $M$ and $N$ are stably isomorphic up to shift. 

The following lemma tells us how N-type resolutions of directly G-linked Gorenstein ideals behave. This is an important tool for proving the main theorem. Before that lets note the following remark which is used in the proof of the lemma. 

\begin{remark}
Let $A$ be a graded Gorenstein ring, and let $p$ be a homogeneous
height-two Gorenstein ideal. Consider an N-type resolution $0\x A(a)\xrightarrow{\iota}E
\xrightarrow{\pi}p\x 0$, then $E^*\cong E(-a)$.
\end{remark}
\begin{proof}
Dualizing with respect to $A(a)$ gives
\[
0\to A(a)\to E^*(a)\xrightarrow{\rho}A
\to\operatorname{Ext}_A^1(p,A(a))\to0.
\]
Since $p$ is Gorenstein ideal of height two,
$\Ext_A^1(p,A(a))\cong(A/p)(b)$
for some $b\in\mathbb Z$. Hence, $\operatorname{im}\rho=p$. So, we get an N-type resolution $0\x A(a)\x E^*(a)\x p\x 0$. Let $V$ be the pullback of $E\to p$ and $E^*(a)\to p$.
Then we have exact sequences
\[
0\to A(a)\to V\to E\to 0,
\qquad
0\to A(a)\to V\to E^*(a)\to 0.
\]
Since $E$ and $E^*(a)$ are MCM over the Gorenstein ring $A$,
both sequences split in the graded category. Hence $E\oplus A(a)\cong V\cong E^*(a)\oplus A(a)$. By graded free cancellation, $E\cong E^*(a)$,
and therefore $E^*\cong E(-a)$.
\end{proof}

\begin{lem}\label{lem:propagation}
Let $A$ be a standard graded Gorenstein ring  of dimension at least two.
Let $q,p,p'$ be homogeneous height-two Gorenstein ideals with
$q:p=p'$ and $q:p'=p$. Choose rank-two $N$-type resolutions
\begin{align*}
 0&\x A(a_0)\x E_0\x q\x0,\\
 0&\x A(a)\x E\x p\x0,\\
 0&\x A(a')\x E'\x p'\x0.
\end{align*}
For every $t\geq1$ there are graded free modules $H_t,H'_t$ and
exact sequences
\begin{align}
 0&\x\Omega_A^tE_0\x\Omega_A^tE\oplus H_t
 \x\Omega_A^{t-1}(E'^*)(a_0)\x0,
 \label{eq:prop1}\\
 0&\x\Omega_A^tE_0\x\Omega_A^tE'\oplus H'_t
 \x\Omega_A^{t-1}(E^*)(a_0)\x0.
 \label{eq:prop2}
\end{align}
Consequently, for any integer $m\geq0$,
\begin{equation}\label{eq:prop-conclusion}
 \bigl[\Omega_A^tE_0\in\mathcal{P}_{i}(A)\text{ for }1\leq t\leq m\bigr]
 \ \Longrightarrow\
 \bigl[\Omega_A^tE,\Omega_A^tE'\in\mathcal{P}_{i}(A)
       \text{ for }0\leq t\leq m\bigr].
\end{equation}
\end{lem}

\begin{proof}
Set $C=A/q$ and $L=p/q$, now these are Cohen-Macaulay of dimension
$\dim A-2$. The inclusion $q\subset p$ lifts to $E_0\to E$, and since
$\Ext_A^1(E_0,A(a))=0$. The resulting cone is exact (Remark \ref{mappingcone}):\[
 0\x A(a_0)\x E_0\oplus A(a)\x E\x L\x0.\]
All its terms before $L$ are MCM $A$-modules, and as $L$ is Cohen-Macaulay of codimension two, dualizing with respect to $A$ gives an exact sequence
\[0\x E^*\x E_0^*\oplus A(-a)\x A(-a_0) \x\Ext_A^2(L,A)\x0\tag{$\ast$}\]

Dualizing the $N$-type sequence for $q$ shows that
$\Ext_A^2(C,A)\cong C(-a_0)$. The change-of-rings spectral
sequence, whose only nonzero row is in degree two, therefore gives
\[
 \Ext_A^2(L,A)\cong\h_C(L,C)(-a_0).
\]
Apply $\h_C(-,C)$ to $0\to L\to C\to A/p\to0$.
The module $A/p$ is MCM over the Gorenstein ring $C$, so
$\Ext_C^1(A/p,C)=0$. Moreover
$\h_C(A/p,C)=(q:p)/q=p'/q$. Consequently
\[
 \h_C(L,C)\cong C/(p'/q)=A/p'.
\]
It follows that the image of $A(-a_0)$ in ($\ast$) is $p'(-a_0)$. Thus $0\x E^*\x E_0^*\oplus A(-a)\x p'(-a_0)\x0$ is exact. Now, pull back the above along
$E'(-a_0)\to p'(-a_0)$. The pullback $V$ fits into
\[
 0\x A(a'-a_0)\x V\x E_0^*\oplus A(-a)\x0.
\]
The dual of this exact sequence splits. 
The other pullback sequence is
$0\to E^*\to V\to E'(-a_0)\to0$.
Dualizing this sequence and using the splitting gives
\[
 0\x E'^*(a_0)\xrightarrow{\iota}
 E_0\oplus A(a)\oplus A(a_0-a')\xrightarrow[]{\rho} E\x0.
\]

Let $\pi_0:P_0\twoheadrightarrow E_0$ be a graded free cover; set
$\pi=\pi_0\oplus\mathrm{id}:P_0\oplus A(a)\oplus
A(a_0-a')\twoheadrightarrow E_0\oplus A(a)\oplus A(a_0-a')$ and let $K=\ker(\rho\pi)$.
The map $K\to E'^*(a_0)$, $z\mapsto\iota^{-1}(\pi(z))$, is surjective
with kernel $\ker\pi=\Omega E_0$. As $P_0\oplus A(a)\oplus
A(a_0-a')\twoheadrightarrow E$ is a free presentation,
$K\cong\Omega E\oplus H_1$ for some graded free $H_1$.
Thus $0\to\Omega E_0\to\Omega E\oplus H_1\to E'^*(a_0)\to0$,
which is $(1)$ for $t=1$.
Taking successive syzygies via the horseshoe lemma gives
$(1)$ for all $t\geq1$; interchanging
$(p,E,a)$ and $(p',E',a')$ gives $(2)$.

Finally, $E,E'$ belongs to $\mathcal{P}_{i}(A)$; as each is either type~(i) or free. Now, by the remark above, we have 
$E^*\cong E(-a)$ and $E'^*\cong E'(-a')$.
For $t=1$, both right terms in $(1)$ and $(2)$
are in $\mathcal{P}_{i}(A)$, and $E'^*$, $E^*$ are in $\mathcal{P}_{i}$; hence, by Lemma~\ref{shortexactDL}, we get the required result.
Inductively, for $t>1$, they are in $\mathcal{P}_{i}(A)$ by the conclusions
already obtained at $t-1$.
Extension closure from Lemma~\ref{lem:closure} proves $(3)$.
\end{proof}

We set $\operatorname{NF}_A(L)=\{\mathfrak p\in\Spec A\mid L_\mathfrak p$ is not a free $A_\mathfrak p$-module\}. The following lemma shows that, given a finite collection of orientable MCM modules, one can choose a regular element \(f\) such that all these modules remain orientable and MCM after reduction modulo \(f\). %

The underlying tool for the lemma is essentially prime avoidance. 

\begin{lem}\label{lem:specialize}
Let $A$ be a standard graded complete intersection of dimension at
least three, and let $\mathcal L$ be a finite collection of graded orientable
MCM modules. Let $f\in A_+$ be homogeneous and $A$-regular. Suppose
$f$ avoids every height-two prime in $\operatorname{NF}_A(L)$ for every
$L\in\mathcal L$. Set $B=A/(f)$. Then each $L/fL$ is orientable
MCM over $B$, with the same rank as $L$.

In particular, fix finitely many modules which up to stability has type~(i) filtrations over $A$, and include all of their factors in $\mathcal L$. Reduction modulo
$f$ carries each of those stable classes into $\mathcal{P}_{i}(B)$.
\end{lem}

\begin{proof} Define $V=\{\mathfrak p\in\Spec A\mid L_\mathfrak p$ is nonfree for some $L\in\mathcal{L}\}$, this is a closed set, since the nonfree locus of a finite module is a closed set. As, $L$ is locally
free at primes of height at most one, height-two primes in $V$ (if any)
are minimal primes of $V$. 
The element $f$ is regular on every MCM module hence, $L/fL$ is MCM over $B$.

If $\overline{\mathfrak p}\in\Spec B$ has height at most one and $\mathfrak p$ is its
inverse image in $A$, then $\height \mathfrak p\leq2$.
If $\height \mathfrak p\leq1$, $L_\mathfrak p$ is free. If $\height\mathfrak p=2$, the
avoidance hypothesis and $f\in\mathfrak p$ again force $L_\mathfrak p$ to be free. Thus $L/fL$
is locally free of the same rank in depth one over $B$.

Let $(\bigwedge^rL)^{**}\cong A(a)$; set $M=\bigwedge_B^r(L/fL)$ and $\delta$ be the composition of $\bigwedge_A^rL\x (\bigwedge_A^rL)^{**}\x A(a)$; define $\phi:M\x B(a)$ by $\phi(\overline{\ell_1}\wedge\cdots\wedge\overline{\ell_r})=\overline{\delta(\ell_1\wedge\cdots\wedge\ell_r)}$.
The established freeness of $L_{\mathfrak p}$ implies that $\phi_{\mathfrak q}$ is an isomorphism whenever $\operatorname{ht}_B\mathfrak q\le1$.
Thus $\ker\phi$ and $\operatorname{coker}\phi$ have grade at least two; dualizing the two short exact sequences associated to $\phi$ gives $\phi^*:B(-a)\xrightarrow{\sim}M^*$.
Dualizing again yields $\bigl(\bigwedge_B^r(L/fL)\bigr)^{**}\cong B(a)$, proving orientability.

Every short exact sequence in a filtration remains exact after reduction, because its quotient is MCM and hence $f$ is regular on it. Stable isomorphisms also reduce to stable isomorphisms. Each reduced rank-two factor is either type (i) or free, and
Lemma~\ref{lem:closure}(a) handles both possibilities.
\end{proof}

Note that for every MCM $A$-module $E$ and every $t\geq0$ one also has $$
 \Omega_B^t(E/fE)\cong(\Omega_A^tE)/f(\Omega_A^tE).$$

The following lemma is an alternative strategy that is similar to the strategy used to prove that the clean monomial ideals are glicci. We replace the clean filtration with a monomial coresolution. 

\begin{lem}\label{lem:irreducible}
Let $S$ and $I$ be as in Theorem~\ref{thm:main}. There is a finite
exact complex of multigraded, hence also graded, $S$-modules
\[
 0\x S/I\x\mathcal W^0\x\mathcal W^1
 \x\cdots\x\mathcal W^\ell\x0,
 \qquad
 \mathcal W^i=\bigoplus_\lambda S/W_{i\lambda}\]
where each $W_{i\lambda}$ is generated by positive powers of
distinct variables. If $S/W_{i\lambda}$ occurs in $\mathcal W^i$, then
\[
 g\leq\height W_{i\lambda}\leq g+i.
\]
\end{lem}

The proof follows from the techniques of \cite{miller2002cohen} and \cite[Proposition~2.1 and Theorem~2.2]{yanagawa2008notes}. Though they present a very general proof, for the sake of completeness, we provide a short proof in the particular setting that we need. 
\begin{proof}
    Write $I=\bigcap_{i=1}^t Q_i$ where $Q_i$'s are the irreducible monomial ideals of the form $(x_{i_1}^{a_{i_1}},...,x_{i_m}^{a_{i_m}})$. Write $[t]=\{1,\ldots,t\}$, for every nonempty subset
$\sigma\subseteq[t]$, set
\[
Q_\sigma=\sum_{j\in\sigma}Q_j,
\qquad
\mathcal{W}^i=\bigoplus_{\substack{\sigma\subseteq[t]\\|\sigma|=i+1}}
S/Q_\sigma
\quad(0\leq i\leq t-1).
\]
Each $Q_\sigma$ is again generated by powers of distinct variables.
For $\sigma\subseteq\tau$, denote the natural quotient map by $
\pi_{\sigma,\tau}:S/Q_\sigma\longrightarrow S/Q_\tau$. Now define the augmentation
$\varepsilon:S/I\x W^0$ by
$\varepsilon(f+I)=(f+Q_j)_{j=1}^t$.
Define $d^i:\mathcal{W}^i\to \mathcal{W}^{i+1}$ by specifying its components:
\[
d^i_{\tau,\sigma}=
\begin{cases}
(-1)^{|\{j\in\sigma:j<a\}|}\pi_{\sigma,\tau},
  &\tau=\sigma\cup\{a\},\ a\notin\sigma,\\
0,&\text{otherwise}.
\end{cases}
\]
From the definitions it is clear that the maps have degree (infact multidegree) 0 and $d^{i+1}d^{i}=0$ and $d^0\epsilon=0$. The exactness of the complex is standard, we leave it for the reader. 

We now modify the constructed complex so that every nonzero differential component $S/W\to S/V$ satisfies $W\subsetneq V$.

First, any multigraded degree-zero map of modules $\varphi:S/W\to S/V$
is determined by $\varphi(1+W)\in(S/V)_0=k$.
If $\varphi\ne0$, then $\varphi(1+W)=c(1+V)$ for some
$c\in k^\times$, and well-definedness implies $W\subseteq V$.
Thus $\varphi$ is a nonzero scalar multiple of the natural
quotient map. In particular, if $W=V$, then $\varphi$ is
an isomorphism.

Suppose that $d^i$ has such an isomorphism component.
Let
$M=S/W$ and separate the corresponding summands: $
\mathcal W^i=M\oplus A$, $
\mathcal W^{i+1}=M\oplus B$, $d^i=
\begin{pmatrix}
u&\alpha\\
\beta&\gamma
\end{pmatrix}$, where $u:M\to M$ is an isomorphism,
$\alpha:A\to M$, $\beta:M\to B$, and $\gamma:A\to B$.
Define homogeneous automorphisms
\[
P=
\begin{pmatrix}
1&-u^{-1}\alpha\\
0&1
\end{pmatrix}
\text{ of }M\oplus A,
Q=
\begin{pmatrix}
1&0\\
-\beta u^{-1}&1
\end{pmatrix}
\text{ of }M\oplus B\text{ then }
\widetilde d^{\,i}:=Qd^iP
=
\begin{pmatrix}
u&0\\
0&\gamma-\beta u^{-1}\alpha
\end{pmatrix}.
\]
At the neighboring positions, replace the maps by $\widetilde d^{\,i-1}=P^{-1}d^{i-1}$, $\widetilde d^{\,i+1}=d^{i+1}Q^{-1}$ (set $d^{-1}=\epsilon$).
These coordinate changes give an exact isomorphic complex. Write
$\widetilde d^{\,i-1}=\binom{p}{r}$ and
$\widetilde d^{\,i+1}=(q\ \ s)$ according to the displayed
direct-sum decompositions. Now, the complex identities give, $up=0$, $qu=0$. Since $u$ is an isomorphism,  $p=q=0$.
Thus the preceding differential has no component into
the first copy of $M$, and the following differential
vanishes on the second copy of $M$.

Consequently, $
0\x M\xrightarrow{u}M\x0$ 
is a exact direct summand of the entire complex. So, deleting it, and changing the maps as above leaves an exact complex. Now, repeat this operation whenever the current complex has
a nonzero component between identical quotients, as there are only finitely many components, the process terminates.
Retain the notation $\mathcal W^\bullet$ for the resulting
finite exact coresolution of $S/I$, denote its differentials
by $\delta^i$, and set $\delta^\ell=0$. Note that we still have $$\mathcal{W}^i=\bigoplus_{\substack{\sigma\subseteq[t]\\|\sigma|=i+1}}
S/Q_\sigma
\quad(0\leq i\leq t-1)$$ and it has all the same properties as the original one.

Bounds on height: Since $I\seq W_{i\lambda}$, we get the lower bound. For, the upper bound let us observe the following; set $Z_j=\ker\delta^j$, with $\delta^\ell=0$, and $Z_0=S/I$, and fix a summand $S/W$ of $\mathcal W^i$.
After renumbering variables, write $W=(x_1^{a_1},\ldots,x_h^{a_h})$, and let $P=\sqrt{W}=(x_1,\ldots,x_h)$ and $m=\prod_{r=1}^h x_r^{a_r-1}$.
Since $(W:m)=P$, the element $\eta=(m+W)/1\in(S/W)_P$ is nonzero and annihilated by $PS_P$.
View $\eta$ as an element of $(\mathcal W^i)_P$, now every nonzero outgoing component $S/W\to S/V$ satisfies $W\subsetneq V$; if $V\nsubseteq P$, then $(S/V)_P=0$.
If $V\subseteq P$, its pure-power form and strict containment imply $x_r^{b_r}\in V$ for some $b_r<a_r$, so $m\in V$.
Thus every outgoing component kills $\eta$, and hence $\eta\in\ker(\delta^i_P)=(Z_i)_P$.
Since this kernel contains a nonzero element annihilated by $PS_P$, we obtain $\operatorname{depth}_{S_P}(Z_i)_P=0$.

Now we prove the upper bound by induction on $i$; We have the short exact sequence $0\x Z_i\x \mathcal{W}^i\x Z_{i+1}\x 0$. Fix a summand $S/W$ of $\mathcal W^i$ and set $P=\sqrt W$, $h=\height_S W$, $e=h-g$. Now $i=0$ is immediate from the definition, so $i>0$. For $j<i$, the inductive hypothesis gives
$\operatorname{ht}_S V\leq g+j$ for every summand $S/V$
of $\mathcal W^j$.
If $V\subseteq P$, then $\dep_{S_P}(S/V)_P
=h-\operatorname{ht}_S V\geq e-j$. If $V\nsubseteq P$, the localization is zero, whose depth
we interpret as $+\infty$.
It follows that
$\operatorname{depth}_{S_P}(\mathcal W^j)_P\geq e-j$. Applying the depth lemma successively to the localization
of the above short exact sequence yields
\[
\begin{aligned}
\operatorname{depth}_{S_P}(Z_{j+1})_P
&\geq
\min\bigl\{
\operatorname{depth}_{S_P}(\mathcal W^j)_P,\,
\operatorname{depth}_{S_P}(Z_j)_P-1
\bigr\}\\
&\geq e-j-1.
\end{aligned}
\]
Now from the above paragraph, we get $0=\dep_{S_P}(Z_i)_P\ge e-i=h-g-i$. Hence, $h\le g+i$.
\end{proof}

Let \(\sW=\{Q_\sigma\mid \sigma\seq [t]\}\) be the finite set of distinct ideals occurring in the resolution above. For every \(W\in\sW\), fix an ordering
of its generators and set $
 W=(y_{W,1},\ldots,y_{W,h_W})$, $
 y_{W,r}=z_{W,r}^{a_{W,r}}$, $h_W=\height_S W$, and for $1\le r\le h_W$, set $Q_{W,r}=(y_{W,1},\ldots,y_{W,r})$,  $P_{W,r}=\sqrt{Q_{W,r}}=(z_{W,1},\ldots,z_{W,r})$. Here the \(z_{W,r}\) are distinct variables for a fixed \(W\), and
\(a_{W,r}\geq1\). Different \(W\)'s need not have compatible orderings.
Since \(h_W\geq g\), every prefix used below is defined. 

The following theorem provides a way to choose a complete intersection \(K\subseteq I\), which is a key step in proving our main result. In the height-three case, after reducing modulo a non-zero element, the Betti numbers of an ideal generated by three variables allow us to obtain rank-two syzygies. However, this approach does not extend to the general case. To overcome this difficulty, we use part (c) of the following theorem. 
\begin{thm}\label{thm:construction}
We keep the same notations as Theorem~\ref{thm:main} and carry the notations $\mathfrak W$, $Q_{W,j}$, $P_{W,j}$ from the previous proof. Set \(c=g-2\).
There exist homogeneous \(f_1,\ldots,f_c\in S\) such that, for
\[
 K_j=(f_1,\ldots,f_j)\subset S,\qquad R_j=S/K_j
 \quad(0\leq j\leq c),\qquad K_0=(0),
\]
the following hold.
\begin{enumerate}
\item[$(a)$] The sequence \(f_1,\ldots,f_c\) is \(S\)-regular and lies in \(I\).
\item[$(b)$] For every \(W\in\sW\) and \(1\leq j\leq c\),
\[
 f_j\in\m Q_{W,j+2},\qquad f_j\notin P_{W,j+1}.
\]
\item[$(c)$] For every \(W\in\sW\) and \(0\leq j\leq c\),
\(J_{W,j}:=Q_{W,j+2}R_j\) is a height-two Gorenstein ideal and has
an exact sequence
\[
 0\x R_j(\beta_{W,j})\xrightarrow{\iota_{W,j}}E_{W,j}
 \xrightarrow{\nu_{W,j}}J_{W,j}\x0
\]
with \(E_{W,j}\) rank-two orientable MCM and $ \Omega_{R_j}^tE_{W,j}\in\mathcal{P}_{i}(R_j)$, for all $0\leq t\leq\ell$.
\end{enumerate}
\end{thm}
We break the proof into steps to keep it readable. 
\begin{proof}
We induct on the construction stage \(j\). The index \(i\) remains
reserved for the degree in the coresolution, and \(t\) denotes a
syzygy order.

For $j=0:$ \(W\in\sW\), \(J_{W,0}:=Q_{W,2}\subset S\) is a complete intersection.
Its Koszul sequence is
\[
 0\x S(-a_{W,1}-a_{W,2})
 \xrightarrow{\binom{-y_{W,2}}{y_{W,1}}}
 S(-a_{W,1})\oplus S(-a_{W,2})
 \xrightarrow{(y_{W,1}\ \ y_{W,2})}Q_{W,2}\x0.
\]
Choose this 
middle module as \(E_{W,0}\) and put
\(\beta_{W,0}=-a_{W,1}-a_{W,2}\).
It is free, and its positive minimal syzygies are zero. Hence holds for every \(0\leq t\leq\ell\).

\smallskip\noindent\textbf{1. Choosing the regular sequence.}
Let \(1\leq j\leq c\), and assume all assertions through stage \(j-1\).
Set
\[
 A=R_{j-1},\qquad \pi_{j-1}:S\twoheadrightarrow A,\qquad
 \mathfrak a_j=\m\left(I\cap\bigcap_{W\in\sW}Q_{W,j+2}\right)\subset S.
\]
The ideal \(\mathfrak a_j\) is not contained in any prime of \(S\) of height
at most \(j+1\). 
Now, for each \(W\in\sW\) and \(1\leq t\leq\ell\), fix a finite
type-(i) filtration representing the stable class of
\(\Omega_A^tE_{W,j-1}\). Let \(\mathcal L_j\) consist of every factor
of these filtrations together with every \(E_{W,j-1}\).
Now, for \(L\in\mathcal L_j\), let
\(\NF_A(L)=\{\mathfrak p\in\operatorname{Spec}A:
L_{\mathfrak p}\text{ is not free}\}\).
This is a closed set, and it does not contain any prime ideal of height 1 in $A$. Hence, there are finitely many height 2 primes of $A$ in it, and these are minimal primes so homogeneous. Now, choose \(\overline{f}_j\in\mathfrak a_j\) outside the following finite collection
of homogeneous primes of \(S\):
\begin{enumerate}
\item[$(i)$] the inverse images under \(\pi_{j-1}\) of \(\Ass A\);
\item[$(ii)$] the primes \(P_{W,j+1}\), for \(W\in\sW\);
\item[$(iii)$] the inverse images of the height-two primes in \(\NF_A(L)\),
for \(L\in\mathcal L_j\).
\end{enumerate}

Since $k$ is infinite, homogeneous prime avoidance yields a
homogeneous element $\overline{f}_j\in \mathfrak a_j$ which avoids all those primes. Set \(f_j=\pi_{j-1}(\overline{f}_j)\in A\) and
\(B=A/(f_j)=R_j\). This proves (a) and (b) at stage \(j\).

\smallskip\noindent\textbf{2. Preserving the chosen layers over \(B\).} By the way, we have chosen $f_j$, and using Lemma~\ref{lem:specialize} we get that for every $L\in\mathcal{L}_j$, $L/f_jL$ is an orientable MCM $B=A/(f_j)$-module. In particular if $L\in\mathcal{P}_i(A)$, then $L/f_jL\in\mathcal{P}_i(B)$.

Now, set \(\overline E_W=E_{W,j-1}/f_jE_{W,j-1}\).
The minimal graded free resolution of \(E_{W,j-1}\), after tensoring with $B$, becomes a minimal graded free resolution for $\overline{E_W}$. Therefore
\[
 \Omega_B^t\overline E_W
 \cong(\Omega_A^tE_{W,j-1})/f_j(\Omega_A^tE_{W,j-1})
 \in\mathcal{P}_{(i)}(B)\qquad(1\leq t\leq\ell).
\]
Also, note that \(\overline E_W\) is a rank-two orientable MCM.

\smallskip\noindent\textbf{3. Constructing the link for $J_{W,j}$.}
Write $f=\pi_{j-1}(f_j)\in A$, so that $B=A/(f)$.
Fix $W\in\mathfrak W$, and let $Q=Q_{W,j+1}$, $y=y_{W,j+2}$, and $T=S/Q$. Note from the construction of $Q_{W,j+1}$ we get $T$ to be a Gorenstein ring, and $T\cong A/QA$. 

For $r<j$, we have $f_r\in Q_{W,r+2}\subseteq Q$; hence $K_{j-1}\subseteq Q$. Since $f_j\in\mathfrak m(Q+(y))=\mathfrak mQ_{W.j+2}$, choose a homogeneous $b\in\mathfrak m$ such that $f_j\equiv yb\pmod Q$.
Writing $\bar y,\bar b$ for their images in $T$, the choice $f_j\notin P_{W,j+1}$ and the equality $\operatorname{Ass}(T)=\{P_{W,j+1}/Q\}$ show that $\bar y\bar b$ is $T$-regular.
Consequently, both $\bar y$ and $\bar b$ are $T$-regular.

Define ideals of $B$ by
\[
H=QB,\qquad J=(Q+(y))B=J_{W,j},\qquad J'=(Q+(b))B.
\]
The congruence $f_j\equiv yb\pmod Q$ gives $B/H\cong T/(\bar y\bar b)$, $B/J\cong T/(\bar y)$, $B/J'\cong T/(\bar b)$.
Now, these quotients are complete intersections; hence, cancellation of the nonzerodivisors $\bar y,\bar b$ gives $(\bar y\bar b):_T\bar y=(\bar b)$ and $(\bar y\bar b):_T\bar b=(\bar y)$.
Through the identification $B/H\cong T/(\bar y\bar b)$, these equalities yield $H:_B J=J'$, $H:_B J'=J$.

Thus $H$ links the desired ideal $J_{W,j}$ to $J'$.

\smallskip\noindent\textbf{4. Applying the Lemma~\ref{lem:propagation} to finish (c).}
Now $J_{W,j-1}=Q_{W,j+1}A$.
As element $f$ is regular on $A$, the natural map $J_{W,j-1}/fJ_{W,j-1}\x Q_{W,j+1}B$ is an isomorphism.

So tensoring the $N$-type resolution $0\x R_{j-1}(\beta_{W,j-1})\x E_{W,j-1}\x J_{W,j-1}\x0$ with $B$ remains exact and gives an $N$-type resolution
\[
0\longrightarrow B(\beta_{W,j-1})
\longrightarrow \overline E_W
\longrightarrow Q_{W,j+1}B
\longrightarrow0,
\qquad
\overline E_W=E_{W,j-1}/fE_{W,j-1}.
\]

By the rank-two $N$-type resolutions for Gorenstein ideals, choose exact sequences
\[
0\longrightarrow B(\beta)\longrightarrow E\longrightarrow J_{W,j}\longrightarrow0,
\qquad
0\longrightarrow B(\beta')\longrightarrow E'\longrightarrow J'\longrightarrow0,
\]
where $E,E'$ are rank-two orientable MCM modules and $J'$ as in step~3.
Step~2 gives $\Omega_B^t\overline E_W\in\mathcal P_{i}(B)$ for $1\le t\le\ell$.
Applying the Lemma~\ref{lem:propagation} to the link $J\sim_H J'$, with linking middle module $\overline E_W$, therefore gives
\[
\Omega_B^tE,\ \Omega_B^tE'\in\mathcal P_{(i)}(B)
\qquad(0\le t\le\ell).
\]
Set $E_{W,j}=E$ and $\beta_{W,j}=\beta$. 
This proves the required assertion for $W$.
Since the same $f_j$ was chosen simultaneously for every $W\in\mathfrak W$, the induction step is complete.

\smallskip\noindent\textbf{5. Completing the induction on \(j\).}
The same \(f_j\) was chosen for every \(W\in\sW\), and Step 4 supplies
the required \(E_{W,j}\) for every \(W\). Thus all assertions hold at
stage \(j\), completing the induction up to \(j=c=g-2\).
\end{proof}

From now on fix the equations of Proposition~\ref{thm:construction}
and set $R=R_c$, $K=(f_1,\ldots,f_c)$, $d=\dim R=n-g+2$. The ring $R$ is standard graded Gorenstein and $d\geq2$.
For $W\in\mathcal W$, set $E_W=E_{W,c}$.
Since $c+2=g$, the construction gives
\begin{equation}\label{eq:final-E}
 0\x R(b_W)\x E_W\x Q_{W,g}R\x0,
 \qquad
 \Omega_R^tE_W\in\Pone(R)\quad(0\leq t\leq\ell).
\end{equation}
Also $K\subseteq I\cap\bigcap_{W\in\mathcal W}W$, so the
coresolution, is an exact complex of
graded $R$-modules.

\begin{prop}\label{prop:all-heights}
For $W\in\mathcal W$, let $h=\heightideal_S W$. Then
\[\label{eq:all-heights}
 (\Omega_R^q(S/W))^*\in\mathcal{P}_{ii}(R)
 \qquad\text{whenever }h-g+2\leq q\leq\ell+2.
\]
In particular, if $S/W$ occurs in $\mathcal W^i$, then
$(\Omega_R^{i+2}(S/W))^*\in\mathcal{P}_{ii}(R)$.
\end{prop}

\begin{proof}
Fix $W$, set $e=h-g$ and
\[
 T_0=R/Q_{W,g}R\cong S/Q_{W,g},\qquad
 u_r=z_{g+r}^{a_{g+r}}\ (1\leq r\leq e),\qquad
 T_r=T_0/(u_1,\ldots,u_r)T_0.
\]
The $u_r$ form a $T_0$-regular sequence because they are powers
of distinct variables not occurring in $Q_{W,g}$. 
Each $T_r$ is Cohen-Macaulay over $R$, with $\depth_R T_r=\dim T_r=d-r-2$.

We prove by induction on $r$ that
\begin{equation}\label{eq:Tr-goal}
 (\Omega_R^qT_r)^*\in\mathcal{P}_{ii}(R)
 \qquad(r+2\leq q\leq\ell+2).
\end{equation}
For $r=0$, \eqref{eq:final-E} and
Lemma~\ref{lem:syz}(c) give
\[
 (\Omega_R^qT_0)^*
 \st(\Omega_R^{q-1}E_W)^*
 =(\Omega_R(\Omega_R^{q-2}E_W))^*.
\]
Here $0\leq q-2\leq\ell$, so
$\Omega_R^{q-2}E_W\in\Pone(R)$ by
\eqref{eq:final-E}. Lemma~\ref{lem:closure}(c) proves the claim.

For $r\geq1$, apply Lemma~\ref{lem:syz} to $0\x T_{r-1}\xrightarrow{\cdot u_r}T_{r-1}\x T_r\x 0$, we get for all $q\geq r+2$
\[
 0\x\Omega_R^qT_{r-1}
 \x\Omega_R^qT_r\oplus H
 \x\Omega_R^{q-1}T_{r-1}(-\deg u_r)\x 0.
\]
All three syzygies are MCM; so, $0\x(\Omega_R^{q-1}T_{r-1})^*(\deg u_r)
 \x(\Omega_R^qT_r)^*\oplus H^*\x(\Omega_R^qT_{r-1})^*\x0$.
Both end terms belong to $\mathcal{P}_{ii}(R)$ by induction, since
$q-1\geq r+1=(r-1)+2$ and $q\leq\ell+2$.
Lemma~\ref{lem:closure}(a) proves \eqref{eq:Tr-goal}.
Take $r=e$ to obtain the required.

If $S/W$ occurs in $\mathcal W^i$, then
$h-g\leq i$ by the height bounds in the Theorem~\ref{thm:construction}. Therefore
$h-g+2\leq i+2\leq\ell+2$, proving the last assertion.
\end{proof}

Now we note the following, lemma, which is the last step need to finish the main theorem. Define $C_0=S/I$ and write the coresolution as
\[
 0\x C_i\x\mathcal W^i\x C_{i+1}\x0
 \quad(0\leq i<\ell),\qquad C_\ell\cong\mathcal W^\ell.\]
These are exact sequences of graded $R$-modules because every
term in the coresolution is annihilated by $K$.

\begin{lem}\label{lem:Ci-depth}
For $0\leq i\leq\ell$, both $\Omega_R^{i+2}C_i$ and
$\Omega_R^{i+2}\mathcal W^i$ is orientable MCM $R$-module.
\end{lem}

\begin{proof}
For a summand $S/W$ of $\mathcal W^i$, with $h=\heightideal_S W$, by the height bound in Theorem~\ref{thm:construction} we get
\[
 \depth_R(S/W)=n-h\geq n-g-i=d-i-2.
\]
 Hence
$\depth_R\mathcal W^i\geq d-i-2$.
Also $\depth_R C_0=d-2$; hence, by the depth lemma, we get
\[
 \depth_R C_{i+1}\geq
 \min\{\depth_R\mathcal W^i,\depth_R C_i-1\}
 \geq d-i-3.
\]
Thus $\depth_R C_i\geq d-i-2$ for all $i$.
Taking $i+2$ syzygies, we get MCM $R$-module. As the ideals involved in the definition of $C_0$ and $\mathcal{W}^i$ have height at least two, by induction, we obtain the orientability.
\end{proof}

Now we are set to finish the proof of Theorem~\ref{thm:main}.

\begin{proof}[Proof of Theorem~\ref{thm:main}]
Set $X_i=(\Omega_R^{i+2}C_i)^*$, $Y_i=(\Omega_R^{i+2}\mathcal W^i)^*$.
By Proposition~\ref{prop:all-heights} and finite direct sums, we get $Y_i\in\mathcal{P}_{ii}(R)$ for all $0\leq i\leq\ell$.

Apply Lemma~\ref{lem:syz}(a) to $0\x C_i\x\mathcal{W}^i\x C_{i+1}\x0$, with $q=i+2$.
It yields
\[
 0\x\Omega_R^{i+3}C_{i+1}
 \x\Omega_R^{i+2}C_i\oplus G_i
 \x\Omega_R^{i+2}\mathcal W^i\x0
\]
for a graded free $G_i$.
Every term is MCM by Lemma~\ref{lem:Ci-depth}, so the dual
sequence is exact:
\[0\x Y_i\x X_i\oplus G_i^*\x X_{i+1}\x0.\]
At the last term, $X_\ell=Y_\ell\in\mathcal{P}_{ii}(R)$.
Descending induction on the exact sequence above, using
Lemma~\ref{lem:closure}, we conclude that $(\Omega_R^2(S/I))^*=X_0\in\mathcal{P}_{(ii)}(R)$. 

Now consider an $N$ type resolution of $J=I/K\seq R=S/K$, $0\x F\x M\x J\x 0$. Then $\Omega_R(M)^*\cong_{st}(\Omega^2_R(S/I))^*$; hence, $\Omega(M)^*\in\mathcal{P}_{ii}(R)$. Hence by Lemma~\ref{lem:closure} we get $(\Omega(\Omega_RM)^*)^*\cong_{st} M\in \mathcal{P}_{i}(R)$. So, by Theorem~\ref{DLstablyGlicci} we get $J$ to be glicci. 
\end{proof}

\section{Reduction of generators in the Gorenstein linkage classes}

The following proposition shows that almost complete
intersection ideals are glicci. On the other hand, we
have seen examples of Cohen-Macaulay ideals of heights
one and two that are not glicci. This motivates the
question of determining the minimum number of generators
among ideals in a given Gorenstein linkage class.
More precisely, let $R$ be a Gorenstein local ring
and $I$ a Cohen-Macaulay ideal. Set
\[
\mathfrak{S}_R(I)=
\{\mu(J)\mid J \text{ belongs to the Gorenstein
linkage class of } I\}.
\]
In this section, we show that if $\operatorname{ht} I\geq 3$,
then $\operatorname{ht} I\leq\inf\mathfrak{S}_R(I)
\leq\operatorname{ht} I+2$. The lower bound is immediate. However, for ideals of
heights one and two over nonregular rings, the
situation is less clear.

\begin{ex}
    Let $R$ be a hypersurface local ring, with an isolated singularity and $\dim R\ge 8$. Now, by Corollary~\ref{herzogrankisolatedsing} for a nonfree MCM $R$-module $M$, we have $2\rank M+1\ge\dim R$. So all rank three MCM $R$-modules are free. Now, note that every Cohen-Macaulay $R$-ideal $J$, with $\operatorname{ht}_RJ=2$ and $\mu(J)\le 4$ is perfect hence licci; This is because $\rank\Omega^2(R/J)\le 3$. 

    So non-glicci height two $R$-ideal cannot be in the Gorenstein linkage class of 4 generated ideal. Example~\ref{isolatedsingularityexample} shows the existence of non-glicci Cohen-Macaulay ideals over $R$. 
\end{ex}

But it is very interesting to see that the result holds in height three and higher without any additional assumptions. 

\begin{prop}
    Let $R$ be a local Gorenstein ring, and let $I$ be a Cohen-Macaulay ideal of height at least two and $\mu(I)=\height I+1$. Then $I$ is glicci.  
\end{prop}
\begin{proof}
    Choose a generating set for $I$ such that $I=(f_1,...,f_n)$ where $n=\mu(I)=\gr I+1$, and $C=(f_1,..,f_{n-1})\seq I$ is a complete intersection. Now set $J=C:I$, then $J\sim I$. Now $\omega_{R/J}\cong I/C$; as $I/C$ is principally generated, we get $J$ to be a Gorenstein ideal. As Gorenstein ideals of height at least two are glicci (Corollary~\ref{Gorensteinidealsareglicci}) we get the required. 
\end{proof}

\begin{thm}\label{reductionofdeviation}
    Let $S$ be a local Gorenstein ring with an infinite residue field. Let $I$ be a proper Cohen-Macaulay ideal with $\height g\ge 3$ and $\type S/I\ge 2$. If $\mu(I)\ge g+3$ then there exists a Cohen-Macaulay ideal $J$ with $\mu(J)\le g+2$ that is in the Gorenstein linkage class of $I$, linked by a sequence
of perfect Gorenstein ideals.

    In particular, $I\sim\cdots\sim I'\sim_GJ'\sim J$, where $\sim$ denote complete intersection links and $\sim_G$ denotes Gorenstein links; i.e., it requires at most one Gorenstein link that is not a complete intersection.   
\end{thm}

We need to note the following lemmas, which are an important tool in the proof. 

\begin{lem}\label{tightbourbaki}
    Let $R$ be a local Gorenstein ring with an infinite residue field. Let $M$ be an orientable MCM $R$-module that is not free and $\rank M\ge 2$. Then there exists an N-type resolution $($\emph{tight N-type resolution}$)$ $0\x R^{\rank M-1}\x M\x I\x 0$ such that $\mu(M)=\mu(I)+\rank M-1$. 
\end{lem}
The proof of this follows from the proof of \cite[Lemma~2.1]{herzog1987maximal}. In fact, they have explicitly mentioned this as a claim and have proven it. We need the following lemma on matrix factorizations by Eisenbud and Peeva.

\begin{lem}\cite[Theorem~2.1.1]{eisenbud_peeva_2016}\label{matrixfactorizations}
Let $(S,\mathfrak m)$ be Gorenstein local, let $f\in\mathfrak m$ be a nonzerodivisor, and set $R=S/(f)$. A matrix factorization of $f$ is a pair of matrices $A,B\in\operatorname{Mat}_q(S)$ satisfying $AB=BA=fI_q$; it is minimal if both matrices have entries in $\mathfrak m$. For a minimal factorization, $L=\operatorname{coker}A$ is maximal Cohen-Macaulay over $R$, with $\operatorname{pd}_S L=1$, $L^*=\operatorname{Hom}_R(L,R)\cong\operatorname{coker}(A^{\mathsf t})$, and $\mu_R(L)=\mu_R(L^*)=q$.
\end{lem}

We are now ready to prove the theorem. Although the constructions may initially seem unintuitive, their motivation becomes clearer by reading the proof backward, starting with the desired exact sequence in Step 3 and the choice of an orientable rank-two MCM module. This naturally leads to the constructions involving canonical modules. We retain the present order, however, as it makes the details easier to follow.

\begin{proof}[Proof of Theorem~\ref{reductionofdeviation}] Again to keep the proof readable, we split it into parts. It suffices to prove the claim in height three; we leave the straightforward reduction to the reader.


\noindent \textbf{1. }Choosing a regular element to cut down height.

Let $q=\type(S/I)\ge2$ 
and choose an $S$-regular sequence $x,y,z\in I$. Let $D$ be the $q\times q$ tridiagonal matrix with diagonal entries $x$, superdiagonal entries $y$, and subdiagonal entries $z$, and set $f=\det D$. Since $f\equiv x^q\pmod{(y,z)}$, the sequence $f,y,z$ is $S$-regular. So, $R=S/(f)$ is Gorenstein, and $f\in I^q\subseteq\mathfrak mI$, so $\mu_R(I/(f))=\mu_S(I)$.\\

\noindent\textbf{2: }Choosing a rank 2 orientable MCM $R$-module.

The identity $D\operatorname{adj}(D)=\operatorname{adj}(D)D=fI_q$ gives a minimal matrix factorization. By the preceding lemma, $L=\coker_S D$ and $L^*=\h_R(L,R)\cong\coker_S D^{\mathsf t}$ are MCM $R$-modules, with $\mu_R(L)=\mu_R(L^*)=q$ and $\operatorname{pd}_S L=\operatorname{pd}_S L^*=1$.

We next show that $L$ is locally free of rank one in codimension one. Fix $\mathfrak p\in\operatorname{Spec}R$ with $\operatorname{ht}\mathfrak p\le1$. Since $f,y,z$ is an $S$-regular sequence, $\overline y,\overline z$ is an $R$-regular sequence. Hence $\mathfrak p$ cannot contain both $\overline y$ and $\overline z$, so at least one of them is a unit in $R_{\mathfrak p}$. Now, $y^{q-1},z^{q-1}\in I_{q-1}(D)$. So $(\overline D)_{\mathfrak p}$ has an invertible $(q-1)\times(q-1)$ submatrix over $R_{\mathfrak p}$. Using the invertible \((q-1)\)-minor as a pivot block, we may reduce
\((\overline D)_{\mathfrak p}\) by row and column operations to
\(\operatorname{diag}(I_{q-1},a)\). Since \(\det D\) vanishes in \(R\)
and these operations preserve the determinant up to a unit, we have
\(a=0\). Taking cokernels gives \(L_{\mathfrak p}\cong R_{\mathfrak p}\).
Taking cokernels therefore gives $L_{\mathfrak p}\cong R_{\mathfrak p}$. Hence, both $L$ and $L^*$ are locally free of rank one in depth one.

  Set $N=L\oplus L^*$ then by Remark~\ref{N+N*isorientable}, we get that $N$ is orientable MCM $R$-module of rank 2. Also note that $N$ has no free summands.\\ 

\noindent\textbf{3.}  Lifting the canonical-module generators to $M$.

Let $m=\mu_S(I)$, $Q=\Omega_R^2(S/I)$, and $M=Q^*$. Since $S/I$ is Cohen-Macaulay, and height $\height_R \overline{I}=2$, and taking $R$-dual to the exact sequence $0\x \Omega^2_R(S/I)\x R^m\x I\x 0$ we get the following exact sequence 
\[
0\x R\x R^m\x M
\xrightarrow{\theta}\omega_{S/I}\x0.
\]

From the definition of $M$ it is clear that $M$ is an orientable MCM $R$-module.   Since $\mu_R(\omega_{S/I})=q$, the same sequence gives $\mu_R(M)\le m+q$. Choose a minimal surjection $\pi:S^q\to\omega_{S/I}$, as  $I\omega_{S/I}=0$ we get that $\pi D=0$. Thus $\pi$ factors through $L=\coker_S D$, giving an $R$-linear surjection $\rho:L\to\omega_{S/I}$.

Let $V=\ker\theta$, then $\operatorname{pd}_R(V)\le 1$, hence $\Ext^1_R(L,V)=0$. So, there exists an $\alpha:L\to M$ satisfying $\theta\alpha=\rho$. Now, the surjection $\rho\otimes_R k$ is an isomorphism because $\mu_R(L)=\mu_R(\omega_{S/I})=q$. Hence, the map $\alpha\otimes_R k$ is injective. Thus the images under $\alpha$ of a minimal generating set of $L$ give $q$ minimal generators of $M$.\\

\noindent\textbf{4.} Constructing the exact sequence and counting generators.

Choose minimal surjections $\beta:R^q\to L^*$ and $\gamma:R^q\to L$. Dualizing the map $\alpha:L\to M$ we get  $\alpha^*:Q=M^*\to L^*$. Since $\beta$ is surjective, the map $(\alpha^*\oplus\beta):Q\oplus R^q\to L^*$ is surjective. Set $K_1=\ker(\alpha^*\oplus\beta)$ and $K_2=\ker\gamma$, so that we get an exact sequence of MCM $R$-modules
\[
0\longrightarrow K_1\longrightarrow Q\oplus R^q
\xrightarrow{(\alpha^*\oplus\beta)}L^*\longrightarrow0,
\qquad
0\longrightarrow K_2\longrightarrow R^q
\xrightarrow{\gamma}L\longrightarrow0.
\]

Taking their direct sum and dual we get an exact sequence of orientable MCM $R$-modules
\[
0\x N=L\oplus L^*
\xrightarrow{\Phi}M\oplus R^q\oplus R^q
\x E=K_1^*\oplus K_2^*\x0,
\qquad
\Phi(u,v)=(\alpha(u),\beta^*(u),\gamma^*(v)).
\]
Now, $\rank_R E
=(m-1)+2q-2=m+2q-3$. Also note that $\beta^*(L)\subseteq\overline{\mathfrak m}R^q$ and $\gamma^*(L^*)\subseteq\overline{\mathfrak m}R^q$. So, $\Phi\otimes k=(\alpha\otimes_R k,0,0)$. By the previous step its rank is q. Hence $
\mu_R(E)=\mu_R(M)+2q-q
=\mu_R(M)+q\leq m+2q$. Also, note that 
$E$ is not free; if not then $L$ is free of rank one which is contradiction to $\rank_R L^*=1$ and $\mu_R(L^*)=q\ge 2$.\\

\noindent\textbf{5.} Gorenstein linkage. 

Consider a tight N-type resolution of $E$, $0\x R^{m+2q-4}\x E\x \overline{J}\x 0$. By the tightness we get $\mu_R(\overline{J})\le 4$, then $\mu_S(J)\le 5$. Consider an N-type resolution $0\x F'\x M=\Omega^2_R(S/I)^*\x \overline{I'}\x 0$, and an N-type resolution $0\x F\x T\x \overline{I}\x 0$. Then by Corollary~\ref{oddCIlinkageclass}, we get $\overline{I}$ and $\overline{I'}$ are in the same odd complete intersection linkage class. By Theorem~\ref{lemma3.4CDH}, we get $\overline{I'}$ and $\overline{J}$ are in the same even Gorenstein linkage class. So $I$ is in the Gorenstein linkage class of $J$.

For the perfectness of linking ideals; we have $I$ and $I'$ are linked by complete intersections. By Remark~\ref{linksin3.4CDH} we get $\overline{I'}\sim\cdots \sim \overline{I''}\sim_G\overline{J'}\sim J$. Again from the same remark we get this linking Gorenstein ideal has an N-type resolution $0\x R\x N\x \overline{K}\x 0$. Now as $\operatorname{pd}_S N=1$, we get that $K$ is a perfect Gorenstein ideal of height 3.  

\end{proof}

\begin{cor}
    Let $R$ be a Gorenstein local ring, and let $I$ be a Cohen-Macaulay ideal with $\height I\ge 3$. Suppose that if $I$ is not glicci, then $\inf \mathfrak{S}_R(I)=\height I+2$.
\end{cor}

It is also interesting to see that this theorem does not hold when we replace Gorenstein linkage by complete intersection linkage. So this leads us to ask the following question.

\begin{question}\label{minimaldeviation} Let $R$ be a local Gorenstein ring, and let $I\subseteq R$ be a proper Cohen-Macaulay ideal
of grade $g$. Set $\Sigma(I)=\{\mu(J)\mid J \text{ belongs to the
complete intersection linkage class of } I\}$. Clearly, $g\leq\inf\Sigma(I)\leq\mu(I)$. If $I$ is not licci. Under what conditions on $I$ is $\inf\Sigma(I)=\mu(I)$?   
\end{question}

\begin{remark}
    If height of $I$ is one, then by Lemma~\ref{CIlinkageclassesofheight1}, we get $\inf\Sigma(I)=\min\{\mu(I),\mu(I^*)\}$. So, this answers the previous question for the height one case.
\end{remark}

\section{Further Questions}
Throughout this paper, several questions arose naturally from the results and methods developed in the preceding sections. For the convenience of the reader, we collect all these questions here in one place. We retain their original numbering and provide references to the relevant sections, where the motivation, known results, and further discussion can be found.  \\

\begin{recalledquestion}{glicciconjecture}
      Let $R$ be a regular (graded local) local ring and let $I$ be a perfect (homogeneous) ideal of positive height. Then $I$ is in the (homogeneous) Gorenstein linkage class of a complete intersection.
\end{recalledquestion}

\begin{recalledquestion}{extendedglicciquestion}
Let $R$ be a Gorenstein local (graded local) ring and let $I$ be a perfect (homogeneous) ideal of positive height. Then $I$ is in the (homogeneous) Gorenstein linkage class of a complete intersection.    
\end{recalledquestion}

\begin{recalledquestion}{whenisCMglicci}Under what additional assumptions on a local Gorenstein
ring $R$ is every Cohen-Macaulay ideal of $R$ glicci?
\end{recalledquestion}

\begin{recalledquestion}{questiongolod}  Is every proper unmixed ideal, over a Gorenstein local ring in the complete intersection linkage class of some ideal $J$, such that the natural projection map $R\x R/J$ is a Golod homomorphism.
    
\end{recalledquestion}

\begin{recalledquestion}{deformationquestion}
     Let $(S,J)$ be a deformation of $(R,I)$; that is, there exists
a sequence $\alpha_1,\ldots,\alpha_n$ which is regular on both $S$ and
$S/J$, together with an isomorphism
\[
R \cong S/(\alpha_1,\ldots,\alpha_n)\qquad\text{under which}\qquad
I \cong
\frac{J+(\alpha_1,\ldots,\alpha_n)}
     {(\alpha_1,\ldots,\alpha_n)}.
\]
Assume that $S$ and $R$ are regular local rings and that $I$ is
$\mathfrak m_R$-primary. If $I$ is glicci, then is $J$ glicci?
\end{recalledquestion}

\begin{recalledquestion}{questiononeheightup}
     Let $R$ be a graded-local Gorenstein ring. Let $I$ be a homogeneous Cohen-Macaulay ideal of $R$. If $f\in\mathfrak m$ is a homogeneous non zero divisor of $R/I$, then is $(I,f)$ homogeneously glicci in $R$?
\end{recalledquestion}

\begin{recalledquestion}{DLextensionvariablequestion}
     Let $B$ be a two dimensional local Gorenstein ring, with an infinite residue field. Let $R=B[[t_1,..,t_n]]$ or $R=B[t_1,..,t_n]$. Is every orientable MCM $R$-module stably isomorphic to a doubly layered $R$-module?
\end{recalledquestion}

\begin{recalledquestion}{Questionmonomial2glicci}
Let $S$ be a polynomial ring over an infinite field,
and let $I$ be a Cohen--Macaulay monomial ideal.
Is $I$ homogeneously $2$-glicci?
\end{recalledquestion}

\begin{recalledquestion}{minimalassquestion}Let $S$ be a polynomial ring over a field.
Does there exist a Cohen-Macaulay monomial ideal
of $S$ with exactly eight or nine associated
primes that is not clean?    
\end{recalledquestion}

\begin{recalledquestion}{minimaldeviation} Let $R$ be a local Gorenstein ring, and let $I\subseteq R$ be a proper Cohen-Macaulay ideal
of grade $g$. Set $\Sigma(I)=\{\mu(J)\mid J \text{ belongs to the
complete intersection linkage class of } I\}$. Clearly, $g\leq\inf\Sigma(I)\leq\mu(I)$. If $I$ is not licci. Under what conditions on $I$ is $\inf\Sigma(I)=\mu(I)$? 
\end{recalledquestion}

\bibliographystyle{alpha}

\bibliography{reference_cleaned}

\end{document}